\documentclass[12pt,reqno]{amsart}

\usepackage{amsmath, amsthm}

\usepackage{geometry}
\usepackage{fullpage}
\usepackage[all]{xy}

\usepackage[mathcal,mathscr]{euscript}
\usepackage[OT2,T1]{fontenc}
\usepackage[english]{babel}
\usepackage{newtxtext}
\usepackage{newtxmath}

\usepackage{etoolbox}
\patchcmd{\section}{\scshape}{\bfseries}{}{}
\makeatletter
\renewcommand{\@secnumfont}{\bfseries}
\makeatother
\numberwithin{equation}{section}
\usepackage{enumitem}

\newtheorem{introtheorem}{Theorem}

\newtheorem{introprop}[introtheorem]{Proposition}
\theoremstyle{definition}
\newtheorem*{introdef*}{Definition}
\newtheorem*{introremark*}{Remark}
\theoremstyle{plain}
\newtheorem{theorem}{Theorem}[subsection]
\newtheorem{proposition}[theorem]{Proposition}
\newtheorem{lemma}[theorem]{Lemma}
\newtheorem{corollary}[theorem]{Corollary}

\newcommand{\thistheoremname}{}
\newtheorem*{genericthm}{\thistheoremname}
\newenvironment{namedtheorem*}[1]
  {\renewcommand{\thistheoremname}{#1}%
   \begin{genericthm}}
  {\end{genericthm}}
\theoremstyle{definition}
\newtheorem{definition}[theorem]{Definition}

\newtheorem{remark}[theorem]{Remark}

\newtheorem{example}[theorem]{Example}

\NewCommandCopy{\proofqedsymbol}{\qedsymbol}

\AtBeginEnvironment{proof}{\renewcommand{\qedsymbol}{\proofqedsymbol}}

\AtBeginEnvironment{introexample}{%
  \pushQED{\qed}\renewcommand{\qedsymbol}{$\lozenge$}%
}
\AtEndEnvironment{introexample}{\popQED}

\AtBeginEnvironment{example}{%
  \pushQED{\qed}\renewcommand{\qedsymbol}{$\lozenge$}%
}
\AtEndEnvironment{example}{\popQED}

\AtBeginEnvironment{examples}{%
  \pushQED{\qed}\renewcommand{\qedsymbol}{$\lozenge$}%
}
\AtEndEnvironment{examples}{\popQED}

\AtBeginEnvironment{introremark}{%
  \pushQED{\qed}\renewcommand{\qedsymbol}{$\lozenge$}%
}
\AtEndEnvironment{introremark}{\popQED}

\AtBeginEnvironment{remark}{%
  \pushQED{\qed}\renewcommand{\qedsymbol}{$\lozenge$}%
}
\AtEndEnvironment{remark}{\popQED}

\AtBeginEnvironment{remarks}{%
  \pushQED{\qed}\renewcommand{\qedsymbol}{$\lozenge$}%
}
\AtEndEnvironment{remarks}{\popQED}

\AtBeginEnvironment{definition}{%
  \pushQED{\qed}\renewcommand{\qedsymbol}{$\triangle$}%
}
\AtEndEnvironment{definition}{\popQED}

\usepackage[usenames, dvipsnames, table,svgnames]{xcolor}
\usepackage[hidelinks,hypertexnames=false]{hyperref}
\usepackage{cleveref}
\usepackage{mathtools}
\usepackage{mathrsfs}

\usepackage{tabularray}

\newcommand{\Z}{\mathbf{Z}}
\newcommand{\R}{\mathbf{R}}
\newcommand{\Q}{\mathbf{Q}}

\newcommand{\Cc}{\mathbf{C}}
\newcommand{\F}{\mathcal{F}}
\newcommand{\Hh}{\mathbb{H}}
\newcommand{\PP}{\mathbf{P}}
\newcommand{\E}{\mathscr{E}}
\newcommand{\e}{\mathrm{e}}

\newcommand{\A}{\mathscr{A}}
\newcommand{\M}{\mathrm{M}}

\renewcommand{\geq}{\geqslant}
\renewcommand{\leq}{\leqslant}
\newcommand{\ii}{\mathrm{i}}
\renewcommand{\Re}{\mathrm{Re}}
\renewcommand{\Im}{\mathrm{Im}}
\renewcommand\d{\mathop{}\!\mathrm{d}} 
\DeclareMathOperator*{\LS}{\mathrm{Ls}}
\DeclareMathOperator*{\LI}{\mathrm{Li}}
\DeclareMathOperator*{\LH}{\mathrm{Lim}} 
\usepackage{tikz}
\usepackage{dsfont}
\usepackage{pgfplots}
\pgfplotsset{compat=1.18}
\usetikzlibrary{decorations.markings}
\usetikzlibrary{decorations.pathreplacing}
\usetikzlibrary{patterns}

\makeatletter
\newcommand*{\defeq}{\mathrel{\rlap{%
                     \raisebox{0.24ex}{$\m@th\cdot$}}%
                     \raisebox{-0.24ex}{$\m@th\cdot$}}%
                     =}
\makeatother

\newcommand{\ord}{\operatorname{ord}}

\usepackage{mathtools}

\DeclarePairedDelimiter\floor{\lfloor}{\rfloor}
\newcommand{\lt}{\ensuremath <}
\newcommand{\gt}{\ensuremath >}
\renewcommand{\le}{\leqslant}
\renewcommand{\ge}{\geqslant}
\newcommand{\ga}{\gamma}
\newcommand{\al}{\alpha}
\newcommand{\be}{\beta}
\newcommand{\de}{\delta}
\newcommand{\vph}{\varphi}
\newcommand{\f}{\frac}
\newcommand{\cs}{\mathscr}

\newcommand{\ol}{\overline}

\newcommand{\bb}{\mathbf}
\newcommand{\dx}{\d x}
\renewcommand{\mod}[1]{\ (\mathrm{mod}\ #1)}
\DeclareMathOperator{\sgn}{sgn}
\DeclareMathOperator{\me}{e}
\DeclareMathOperator{\Arg}{arg}

\definecolor{myblue}{HTML}{1F4C6F}
\definecolor{mypink}{HTML}{D83079}
\definecolor{mypurple}{HTML}{7C3E74}
\usetikzlibrary{plotmarks}
\usepackage{ninecolors}

\begin{document}
\date{\today}
\title[Zeros of Eisenstein series]{Geodesic clustering of zeros\\[1mm] of Eisenstein series for congruence groups}
\author[S.~Carrillo]{Sebasti\'an Carrillo Santana}
\author[G.~Cornelissen]{Gunther Cornelissen}
\author[B.~Ringeling]{Berend Ringeling}
\address{\parbox{\linewidth}{\normalfont Mathematisch Instituut, Universiteit Utrecht, Postbus 80.010, 3508 TA Utrecht, Nederland \\
Max-Planck-Institut f\"ur Mathematik, Postfach 7280, 53072 Bonn, Deutschland \\ \vspace*{-3mm} \mbox{ } }}
\email{g.cornelissen@uu.nl}
\address{\normalfont Mathematisch Instituut, Universiteit Utrecht, Postbus 80.010, 3508 TA Utrecht, Nederland}
\email{s.carrillosantana@uu.nl}
\address{\normalfont CRM / Universit\'e de Montr\'eal, P.O. Box 6128, Centre-ville Station
Montr\'eal (Qu\'ebec) H3C 3J7, Canada}
\email{b.j.ringeling@gmail.com}

 \subjclass[2010]{primary: 11F11 ; secondary: 11J91} 
 \keywords{\normalfont Modular form, Congruence group, Eisenstein series, Zeros, Transcendence, Hausdorff convergence, Equidistribution}
\thanks{We thank Walter Bridges, Jakub Byszewski, Emmanuel Kowalski, P\"ar Kurlberg, Johannes Linn, Zeev Rudnick and Don Zagier for useful discussions. G.C.\ thanks the Max-Planck-Institut f\"ur Mathematik in Bonn for excellent conditions during work on this paper.
}

\begin{abstract} \noindent 
We consider a set of generators for the space of Eisenstein series of varying even weight $k$ for any fixed congruence group $\Gamma$ and study the set of all of their zeros taken for all $\Gamma(1)$-conjugates of $\Gamma$ in the standard fundamental domain for $\Gamma(1)$. We describe (a) an upper bound $\kappa_\Gamma + O(1/k)$ for their imaginary part;  (b) a finite configuration of geodesics segments to which all zeros converge in Hausdorff distance as $k \rightarrow \infty$; (c) a finite set containing all algebraic zeros for all weights. 
The bound in (a) depends on the (non-)vanishing of a new generalization of Ramanujan sums.  The proof of (b) originates in a method used to study phase transitions in statistical physics. The proof of (c) relies on the theory of complex multiplication.
The results can be made quantitative for specific groups. For $\Gamma=\Gamma(N)$ with $4 \nmid N$, $\kappa_\Gamma=1$ and the zeros tend to the unit circle, whereas if $4 \mid N$, $\kappa_\Gamma=2$ and the limit configuration includes parts of vertical geodesics and circles of radius $2$. 
{In both cases, for sufficiently large $k$, the only algebraic zeros are at $\mathrm{i}$ and $\exp(2\pi \mathrm{i}/3)$, and, for $4 \mid N$, at $2\mathrm{i}$.}

For $\Gamma(N)$ with $N$ odd, we use finer estimates to prove a trichotomy for the exact `convergence speed' of the zeros to the unit circle, as well as angular equidistribution of the zeros as $k \rightarrow + \infty$, {implying weak convergence of the normalized zero counting measure on $Y(1)$ to the image of the argument measure on the lower arc of the fundamental domain}.
\end{abstract}

\maketitle

\vspace*{-1cm}

{\footnotesize
\setcounter{tocdepth}{1}
\renewcommand{\contentsname}{\centering \textnormal{\textbf{Contents}}}
\tableofcontents
}
\section{Introduction}
We will be concerned with properties of zeros of a set of modular forms of varying even weight spanning the space of Eisenstein series. Contrary to most of the existing literature, our goal is to find a general framework applicable to all congruence groups; a small survey of related results can be found {in the next subsection of this introduction.} 

\begin{figure}
\centering
\begin{minipage}{.5\textwidth}
\mbox{ } \vspace*{-1mm}
  \centering
\includegraphics[width=46mm]{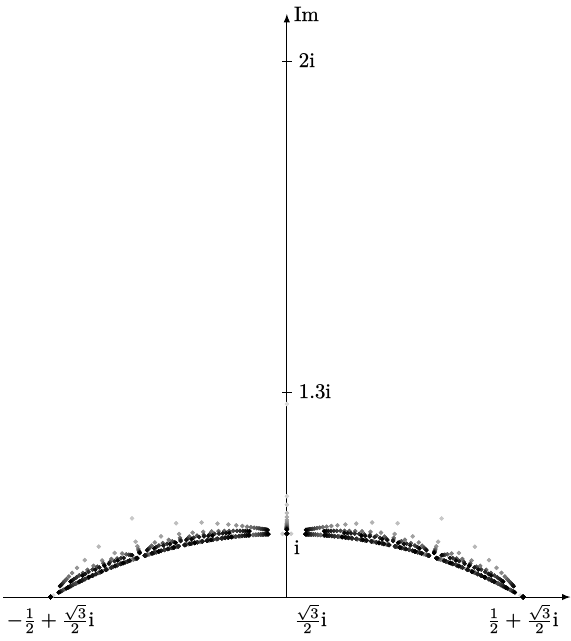} 
\end{minipage}%
\begin{minipage}{.5\textwidth}
  \centering
\includegraphics[width=45mm]{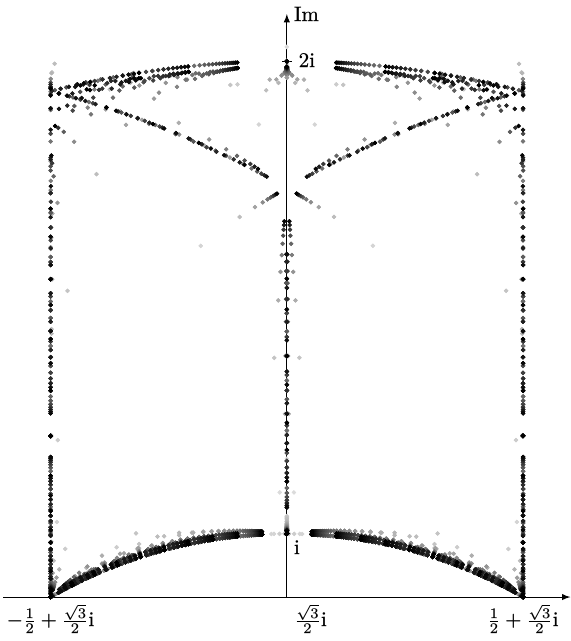}
\end{minipage}
\caption{Every dot represents one of the more than 1000 zeros of all Eisenstein series for $\Gamma(N)$ of even weight in $[4,K]$ in the closure of $\F$, where darker points correspond to higher weight. On the left, $N=3, K=80$; on the right, $N=8, K=60$.} \label{manyzeros}
\end{figure}

Our choice of spanning set and the detailed results will be described below, but we will start by displaying two examples. For our choice of generating Eisenstein series  for a principal congruence group $\Gamma(N) \leq \Gamma(1) =\mathrm{SL}_2(\Z)$ as in Formula \eqref{defeis} below, the union of all zeros is $\Gamma(1)$-invariant, so we can plot the zeros in the (closure of the) standard fundamental domain $\F$ for $\Gamma(1)$. On the left of Figure \ref{manyzeros}, we do this for $\Gamma(3)$, where the zeros have small imaginary part, and, as the weight increases, tend to the unit circle. The right of Figure \ref{manyzeros} concerns the case of the group $\Gamma(8)$, where {additional clustering patterns appear,}  
such as zeros of larger imaginary part clustering around parts of circles of radius $2$ and parts of the vertical boundary, as well as of the imaginary axis below $\ii \sqrt{3}$ and around $2\ii$---in particular, the zeros are all close to segments of geodesics, {a phenomenon we dub `geodesic clustering'}. 

\subsection*{Related results} { \mbox{ }} Zeros of families of modular forms have been studied in various contexts. We cannot possibly survey the over 100 references we found, so we will make a selection. One of the historically first references is the last joint paper of Hardy and Ramanujan \cite{HardyRamanujan}, in which they study asymptotics of inverses of modular forms (compare \cite{HeimNeuhauser}). For Eisenstein series, the problem of locating the zeros was picked up by Wohlfahrt \cite{Wohlfahrt}, (Robert) Rankin \cite{Rankin}, and finally (Fenny) Rankin and Swinnerton-Dyer \cite{RSD}, who proved that all zeros  of Eisenstein series of level $N=1$ in  {the standard fundamental domain}   
are located \emph{precisely on} the lower {boundary} arc 
The result may be refined to show interlacing of the zeros  \cite{Nozaki} and equidistribution of the arguments. 
The result has been generalized, e.g.\ for other Poincar\'e series (e.g., \cite{RankinPoincare}), for other specific groups (e.g., \cite{Brazelton}), for quasimodular forms (e.g., \cite{DukeJenkins}, \cite{GunO}, \cite{BJW}) and for Drinfeld modular forms (e.g., \cite{GC}, \cite{Gekeler}). 

In a seemingly different direction, seminal work of Rudnick \cite{Rudnick} showed that zeros of cuspidal Hecke eigenforms of level $1$ equidistribute in $\Hh$ for the hyperbolic measure. This is proven by reduction to a mass equidistribution statement in the theory of automorphic forms; see also \cite{HolowinskySound} and \cite{Mato}. In this setup, Ghosh and Sarnak \cite{GS} have studied the distribution of zeros of such forms that lie exactly on the `vertical' geodesics $\Re(z)\in \{-1/2,0\}
$; their number grows at least like $k^{1/4-\epsilon}$ (see also \cite{Mato2}) but $\Im(z)$ is unbounded as $k \rightarrow + \infty$. 
Very recently, Rudnick \cite{RudnickFaber} has shown that zeros of extremal modular forms (i.e., that have the maximal number of trailing zero Fourier coefficients for their weight $k$) accumulate, up to an error of $O(1/k)$, near a finite set of vertical lines in the standard fundamental domain. This appears to be the first example of a geodesic clustering result for such zeros tending to, but not lying on, finitely many geodesics. Since then, further such results on the clustering of zeros close to  geodesics in the upper half plane have been obtained, see, e.g., \cite{Kimmel}, \cite{Raveh}.

\subsection*{Notation and Terminology.} We will use the shorthand notations $\e(z) \defeq \me^{2 \pi \ii z}$, $q \coloneq \e(z)$, $q_N\coloneq \e(z/N)$, $\zeta_N \defeq \e(1/N)$ and $\rho \defeq \zeta_3$. We will write $\Hh \defeq \{ z \in \Cc \colon \Im(z)>0 \}$ for the complex upper half plane. We denote by $\{x\}$ the fractional part of $x$. 
We define the standard fundamental domain $\F$ for $\Gamma(1)$ as 
\begin{align*} &\F \defeq \{ z \in \Cc \colon |z|>1 \mbox{ and }-1/2 \leq \Re(z) < 1/2 \} \cup \A \\ &\mbox{ with } \A \defeq \{ \me^{i \theta} \mbox{ with } \pi/2 \leq \theta \leq 2 \pi/3 \}, \end{align*}
where we will call $\A$ the `(left) lower arc'. The `right lower arc', denoted $\A'$, is the negative complex conjugate of $\A$.   

Sets of zeros of functions (in $\F$ or its closure) will always be considered as multisets, with their cardinality keeping track of the multiplicities, and adhering to the convention that the multiplicity of $\ii$ is divided by $2$ and the multiplicity of $\rho$ is divided by $3$, so that, for example, the set of zeros of a non-cusp form of weight $k$ for $\Gamma(1)$ in $\F$ has cardinality $k/12$. 

If $\Gamma$ is a congruence group and $N$ is the smallest integer such that $\Gamma(N) \leq \Gamma$, we call $N$ the \emph{level} of $\Gamma$ and write $h_\Gamma\defeq [\Gamma:\Gamma(N)]$. Also, let $N_\Gamma = \Gamma(1)/\mathrm{Nor}_{\Gamma(1)}(\Gamma)$, where $\mathrm{Nor}_{\Gamma(1)}(\Gamma) \leq \Gamma(1)$ is the normalizer of $\Gamma$ in $\Gamma(1)$; notice that $N_\Gamma$ is finite, since $\Gamma$ is of finite index in $\Gamma(1)$. 
Once a level $N$ has been fixed, we will sometimes use the shorthand notation $u_j\defeq \cos(2\pi j/N)$ for integers $j$.

\subsection*{Which Eisenstein series?} Since, in general, the vector space $\E_k(\Gamma)$ of Eisenstein series for a congruence group $\Gamma$ (i.e., the orthogonal complement of the space of cusp forms for the Petersson inner product) is more than one-dimensional, any clustering property of zeros of a subset of $\E_k(\Gamma)$ will be a property of the specific set, not the space $\E_k(\Gamma)$.  
For a principal congruence group $\Gamma(N)$ and weight $k \geq 3$ we will use the following set of Eisenstein series: 
\begin{equation} \label{defeis} E_k^{a,b}(z) = E_k^{(a,b)}(z) \defeq \sum_{(0,0) \neq (m,n) \in \Z^2} \frac{\e(an-bm)}{(mz+n)^k} \end{equation} 
with $(a,b)$ running over all values in $(N^{-1}\Z/\Z)^2.$ The series \eqref{defeis} are specialisations of Eisenstein--Kronecker double series \cite[\S 3.1]{BrunaultZudilin}, and were used, for example, by Kato in in the construction of modular Euler systems \cite[\S 3]{Kato}, by Brunault and Zudilin in their work on modular regulators  \cite{BrunaultZudilin} and by Brunault in connection with the theory of toric modular forms \cite{Brunault}. 
The set $\{E_k^{a,b}\}$ generates the space $\E_k(\Gamma)$. In a preliminary section, we explain the relation to other generating sets given by Hecke \cite{Hecke} (through $(\Z/N\Z)^2$-Fourier transform; see \eqref{fexp}), we give some standard transformation formulae and relations (\eqref{transfo}, Lemma \ref{eisrel}) and {we prove that the $\{E_k^{a,b}\}$ with $(a,b)$ running over primitive pairs up to sign form a basis, see Proposition \ref{basis}.}   
 If $\Gamma$ is a general congruence group of level $N$, we consider the set of `traces'  
\begin{equation} \label{traces} E^{a,b}_{k,\Gamma} = E^{(a,b)}_{k,\Gamma} \defeq \sum_{g \in \Gamma(N) \backslash \Gamma} E_k^{(a,b)g} \end{equation} with $(a,b)$ running over $(N^{-1}\Z/\Z)^2$; these form a generating set for the space of Eisenstein series for $\Gamma$; see Proposition \ref{Tab}. We prove that none of these modular forms is zero at any cusp.  

 Turning to the zeros, given a congruence group $\Gamma$, an even weight $k \geq 4$ and an index $(a,b) \in (N^{-1}\Z/\Z)^2$, we let 
$Z_{k,\Gamma}^{a,b}$ denote the multiset of zeros  of $E_{k,\Gamma}^{a,b}$ in $\F$. This is not necessarily $\Gamma(1)$-invariant, but by taking into account all conjugate groups, we define the $\Gamma(1)$-invariant multiset of all zeros in even weight $k \geq 4$ as
$$ Z_{k,\Gamma} \defeq \bigcup_{\substack{(a,b) \in (N^{-1}\Z/\Z)^2 \\ \gamma \in N_\Gamma}} \hspace*{-5mm} Z^{a,b}_{k,\gamma^{-1}\Gamma\gamma}. $$
The union of $Z_{k,\Gamma}$ over all weights $k \in 2\Z_{\geq 2}$ is denoted by $Z_\Gamma$. 
We use corresponding `overlined' notation $\overline Z_{k,\Gamma}^{a,b}, \overline Z_{k,\Gamma}, 
\overline Z_\Gamma$ for the sets of zeros in the closure $\overline \F$ of $\F$; these are not the closures of the corresponding `unoverlined' sets.

\subsection*{Downward imaginary concentration} Our first basic result gives a zero-free region. To formulate it,  we use the following generalization of Ramanujan sums. 

\begin{introdef*} 
For $\Gamma$ a congruence group of level $N$ and $(a,b) \in (N^{-1}\Z/\Z)^2$, we define the \emph{Kluyver sum}\footnote{We name these sums after Jan C.~Kluyver, who introduced `Ramanujan sums' before Ramanujan in \cite[p.~410]{Kluyver}, including a proof of the two basic identities \eqref{kluyver}, later also found by and attributed to H\"older \cite[(13)]{Holder}, (von) Sterneck, and Ramanujan \cite{RamanujanTrig}.} as 
 $$ \rho_\Gamma^{a,b}(j) \defeq 
\ 2 \hspace*{-7mm} \sum_{\left( \begin{smallmatrix} \alpha & \beta \\ \gamma & \delta \end{smallmatrix} \right) \in \Gamma(N) \backslash \Gamma} \hspace*{-7mm} \cos(2 \pi(\alpha a + \gamma b)j). $$ 
 This generalizes the usual Ramanujan sums, that occur (up to a factor $2N$) for the group $\Gamma_0(N)$, as well as Kloosterman sums, that occur for conjugates of $\Gamma_0(N)$, see Example \ref{gamma0RS}. The most natural and general definition of a Kluyver sum is associated to a mod-$N$ representation of a finite group, cf.\ Remark \ref{anykluyver}. 
Furthermore, define constants 
\begin{equation*}  \kappa^{a,b}_\Gamma  \defeq \min \{ j \geq 1 \colon \rho_\Gamma^{a,b}(j) \neq 0\} \mbox{ and }
 \kappa_\Gamma \defeq \max \{ \kappa^{a,b}_{\gamma^{-1}\Gamma\gamma} \colon (a,b) \in (N^{-1}\Z/\Z)^2, \gamma \in N_\Gamma \}.
\end{equation*}
\end{introdef*}

\begin{introtheorem} \label{mainIm} For $\Gamma$ a congruence group of level $N$, $\kappa_\Gamma$ is finite; in fact, $\kappa_\Gamma \leq 2 h_\Gamma$ and $\kappa_\Gamma$ divides $N$. For an even weight $k \geq 4$, we have an upper bound 
 \begin{equation}
 \label{imup} \Im(z) < \kappa_\Gamma +O(1/k) \mbox{ for all } z \in Z_{k,\Gamma}, \end{equation} where the implied constant is independent of the weight $k$. 
The bound in \eqref{imup} is optimal, in the sense that there exist zeros $z_i \in Z_{\Gamma}$ ($i=1,\dots$) such that $\lim \Im(z_i) = \kappa_\Gamma$.
Furthermore,
\begin{enumerate} 
\item $\kappa_{\Gamma(N)}= 1$ if $N$ is not divisible by $4$, and $\kappa_{\Gamma(N)}=2$ if $N$ is divisible by $4$; 
\item for Hecke's theta group $\Lambda$, $\kappa_\Lambda=2$; 
\item $\kappa_{\Gamma_1(N)}=\kappa_{\Gamma_0(N)} = N$.
\end{enumerate} 
\end{introtheorem}

A zero-free region was given in level $1$ by R.~Rankin \cite{Rankin}, using the Fourier expansion of Eisenstein series (a result later superseded by \cite{RSD}). Our proof uses an arbitrary finite truncation of the Eisenstein series, and is more similar to the argument in \cite{RSD}, but with more control over the remainder term. The optimality proof will follow after Theorem \ref{mainconfig} below. 

\subsection*{Transcendence} Since the Eisenstein series we consider have Fourier coefficients in a fixed cyclotomic field, the $j$-invariants of their zeros in $\Hh$ are algebraic, and hence a classical theorem of T.~Schneider implies that their zeros in $\Hh$ are at CM points. The next result restricts these zeros in $\F$ even further. 

\begin{introtheorem} \label{tran} For a congruence group $\Gamma$, there are only finitely many algebraic zeros $z \in Z_\Gamma \cap \overline \Q$, and, for sufficiently large weight (depending on $\Gamma$), all of these belong to an imaginary quadratic field of discriminant $D$ with $|D| \leq 4 \kappa_\Gamma^2$. Furthermore, for  sufficiently large weight, 
\begin{enumerate}
\item if $\kappa_\Gamma=1$, then $z \in \{\ii, \rho \}$; 
\item if $\kappa_\Gamma=2$, then $z \in \{ \ii, \rho, \sqrt{-2}, \sqrt{-3}, \frac{-1+\sqrt{-7}}{2}, \frac{-1+\sqrt{-11}}{2},\frac{-1+\sqrt{-15}}{2}, \frac{-1+\sqrt{-15}}{4} {,2\ii} \}$. 
\end{enumerate}
\end{introtheorem} 

The method of proof is to use the theory of complex multiplication (akin to an argument of Kohnen) to show that if an Eisenstein series has such an algebraic zero of discriminant $D$, there must be a (possibly different) Eisenstein series {(on a possibly different  
congruence group)} with an algebraic zero of imaginary part $\sqrt{|D|}/2$, so we can invoke Theorem \ref{mainIm}. 

\subsection*{Explicit bounds} For a given congruence group, it is possible to find a suitable implied constant in Theorem \ref{mainIm} and use this  to quantify the meaning of `sufficiently large weight' in Theorem \ref{tran}, for example as follows. 

\begin{introprop} \label{explicit} For $N \geq 1$ not divisible by $4$ and $k\ge 2(\log(13N))^2$, all $z \in Z_{k, \Gamma(N)}$ satisfy $\Im(z) \leq 1 +\log(13N)/k$ and the only algebraic such $z$ are $\ii$ and $\rho$. 
\end{introprop}

We can use irrationality measures of $\mathrm{arctan}(\sqrt{n})/\pi$ for integer $n$ to eliminate many of the extra possible algebraic zeros for $\Gamma(N)$ even for $N$ divisible by $4$ (but now without an explicit lower bound on the weight).  

\begin{introprop} \label{evenN} For all $N \geq 1$ and $k$ sufficiently large, $Z_{k, \Gamma(N)} \cap \overline \Q {\subseteq} \{\ii, \rho {, 2\ii}\}$. 
\end{introprop}

It is also possible, using class number estimates, to upper bound the total number (with multiplicity) of algebraic zeros of any Eisenstein series of sufficiently large even weight for any congruence group of level $N$ polynomially in $N$, cf. Corollary \ref{absupper}.

\subsection*{Geodesic clustering of zeros} 

Theorem \ref{mainIm} says nothing about zeros below the line $\Im (z) = \kappa_\Gamma$. We now study the limiting behaviour of the zeros as the weight tends to infinity. The results are cleaner to formulate by using  the set of zeros in the closure of the fundamental domain $\F$. Recall that (hyperbolic) geodesics in the upper half plane are semicircles with center on the real axis, or vertical half lines given by complex numbers with fixed real part. We will call any non-empty connected subset of such a geodesic a \emph{geodesic segment}. 

\begin{figure}
\centering
\begin{minipage}{.4\textwidth}
\mbox{ } \vspace*{-1mm}
  \centering
  \includegraphics[width=.75\linewidth]{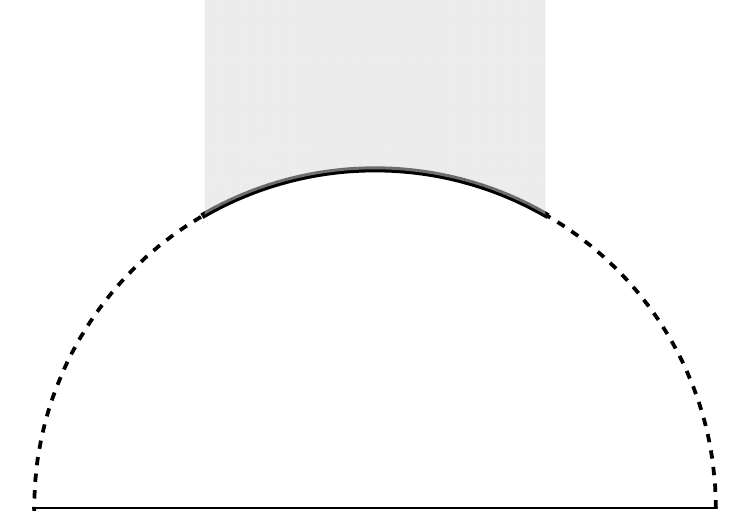}
\end{minipage}%
\begin{minipage}{.5\textwidth}
  \centering
  \includegraphics[width=1.05\linewidth]{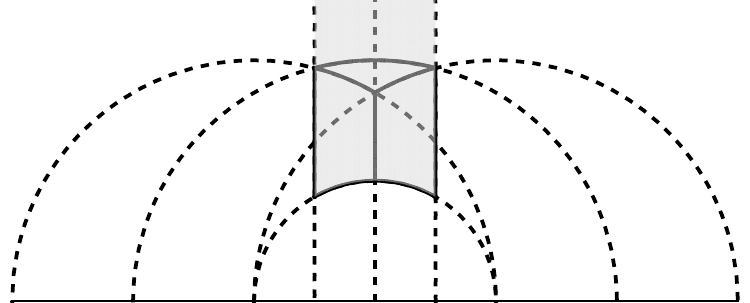}
\end{minipage}
\caption{Geodesic completion of  $\overline Z_{\infty,\Gamma(N)}$ (dashed lines), fundamental domain $\F$ (shaded region), and geodesic segments forming the limit $\overline Z_{\infty,\Gamma(N)}$ of the zeros in the closure of $\F$ as the weight tends to infinity (solid lines); left: $N$ not divisible by $4$, right: $N$ divisible by $4$. These theoretical pictures should be compared with the data from Figure \ref{manyzeros}.} \label{anharmpic}
\end{figure}

\begin{introtheorem} \label{mainconfig}
For any congruence group $\Gamma$, as the weight $k \in 2 \Z_{\geq 2}$ increases, the set $\overline Z_{k,\Gamma}$ converges in Hausdorff distance to a finite union $\overline Z_{\infty,\Gamma}$ of geodesic segments. More precisely, 
{$$ \overline Z_{\infty,\Gamma} = \hspace*{-9mm} \bigcup_{\substack{(a,b) \in (N^{-1}\Z/\Z)^2 \\ \gamma \in N_\Gamma}}   \bigcup_{\substack{(m_1,n_1) \neq (m_2,n_2) \\ \mathrm{in}\, I^{a,b}_{\gamma}}} \hspace*{-5mm} \{ z \in \overline \F \colon |m_1z+n_1|=|m_2z+n_2| \leq |m_3 z + n_3|, \ \forall (m_3,n_3) \in J^{a,b}_{\gamma} \}$$}
where  the finite index sets are defined by declaring $\mathbb I \defeq (\Z_{> 0} \times \Z) \cup (\{ 0 \} \times \Z_{>0})$ and
\begin{align*} 
I^{a,b}_{\gamma} &\defeq \{ (m,n) \in \mathbb I \colon m^2+n^2 \leq 2 (\kappa^{a,b}_{\gamma^{-1}\Gamma\gamma})^2 \mbox{ and } \rho_{\gamma^{-1}\Gamma\gamma}^{(Na,Nb),(n,-m)}(1) \neq 0\}; \\
J^{a,b}_{\gamma} &\defeq\{ (m,n) \in \mathbb I \colon m^2+n^2 \leq 7 (\kappa^{a,b}_{\gamma^{-1}\Gamma\gamma})^4 \mbox{ and } \rho_{\gamma^{-1}\Gamma\gamma}^{(Na,Nb),(n,-m)}(1) \neq 0\}
\end{align*}
in terms of the more general Kluyver sums as in Definition \ref{defgenram}. 
\end{introtheorem}

For the proof, we adapt a method of Sokal used to study Lee--Yang zeros (that describe phase transitions in statistical physics): this describes the limit in terms of complex semi-algebraic sets (i.e., equalities and inequalities of absolute values of holomorphic functions). 
The method, based on the fact that zeros of a meromorphic function $f$ are detected by non-harmonicity of $\log |f|$, needs to be adapted to our case to permit infinite (instead of finite) sums of powers of reciprocals of integer linear polynomials (instead of arbitrary meromorphic functions) in a compact (instead of arbitrary) subset $ \overline \F_\Gamma \defeq \{ z \in \overline \F \colon \Im(z) \leq \kappa_\Gamma\}.$  
In principle, the method allows us to find the limit set exactly. 

\begin{introprop} \label{introcinfty} For $\Gamma(N)$ with $N$ not divisible by $4$, $\overline Z_{\infty, \Gamma(N)} = \A \cup \A'$,
seen on the left in Figure \ref{anharmpic}. 
For $\Gamma(N)$ with $N$ divisible by $4$, $\overline Z_{\infty, \Gamma(N)}$ consists of the boundary of $\overline \F$ below the circle of radius $2$ centered at the origin, the part inside $\overline \F$ of that circle, the part inside $\overline \F$ of the two circles of radius $2$ centered at $\pm 1$ for which $\Im(z)\geq \sqrt{3}$, and the totally imaginary numbers $\ii y$ with $1 \leq y \leq \sqrt{3}$, seen on the right in Figure \ref{anharmpic}.  
\end{introprop} 
The continuation of all occurring geodesic segments to $\Hh$ forms a configuration with interesting properties; see, e.g., Remark \ref{curproj}. 
An example pertaining to the group $\Gamma^0(N)$, is found in Proposition \ref{Gamma0opt}, showing that certain segments occur only for powerful divisors of the level.  

\subsection*{Convergence speed and equidistribution} 

For a principal congruence group with odd level for all but a finite number (in the level) of elements of $Z_k\defeq Z_{k,\Gamma(N)}$, we know their location to even greater accuracy and can determine the exact `convergence speed' to the lower arc $\A \cup \A'$. To formulate the result, we use the Euclidean metric $d_E$ and our convention on how to count cardinalities of multisets of zeros.

\begin{introtheorem} \label{main2}
Let $N$ be an odd level; then for all even $k \geq 4$, we have $\#Z_k =kN^2/12.$  For sufficiently large even weight $k$, we have: 
\begin{enumerate}
\item \textup{(Many zeros approach the arc at moderate speed)} There exists a subset $\widehat{Z}_k^{\textup{off}} \subseteq Z_k$ with $$\#\widehat{Z}_k^{\textup{off}} \ge \frac{(N-1)^2k}{12}  - c_{\textup{off}}\, N^2 \log N$$ for some absolute constant $c_{\textup{off}}\gt 0$ such that  $\widehat{Z}_k^{\textup{off}}$ consists of only transcendental zeros that are at Euclidean distance $d_E(\widehat{Z}_k^{\textup{off}}, \A \cup \A') \asymp 1/k$ from the lower arc $\A \cup \A'$;
\item \textup{(Zeros on the arc)} On the arc $\A$, there are at least $$\#(Z_k \cap \A) \ge  \frac{k(2N-1)}{12} - c_{\textup{on}}\, N \log N$$ zeros for some absolute constant $c_{\textup{on}}\gt 0$; 
\item \textup{(Few exceptional zeros approach $\ii$ `arc-wise' or $\rho$ `vertically' at rapid speed)}. 
\begin{enumerate}
\item  If $k \equiv 2, 6 \mod 8$ 
then 
there is at least one zero $z_k \neq \ii$ for which $|z_k| = 1$ and $d_E(z_k,\ii) = O(c_\ii^{-k})$ for some absolute constant $c_\ii>0$.  
\item If $k \equiv 0 \mod {12}$ then 
there is at least one zero $z_k \neq \rho$ for which $\Re\, z_k = -1/2$ and $d_E(z_k,\rho) = O(c_\rho^{-k})$ for some absolute constant $c_\rho>0$. 
\end{enumerate}
\end{enumerate}
\end{introtheorem}

 The method of proof of Theorem \ref{main2} is a combination of sign changing arguments and Rouch\'e's theorem applied to find zeros of $E_k^{a,b}$ close to those of a suitable main term, where we might miss or overcount a limited number of zeros that are close to the intersection of different circles in the corresponding configuration. In Proposition \ref{zerosat}, we study when $\ii$ or $\rho$ are themselves in $Z_k^{a,b}$ for sufficiently large $k$. 
The precise analytic estimates that are used in the proof also imply the following equidistribution result. 

 \begin{introprop} \label{equidist}
For fixed odd level $N$, the set of arguments 
$ \{ \arg(z) \colon z \in Z_{\Gamma(N)}, \Re(z)<0\}$ equidistributes in $\arg \A$ in the sense that 
\[ \sup_{\frac{\pi}{2} \leq \theta_1 \leq \theta_2 \leq \frac{2\pi}{3}} \Big| \frac{\#(\arg(Z_{k}) \cap [\theta_1,\theta_2])}{\# (Z_{k}\cap \{\Re(z)<0\})} - \frac{6}{\pi} (\theta_2-\theta_1)\Big| = O\Big(\frac{{1 + \log N}}{k}\Big) \to 0 \mbox{ as }k\to\infty\]
 with an absolute implied constant.
 \end{introprop}

This implies that the normalized argument 
counting measure (with multiplicities and weights) for $Z_k\cap \{\Re(z)<0\}$ converges weakly to the restriction of the Haar measure on $\A$, {and, combined with Hausdorff convergence, that the normalized zero counting measure on $Y(1)$ converges to the argument measure restricted to the image of $\A$ in $Y(1)$, cf. Corollary \ref{cortotconv}}. 
 
\begin{introremark*} Cuspidal modular forms relate to interesting divisors on moduli spaces (see, e.g., \cite{vdGK}); our work shows that non-cuspidal forms lead to intriguing configurations of geodesics, but only `in infinite weight', in the following sense. Consider the corresponding open (finite volume, non-compact) Riemann surface $Y_\Gamma \defeq \Gamma \backslash \Hh$ as a Riemannian manifold for the hyperbolic metric.   Then the collection of zeros of Eisenstein series, considered as divisor in $Y_\Gamma$, approaches, as the weight increases, a compact subset of a finite configuration of geodesics in $Y_\Gamma$.
\end{introremark*} 

\section{Preliminaries on Eisenstein series for $\Gamma(N)$} 

We discuss the various definitions of Eisenstein series for a principal congruence group, their relations, Fourier expansions, and the basis problem for the Eisenstein space.  

\subsection{The space of Eisenstein series} The dimension of the space of $\E_k(N)$ of Eisenstein series for $\Gamma(N)$  of even weight $k>2$ and level $N$ equals the number of cusps of $\Gamma(N)$, given as $$\sigma(1)=1, \sigma(2)=3, \mbox{ and } \sigma(N)\defeq \frac{N^2}{2} \prod\limits_{p \mid N} \left(1-\frac{1}{p^2}\right) \mbox{ for } N\geq 3.$$  The Eisenstein series of Hecke defined by 
\begin{equation} \label{hecke} G_k^{A,B}(z) \defeq \sum_{(m,n) \equiv (A,B)\, \mathrm{mod}\, N} \frac{1}{(m z + n)^k}, \end{equation} 
form a spanning set when $(A,B)$ runs through representatives modulo $N$, and a basis under the further `primitivity' restrictions that $(A,B,N)=1$ and by taking one representative   per equivalence class $(A,B) \sim (-A,-B)$ \cite[\S 1]{Hecke}. This last  index set corresponds to the set of cusps of $\Gamma(N)$; $(A,B)$ corresponds to the point $(A:B)$ in the orbit space $\Gamma(N) \backslash \PP^1(\Q)$. 

In this paper, we will focus instead on the alternative spanning set $\{E_k^{a,b}\}$ for $\E_k(\Gamma(N))$ given by Equation \eqref{defeis} with $(a,b)$ running through $(N^{-1}\Z/\Z)^2$. 
The two sets of generators are connected by a matrix given by a Finite Fourier Transform over the enumerating group $(\Z/N\Z)^2$: 
\begin{align} \label{fexp} E_k^{a,b}(z) &= \sum_{0 \leq A,B \leq N-1} \sum_{(m,n) \equiv (A,B)\,  \mathrm{mod}\, N} \e(aB-bA) \frac{1}{(mz+n)^k} = \sum_{0 \leq A,B \leq N-1} \chi_a(B) \overline{\chi_b(A)} \, G_k^{A,B}, \end{align} 
where $\chi_a(x)\defeq\e(a \cdot x)$ runs over the characters of $\Z/N\Z$ with varying $a$. By Fourier inversion,
\begin{equation*} G_k^{A,B}(z) = {N^{-2} } 
\sum_{a,b} \overline{\chi_a(B)} \chi_b(A) E_k^{a,b}(z), \end{equation*}
 
One may try to identify a suitable set of indices $I \subset (N^{-1}\Z/\Z)^2$ such that $\{E_k^{a,b} \colon (a,b) \in I\}$ forms a \emph{basis} for $\E_k(\Gamma(N))$. For Hecke's series $G_k^{A,B}$, the set of indices can be chosen to consist of `primitive' pairs, up to sign. 
{The same holds for $E_k^{a,b}$, as the following proposition states.}

 {

 \begin{proposition}
\label{basis}
The space $\E_k(\Gamma(N))$ has basis given by $E_k^{a,b}$ where $(a,b) \in (N^{-1}\Z/\Z)^2$ runs through the space of \emph{primitive} pairs $(a,b)$, i.e., such that $(Na,	Nb,N)=1$, up to sign. 
 \end{proposition}
 
\begin{proof} To simplify notation, write $v=(a,b)$, so $E_k^{a,b}=E_k^v$. We have seen above that $\E_k(\Gamma(N))$ is spanned by $E^v$ for $v$ ranging over the entire index space $(N^{-1}\Z/\Z)^2$, and that the dimension of $\E_k(\Gamma(N))$ equals the number of primitive index pairs up to sign. Thus, it suffices to show that every $E_k^v$ is a linear combination of the series indexed by primitive pairs.

For every positive divisor $d$ of $N$, define the following sub-vector space of  $\E_k(\Gamma(N))$: \[V_d=\operatorname{span}_{\Cc} \{E_k^v:v\in(N^{-1}\Z/\Z)^2,\ \ord(v)=d\}.\]
where $\ord(v)$ is the (additive) order of $v$ in $(N^{-1}\Z/\Z)^2$. The primitive indices are precisely those of order $N$, so it suffices to show that $V_d\subseteq V_N$ for every $d\mid N$, i.e., that $V_d$ filters $\E_k(\Gamma(N))$ ordered by division.
The result follows if we prove that $V_d\subseteq V_{pd}$ for every prime $p$ with $pd\mid N$; then successively multiplying $d$ by prime factors until reaching $N$ completes the proof.
	
Suppose $\ord(v)=d$. We recall the distribution relation \cite[Lemma~3.7(2)(iii)]{Kato} 
	\begin{equation}\label{eq:distr} \sum_{\substack{x\in((pd)^{-1}\Z/\Z)^2\\px=v}}E_k^x =p^{2-k}E_k^v. \end{equation}

We study the solutions $x\in((pd)^{-1}\Z/\Z)^2$ to the equation $px=v$. First of all, \begin{equation} \label{ordord} d = \ord(v) = \ord(px)\overset{(\ast)}{=}\ord(x)/\gcd(\ord(x),p) \mbox{ divides }\ord(x). \end{equation}  

Thus, if $p$ divides $d$, then $p$ divides $\ord(x)$, so $\gcd(\ord(x),p)=p$, and by equality $(\ast)$ in \eqref{ordord}, every solution $x$ to $px=v$ has order $\ord(x)=dp$. Then \eqref{eq:distr} implies $E_k^v\in V_{pd}$.
	
We are left to consider the case where $p$ does not divide $d$. Then one solution to $px=v$ is given by $x_0=tv$, where $t$ is an integer such that $pt = 1\mbox{ mod }d$, and $\ord(x_0)=d$ since $t$ is coprime to $d$. If $x_1 \neq x_0$ is any other solution, then $p(x_0-x_1)=0$, so all such solutions have order $pd$ (again, since $p$ and $d$ are coprime). Now \eqref{eq:distr} implies 
	\begin{equation}\label{eq:toiter}
		E_k^{tv}-p^{2-k}E_k^v\in V_{pd}.
	\end{equation}
	Choose $h\ge1$ such that $t^hv=v$.
	Applying \eqref{eq:toiter} successively to
	$v,tv,\ldots,t^{h-1}v$, we obtain
	\[
	(1-p^{(2-k)h})E_k^v
	=E_k^{t^hv}-p^{(2-k)h}E_k^v\in V_{pd}.
	\]
The constant $1-p^{(2-k)h}$ is non-zero for $k\ge4$, so we conclude
	that indeed $E_k^v\in V_{pd}$.
\end{proof}
}

\subsection{Fourier expansion, relations} 

Consider the `modular action' of a matrix $\left( \begin{smallmatrix} \alpha & \beta \\ \gamma & \delta \end{smallmatrix} \right) \in \Gamma(1)$ on a holomorphic function $f \colon \Hh \rightarrow \Cc$ given by 
\begin{equation} \label{modact} f|_{ \left( \begin{smallmatrix} \alpha & \beta \\ \gamma & \delta \end{smallmatrix} \right) }(z)\defeq f \left( \frac{\alpha z+ \beta}{\gamma z + \delta} \right) (\gamma z + \delta)^{-k}. \end{equation} 
We find the following transformation formula for our Eisenstein series: 
\begin{equation} \label{transfo} E_k^{(a,b)}|_\gamma(z) = E_k^{(a,b)\gamma}(z) \mbox{ for all } \gamma \in \Gamma(1). 
\end{equation} 
where we use the notation for the exponent $a,b$ as a row vector $(a,b)$ and we consider the usual action of multiplication on the right on row vectors by matrices in $\Gamma(1)$. (We are using $\gamma$ to denote both an entry of a matrix, as well as a matrix, but this should not lead to confusion.) This immediately implies the three first formulae in the next lemma, and the final one follows by direct calculation.

\begin{lemma} \label{eisrel} For all $k \geq 3$, the following relations hold for all $z \in \Hh$: 
\begin{align} 
& E_k^{-a,-b}(z) = (-1)^k E_k^{a,b}(z) \label{eisrel1} \\ 
& E_k^{a,b}(z+1) =E_k^{a,a+b}(z) \label{eisrel2}  \\
& E_k^{a,b}(-1/z) = z^k  E_k^{-b,a} (z) \label{eisrel3} \\
& \overline{E_k^{a,b}(z)} = E^{a,-b}_k (- \overline{z} ) \label{eisrel4} 
\end{align} 
\end{lemma} 

The Fourier expansion at the infinite cusp is given as follows. 

\begin{proposition}\label{eisfour} For $(a,b) \in (N^{-1}\Z/\Z)^2$ and positive integers $k$ and $n$, define a divisor sum as $$\sigma_{k-1}[a,b](n)\defeq \sum_{\substack{j \in \Z, j \mid n \\ j \equiv aN\, \mathrm{mod}\, N}} \hspace*{-3mm} \mathrm{sgn}(j) e(bn/j)  j^{k-1};$$ then, with $q_N=\e(z/N)$, the following hold. 
\begin{enumerate}
\item \label{eisi} Let $k \geq 3$; then 
$$
\frac{(k-1)!}{(2 \pi \ii)^k} \,   E_k^{a,b}(z) =  - \frac{1}{k} B_k(a) +\frac{1}{N^{k-1}} \sum_{n>0} \sigma_{k-1}[a,b] \, q_N^n.
$$ 
\item \label{eis1a} $ \widetilde E_k^{a,b} \defeq 
\frac{-k!}{(2 \pi \ii)^k} \,   E_k^{a,b}(z) \in \E_k(\Gamma)\cap \Q(\zeta_N)[\![q_N]\!]$ \textup{(}as Fourier series in $q_N$\textup{)}. 
\item \label{eisii} For even $k\geq 4$, $E_k^{a,b}$ does not vanish at any cusp. In particular, $E_k^{a,b}$ is never identically zero. 
\end{enumerate} 
\end{proposition}

\begin{proof} \mbox{ }

\eqref{eisi}  This follows easily by using the identity \[\sum_{n \in \Z} \frac{\e(a n)}{(z + n)^k} =  (-1)^{k}\frac{(2\pi \ii)^{k}}{(k-1)!} \sum_{j = 0}^\infty (j-\{ a \}+1)^{k-1} \e(z(j-\{a\}+1) \] for $k \geq 2$ even and $z \in \Cc \setminus \Z$. Alternative equivalent expressions can be found in \cite[Lemma 34]{BrunaultZudilin} and \cite[Prop.~3.10]{Kato}.

\eqref{eis1a} follows from \eqref{eisi}, using that $B_k(x) \in \Q[x]$ and $e(bn/j) = \zeta_N^{Bn/j} \in \Q(\zeta_N)$ for $j \mid n$. 

\eqref{eisii} First notice that the non-vanishing of $E_k^{a,b}$ at all cusps is equivalent to the non-vanishing of $E_k^{a,b}|_\gamma$ for all $\gamma \in \Gamma(1)$. By the transformation formula \eqref{transfo}, this is equivalent to the non-vanishing of the constant Fourier coefficient of $E_k^{(a,b)\gamma}$ at the infinite cusp. The result now follows from the fact that the Bernoulli polynomial $B_k$ has no rational roots for $k$ even. This was first proven by Inkeri \cite{Inkeri}; an easy proof that such $B_k(x)$ does not even have $2$-adic roots is given in \cite[Theorem 4]{PatelSiksek}. 
\end{proof}

\subsection{The norm modular form} Results concerning zeros of all Eisenstein series of given even weight $k$ can be conveniently formulated by using the following \emph{norm modular form}: 
\begin{equation} \mathcal N_k \defeq \prod_{(a,b) \in (N^{-1}\Z/\Z)^2} E_k^{a,b}, \end{equation} 
By Proposition \ref{eisfour}, this product is non-zero at the cusp $\infty$, and, in particular, is not the zero function. By the relations in Lemma \ref{eisrel}, it follows that $\mathcal N_k$ is a modular form for $\Gamma(1)$, of weight $kN^2$. Note that by the same lemma, the factor $E_k^{0,0}$ is itself also a modular form for $\Gamma(1)$, and is in fact the standard Eisenstein series $N^{-k} G_k^{0,0}$ of weight $k$ for $\Gamma(1)$, whose zeros where already studied in \cite{RSD}. 

We may talk about the `union of all zeros of all $E_k^{a,b}$' to not just mean the set $Z_{k,\Gamma(N)}$, but the corresponding $\Q$-divisor $\mathrm{div}(\mathcal N_k)$ on the modular curve $X(1)$, whose degree, i.e., sum of the multiplicities of the zeros counted with the correct weights at elliptic points, by the valence  formula, equals $kN^2/12$. In our convention, this is also the meaning of the symbol $\# Z_{k,\Gamma(N)}$. 

\begin{remark}
If $k \geq 3$ is \emph{odd}, then the Bernoulli polynomial $B_k(x)$ has a rational root exactly at $0,1/2$ and $1$, see  \cite{Inkeri}. Thus, $E_k^{a,b}$ is nonvanishing at $\infty$ except if $a=0$ (so $A=0$), or $N$ is even and $a=1/2$ (so $A=N/2$). For these last two values, in fact, $E_k^{0,0}$ and $E_k^{1/2,1/2}$ (for $N$ even) are identically zero, as follows from the first identity in Lemma \ref{eisrel}. The correct `norm modular form' in these cases is given by leaving out the identically zero Eisenstein series from the defining product. 
\end{remark}

 \section{Kluyver sums} \label{sectionkluyver} 
 
In this section, we lay the foundations for our work by discussing general Kluyver sums. 
 
\subsection{Definitions and examples}  
 
 \begin{definition} Let $N \geq 1$, let $G$ be a subgroup of $\mathrm{SL}_2(\Z/N\Z)$, and let $v,w \in (\Z/N\Z)^2$ be two row vectors. Then, for an integer $j \in \Z$, we define the \emph{Kluyver sum} as 
 \begin{equation} \label{verygenram} \hat\rho_G^{v,w}(j) \coloneq \sum_{g \in G} \zeta_N^{j vgw^{\top}}. \end{equation} 
 We also define the corresponding real version 
  \[
  \rho_G^{v,w}(j) \coloneq \hat\rho_G^{v,w}(j)+\hat\rho_G^{v,w}(-j) = 2 \sum_{g \in G} \cos(2 \pi j vgw^{\top}/N). \qedhere
  \]
\end{definition} 

\begin{remark} \label{anykluyver} 
The definition naturally fits into the following general framework: let $N \geq 1$, $r \geq 1$, $\pi \colon G \rightarrow \mathrm{GL}_r(\Z/N\Z)$ a mod-$N$ representation of a finite group $G$, and let $v,w \in (\Z/N\Z)^r$ be two row vectors. Then, for an integer $j \in \Z$, we define the \emph{Kluyver sum} as 
 \begin{equation}  \hat\rho_\pi^{v,w}(j) \coloneq \sum_{g \in G} \zeta_N^{j v\pi(g)w^{\top}}. \end{equation} 
and \emph{mutatis mutandis} for the real version. Here, we only need the $2$-dimensional case. 
\end{remark}

We collect some basic properties of these sums. 

\begin{lemma} \label{grl} Let $N \geq 1$, let $G$ be a subgroup of $\mathrm{SL}_2(\Z/N\Z)$, and let $v,w \in (\Z/N\Z)^2$ be two row vectors. Then
\begin{enumerate}
\item \label{grl1} For any $h \in \mathrm{SL}_2(\Z/N\Z)$,  $\hat\rho_{h^{-1}Gh}^{v,w}(j) = \hat\rho_{G}^{vh^{-1},wh^{\top}}(j),$ and $\rho_{h^{-1}Gh}^{v,w}(j) = \rho_{G}^{vh^{-1},wh^{\top}}(j).$ 
\item \label{grl2} $\hat\rho_G^{v,w}(j)$ is non-zero for some $j$ with $|j| \leq |G|$ and $\rho_G^{v,w}(j)$ is non-zero for some $j$ with $|j| \leq 2|G|$.  
\item \label{grl3} $\hat\rho_G^{v,w}(j)$ and $\rho_G^{v,w}(j)$ are non-zero for some $j$ dividing $N$.
\end{enumerate}

\end{lemma} 

\begin{proof}
The proof of \eqref{grl1} is a calculation. 

For proving \eqref{grl2}, we use the following general fact: if, given $n$ complex numbers $a_1,\dots,a_n$, the sum $\sum_{i=1}^n a_i^j$ vanishes for every $j=1,\dots,n$, then all the numbers $a_i$ themselves vanish. This is proven in \cite[p.\ 1]{Lenard} using determinants, but also follows from the Newton formulas, which imply that if all these power sums are zero, then also the elementary symmetric functions of the numbers are so, and hence the polynomial satisfied by the numbers is a monomial of degree $n$, implying that the $a_i$ are all zero themselves. This immediately implies the result, noticing that $\hat\rho_G^{v,w}(j)$ is the sum of $|G|$ powers, and $\rho_G^{v,w}(j) = \hat\rho_G^{v,w}(j) + \hat\rho_G^{v,w}(-j)$ is the sum of $2|G|$ powers (of roots of unity). 

For proving \eqref{grl3}, we notice that we can rewrite any sum $S_N(j) = \sum_{i \in I} \zeta_N^{j \cdot a_i}$ ($I$ a finite index set, $a_i$ integers) as follows: 
$$ S_N(j) = S_{N/(N,j)}(j/(N,j)) \stackrel{\ast}{=} S^*_{N/(N,j)}(1) =S^*_N((j,N)), 
$$
where $S^*_N(j)$ is $S_N(j)$ with $\zeta_N$ changed to a suitable other primitive $N$-th root of unity; indeed, in $\ast$, we have used that $N/(N,j)$ and $j/(N,j)$ are coprime, and raising a primitive root of unity to a power coprime to its order just changes the primitive root of unity into another such. Since (non)-vanishing of $S_N(j)$ does not depend on the choice of primitive $N$-th root of unity, we see that if non-vanishing happens at $j$, then also at $(j,N)$, which is a divisor of $N$. 
\end{proof} 

We now consider a short list of examples that will be important later. 

\begin{example} \label{gamma0RS}
We consider the case where 
$$G = \{\left( \begin{smallmatrix} \alpha & \beta \\ 0 & \alpha^{-1} \end{smallmatrix} \right) \colon \alpha \in (\Z/N\Z)^*, \beta \in \Z/N\Z \},$$ the `Borel subgroup' in $\mathrm{SL}_2(\Z/N\Z)$.
Taking $v=(A,B)$ and $w=(w_1,w_2)$, we get   
\begin{equation} \label{formgamma0RS} 
\hat\rho_G^{v,w}(j)  =  \sum_{\substack{\alpha \in (\Z/N\Z)^* \\ \beta \in \Z/N\Z}} \zeta_N^{\alpha Aw_1j + \alpha^{-1} B w_2 j + A w_2 \beta j} 
= S(Aw_1j,Bw_2j;N) \cdot S_N(Aw_2j), 
\end{equation}  
where $$S(a,b;N) \defeq \sum_{\alpha \in (\Z/N\Z)^*} \zeta_N^{\alpha a + \alpha^{-1} b}$$ is a classical Kloosterman sum, and 
$$ S_N(a) \defeq \sum_{\beta \in \Z/N\Z} \zeta_N^{a \beta} = \left\{ \begin{array}{ll} N & \mbox{ if } N \mid a, \\ 0 \mbox{ otherwise.} \end{array} \right. $$ 

Since the right hand side of \eqref{formgamma0RS} is invariant under $j \mapsto -j$ (by changing, in the summation, both $\alpha$ and $\beta$ to their negations), we have 
$ \rho_G^{v,w}(j) = 2 \hat\rho_G^{v,w}(j). $

An interesting case arises when $w_2=0$; then
$$ \hat\rho_G^{v,(w_1,0)}(j) = N  \sum\limits_{\alpha \in (\Z/N\Z)^*} \zeta_N^{j \alpha A},$$ which is a multiple of a classical Ramanujan sum
\begin{equation} \label{deframsum} c_n(q) \defeq \mathrm{Tr}^{\Q(\zeta_n^q)}_{\Q} ( \zeta_n^{q}) =  \sum\limits_{\alpha \in (\Z/n\Z)^*} \zeta_n^{\alpha q}
\end{equation} 
and then we find the evaluation 
\begin{equation} \label{kluif} \rho_G^{v,(w_1,0)}(j) = 2 \hat\rho_G^{v,(w_1,0)}(j) = 2N \mu\left(\frac{N}{(N,Aj)} \right)  \frac{\varphi(N)}{\varphi(N/(N,Aj))} \end{equation}
independent of $B$ and $w_1$. For this last evaluation, we used the first of the following two general identities for Ramanujan sums---that we will use repeatedly---
 \begin{equation} \label{kluyver} c_n(q) = \mu(n/(n,q))\frac{\varphi(n)}{\varphi(n/(n,q))} = \sum\limits_{d \mid (n,q)} \mu(n/d)d, \end{equation} see \cite[p.~410]{Kluyver}. 
\end{example}

\begin{example} \label{gamma1RS}
If 
$$G = \{\left( \begin{smallmatrix} 1 & \beta \\ 0 & 1 \end{smallmatrix} \right) \colon \beta \in \Z/N\Z \},$$ the `unipotent radical' in $\mathrm{SL}_2(\Z/N\Z)$, taking $v=(A,B)$ and $w=(w_1,w_2)$, we get from the previous calculation that    
$$
\hat\rho_G^{v,w}(j) =  \zeta_N^{Aw_1j +  B w_2 j} S_N(A w_2 j) =  \left\{ \begin{array}{ll} N \zeta_N^{Aw_1j +  B w_2 j} & \mbox{ if } N \mid Aw_2j, \\ 0 \mbox{ otherwise.} \end{array} \right. , $$
and
\[
\rho_G^{v,w}(j) =  \left\{ \begin{array}{ll} 2N \cos(2 \pi (Aw_1 +  B w_2) j/N) & \mbox{ if } N \mid Aw_2j, \\ 0 \mbox{ otherwise.} \end{array} \right. \qedhere
\]
\end{example} 

\begin{example} \label{LambdaRS}
For the `inversion' group 
$ G = \langle  
\left( \begin{smallmatrix} 0 & 1 \\ 1 & 0 \end{smallmatrix} \right) \rangle \leq \mathrm{SL}_2(\mathbf F_2), $
of order two, with $v=(A,B)$ and $w=(w_1,w_2)$, we find  
$ \hat\rho_G^{v,w}(j) =  (-1)^{j(Aw_1+Bw_2)} + (-1)^{j(Aw_2+Bw_1)}, $
and $\rho_G^{v,w}(j) = 2 \hat\rho_G^{v,w}(j)$. 
\end{example} 

A more exotic case is the following. 

\begin{example} \label{exotic} let $N=p \neq 2$ be a prime. There is a subgroup $G$ of $\mathrm{SL}_2(\mathbf F_p)$ isomorphic to  the quaternion group $Q_8$, generated by matrices $\left( \begin{smallmatrix} 0 & 1 \\ -1 & 0 \end{smallmatrix} \right)$ and   $\left( \begin{smallmatrix} \xi & \psi \\ \psi & -\xi \end{smallmatrix} \right)$
for any $\xi, \psi \in \mathbf F_p$ satisfying $\xi^2+\psi^2 =-1$. Choosing $w=(1,0)$ gives a finite list $$ \{ g(1,0)^{\top} \colon g \in G \} = \{ \pm (1,0)^{\top}, \pm (0,1)^{\top}, (\pm \xi, \psi)^{\top}, \pm (\psi, -\xi)^{\top} \}.$$ For $p=3$, we can choose $\xi=\psi=1$ and verify numerically that none of sums $\hat\rho_G^{v,(1,0)}(1)$ for the 81 choices of $v$ vanish. Since the above list is invariant under sign change, none of the sums $\rho_G^{v,(1,0)}(1)$ vanishes either. 
\end{example} 

\subsection{(Non-)vanishing of Kluyver sums} \label{detkappa} 
We now study the non-vanishing of Kluyver sums. For this, we define the following invariants, which by Lemma \ref{grl}\eqref{grl2}, are well-defined integers. 

\begin{definition} \label{defkappafinite} Let $N \geq 1$, let $G$ be a subgroup of $\mathrm{SL}_2(\Z/N\Z)$, and let $v,w \in (\Z/N\Z)^2$ be two row vectors. Then we define
\begin{align*} 
\kappa_G^{v,w} &\defeq \min \{ j>0 \colon \rho_G^{v,w}(j) \neq 0 \}, \\
 \kappa_G &\defeq \max \{ \kappa_{G}^{v,w} \colon v,w \in (\Z/N\Z)^2 \}. \qedhere
\end{align*}
\end{definition} 

We notice the following property that will be useful when applying the theory to congruence groups. It shows that in the definition of $\kappa_G$, we can replace one running variable by looping over the conjugates of $G$. 

\begin{lemma} \label{conjvect} $\kappa_G = \max \{ \kappa_{h^{-1} G h}^{v,(1,0)} \colon v \in (\Z/N\Z)^2, h \in \mathrm{SL}_2(\Z/N\Z) \}$.
\end{lemma} 

\begin{proof} We use interchangingly elements mod $N$ and their representatives in $\Z$.
It suffices to prove that, for fixed $j$, the values taken by $\rho_G^{v,w}(j)$ with $v,w$ running over all vectors modulo $N$, and the values taken by $\rho_{h^{-1}Gh}^{v,(1,0)}(j)$ with $v$ running over all such vectors and $h$ over the integer matrices of determinant $1$ modulo $N$, are the same. For this, it suffices to prove that the sets of exponents of $\zeta_N^j$ that occur in the corresponding $\rho_*^*$ are the same.  

Let $v,w$ run over all vectors. Write $w = \lambda \tilde w$, where $\tilde w = (\tilde w_1,\tilde w_2)$ is primitive modulo $N$, i.e., $(\tilde w_1,\tilde w_2, N)=1$. By the Euclidean algorithm, there exists a integer matrix $h$ of determinant $1$ such that  $\tilde w^{\top} = h(1,0)^{\top}.$  This means that the occuring exponents in $\rho_G^{v,w}$ are
\begin{align*} \{ v g w^{\top} & \colon g \in G, v,w \in (\Z/N\Z)^2 \} \\ &= \{ (\lambda v h) h^{-1}gh(1,0)^{\top} \colon g \in G, \lambda \in \Z/N\Z, v \in (\Z/N\Z)^2, h \in \mathrm{SL}_2(\Z/N\Z) \} \\ &= \{ \tilde v \tilde g (1,0)^{\top} \colon \tilde g \in h^{-1}Gh, \tilde v \in (\Z/N\Z)^2, h \in \mathrm{SL}_2(\Z/N\Z) \}, 
\end{align*}  
which are the occurring exponents in $\rho_{h^{-1} G h}^{\tilde v,(1,0)}$.
\end{proof}  

\begin{remark} We can define similar invariants $\hat \kappa_G^{v,w}$ and $\hat 
 \kappa_G$, using $\hat\rho$ instead of $\rho$ in the above definitions. We will also calculate these invariants. 
\end{remark} 

We notice that it follows from Lemma \ref{grl} that all these invariants divide the level $N$ and and upper bounded by $2|G|$. Which of these bounds is `better' depends on the situation at hand. 

\begin{example} \label{kappatrivialgroup} 
For the trivial group $G = \{1\}$, $\hat \kappa_{\{1\}}^{v,w}= 
 \hat \kappa_{\{1\}} = 1$. However, since $\rho^{v,w}_{\{1\}} = 2 \cos(2 \pi j vw^{\top}/N)$, 
$$\kappa_{\{1\}}^{v,w}= \left\{ \begin{array}{ll} 2 & \mbox{ if } vw^{\top} \in \pm N/4 + N\Z, \\ 1 & \mbox{ otherwise.} \end{array} \right. $$ 
In particular, the value is $1$ if $N$ is not divisible by $4$. This leads to 
\[
\kappa_{\{1\}} = \left\{ \begin{array}{ll} 2 & \mbox{ if } 4 \mid N,\\ 1 & \mbox{ if } 4 \nmid N. \end{array} \right.  \qedhere 
\]
\end{example} 

In the above example, for any $N \geq 2$, the upper bound $2|G|=2$ beats the upper bound $N$. In the other cases below, the second bound wins.

\begin{example}  \label{kappafinitegamma1} 
If $G$ is as in Example \ref{gamma1RS}, $\hat \rho_G^{v,w}(j)$ is non-vanishing if $N \mid A w_2 j$, and the smallest $j$ for which this happens is $\hat \kappa_G^{v,w} = M$ with $M \defeq N/(N,Aw_2)$. Choosing $A=w_2=1$ gives $M=N$, so the maximum is reached: $\hat \kappa_G = N$. 

For $\rho_G^{v,w}(j)$ to not vanish, we need both $M \mid j$ and $Aw_1 j + B w_2 j \notin \pm N/4 + N \Z.$ Writing $j = M  \tilde \jmath$ and $L=N/M$, we need $Aw_1 \tilde \jmath + B w_2 \tilde \jmath \notin \pm L/4 + L \Z$. If $4 \nmid L$, this is always true with $\tilde \jmath=1$ and $\kappa_G^{v,w} = M$. If $4 \mid L$,  and the condition fails for $\tilde \jmath = 1$, then it holds for $\tilde \jmath = 2$ (this is because, in general, if $x \pm L/4$ is divisible by $L$, then $2x \pm L/4$ cannot be). 
First of all, if $M=N$ (which happens, e.g., for $A=w_2=1$), then $L=1$ and $\tilde \jmath = 1$ suffices. Otherwise, $M \leq N/2$ and then even if $\tilde \jmath = 1$ does not work, $\tilde \jmath = 2$ will, and then still $j \leq N$. In conclusion,  $\kappa_G = N$ in this case, too. 
\end{example}

\begin{example} \label{kappafinitelambda} 
If $G$ is as in Example \ref{LambdaRS}, the constants are 
\[
\kappa_{\langle \left( \begin{smallmatrix} 0 & 1 \\ 1 & 0 \end{smallmatrix} \right) \rangle}^{v,w}= \hat \kappa_{\langle \left( \begin{smallmatrix} 0 & 1 \\ 1 & 0 \end{smallmatrix} \right) \rangle}^{v,w}= \left\{ \begin{array}{ll} 2 & \mbox{ if } v,w \in \{(0,1),(1,0)\}, \\ 1 & \mbox{ otherwise.} \end{array} \right.\qedhere
\]
\end{example} 

\begin{example} \label{kappafinitegamma0} We now assume $G$ is as in Example \ref{gamma0RS}. It appears to be quite subtle to determine in general the smallest $j>0$ for which the sum in \eqref{formgamma0RS} does not vanish. We already observed that all constants $\kappa_G^*,\hat \kappa_G^*$ divide $N$.  

If we look in greater detail, the second factor in \eqref{formgamma0RS}, $S_N(Aw_2j)$, is non-zero if and only if $$j =  M \tilde \jmath \mbox{ with } M \defeq N/L, L\defeq (A w_2,N)$$ for some integer $\tilde{\jmath}$, and this is a minimal requirement that $j$ should satisfy. 

For the entire Kluyver sum $\rho_G^{v,w}$ in \eqref{formgamma0RS} to be non-zero, we then further require that the Kloosterman sum $$S(A w_1 \tilde{\jmath}, B w_2 \tilde{\jmath}, L) \neq 0$$ is non-vanishing. We refer to \cite{IK} for general properties of Kloosterman sums and to Sali\'e \cite[p.\ 97]{Salie} for a study of non-vanishing of such sums. By using the multiplicative properties of the sums, the problem can reduced  to the case where the modulus is a prime power dividing $L$ (but one of the occurring first terms in the Kloosterman sum will change in this process): if $L=\prod p_i^{e_i}$ then \begin{equation} \label{multK} S(a,b;L) = \prod S(a,b_i,p_i^{e_i}) \end{equation}  for some $b_i$ determined inductively from the formula $S(a,\beta_1;r) S(a,\beta_2;s) = S(a,b; rs)$ for coprime $r,s$, where $b=\beta_1 s^2 + \beta_2 r^2$. The trivial fact that then $bj=\beta_1 j s^2 + \beta_2 j r^2$ for any $j$, implies that if Equation \eqref{multK} holds, then for any $j$, 
\begin{equation} \label{multKj} S(ja,jb;L) = \prod S(ja,jb_i,p_i^{e_i}).  \end{equation} 
Also, applying the formula to $L=rs$ with $r=p^{e_i}$ and $s$ coprime to $p$, we see that $\beta_1 s^2 + \beta_2 p^{2 e_i} = b$, and hence $(b,p^{e_i})=(\beta_1,p^{e_i})$ in the inductive process.  

First of all, we notice that if $L=p$ is prime, any Kloosterman sum $S(a,b;p) \in \Z[\zeta_p]$ reduces to $-1$ modulo $p$, hence is always non-zero. 
If $L=p^e$ for $p$ prime and $e \geq 2$, let $a=p^{e_a} a',b = p^{e_b} b',\tilde \jmath = p^{e_{\tilde \jmath}} \tilde \jmath'$ with $a'b',\tilde \jmath'$ coprime to $p$, and without loss of generality, assume $e_a\leq e_b \leq e$. Then $S(\tilde \jmath a , \tilde \jmath b; p^e) = p^{e_{\tilde \jmath} +e_a}S(a', p^{e_b-e_a} b', p^{e-e_a-e_{\tilde \jmath} }).$ This is non-zero (by the previous observation) if $ e_{\tilde \jmath}  = e-1-e_a$, or, in general, 
\begin{equation} \label{beste} e_{\tilde \jmath}  = e-1-\min(e_a,e_b).\end{equation} 
We now analyse whether a smaller value of $e_{\tilde \jmath}$ is possible. If we choose $e_{\tilde \jmath}< e-1-\min(e_a,e_b)$, then the Kloosterman sum is a non-zero multiple of $S(a',p^{e_b-e_a} b', p^m)$ for some $m \geq 2$. This sum is only non-zero if $e_b=e_a$ \cite[(22)]{Salie}. In that case, for $p \neq 2$, the sum is non-zero precisely if $a'b'$ is a quadratic residue modulo $p$. For $p=2$, there is a more complicated case distinction depending on $a'b'$ modulo $8$, see \cite[\S 3]{Salie}, but the problem is decidable in terms of such congruences. The conclusion is that, although for specific choices of $a,b$ (i.e., $A,B,w_1,w_2,N$), a smaller value is possible, there are always such choices for which \eqref{beste} is the smallest possible value.  

Collecting all these local terms, we see that $\tilde \jmath = \prod_{p_i^2 \mid L} p_i^{e_i-1}/(a,b,p_i^{e_i-1})$ makes  all Kloosterman sums non-vanishing. \end{example}

A succinct way of rewriting the conclusion of the above discussion is the following result. 

\begin{proposition} \label{kloosbound} For $G$ as in Example \ref{gamma0RS}, $(A,B), (w_1,w_2) \in (\Z/N\Z)^2$, we have 
\begin{equation} \label{kappagamma0upper} \kappa_G^{(A,B),(w_1,w_2)} \leq \frac{N}{\ell \cdot \mathrm{rad} 
((N,Aw_2)/\ell)}\ \mbox{ where } \ \ell \defeq (N,Aw_1,Aw_2,Bw_2)
\end{equation} and the following upper bound is reached in various families: 
\begin{equation} \label{kappagamma0noconj} 
\left\{ \begin{array}{l} \max\{ \kappa_G^{(A,B),(1,0)} \colon A,B \in [0,N-1] \}  = N/\mathrm{rad}(N); \\
\kappa_G = \max\{ \kappa_G^{v,w} \colon v,w \in (\Z/N\Z)^2 \} = N. 
\end{array} \right. 
\end{equation} 
\end{proposition} 

\begin{proof}
The first maximum is reached, e.g., by $A=1$, and the second by $A=w_2=1$. 
\end{proof} 

\begin{remark}  \label{remkappagamma0noconj} In the special case where $w=(1,0)$, the sum $S_N(Aw_2j)$ is always non-vanishing, the Kloosterman sum simplifies to a Ramanujan sum that we evaluated in \eqref{kluif}, and this immediately leads to a sharp result in that case, as follows. 
Set $d=(N,A)$ and $A'=A/d$, so that $N=Md$ and $A=A'd$. 
Using the identity in \eqref{kluif}, we find 
\begin{align*} c_N(jA) &= \mu(N/(N,jA)) \frac{\varphi(N)}{\varphi(N/(N,jA))} = \mu(M/(M,jA')) \frac{\varphi(N)}{\varphi(M/(M,jA'))} \\ &= \frac{\varphi(N)}{\varphi(M)} c_M(jA'), \end{align*} 
where now $M$ and $A'$ are coprime. Applying the same formula to $c_M(A'j)$ using this coprimeness, we find that 
\begin{equation} \label{cmj} c_N(jA) = c_M(jA') = c_M(j). \end{equation} We conclude that $c_N(jA)=0$ if and only if $c_M(j)=0$, if and only if $\mu(M/(M,j))=0$. 
Recalling that the M\"obius function is nonzero only on the squarefree numbers, for non-vanishing we need $j$ to be such that dividing by $(j,M)$ `removes' the powerful part of $M$, and this happens for the following smallest value of $j$: 
$$ \kappa_{G}^{(A,B),(1,0)} = M/\mathrm{rad}(M) \mbox{ where } M=N/(N,A).$$ 
This matches the upper bound in Proposition \ref{kloosbound} and the largest value is obtained when $(N,A)=1$, see \eqref{kappagamma0noconj}. Notice that this bound is not valid if we also allow conjugate subgroups of $G$; then $\kappa_G=N$ is the optimal upper bound. 
\end{remark}

\section{General congruence groups} A subgroup $\Gamma \leq \Gamma(1)$ is a congruence group of level $N$ if and only if there is an inclusion of group $\Gamma(N) \trianglelefteq \Gamma \leq \Gamma(1)$ and $N$ is minimal with this property. Such congruence groups of level $N$ correspond to finite subgroups $G \leq \mathrm{SL}_2(\Z/N\Z)$ by, given $\Gamma$, setting $G=\Gamma(N) \backslash \Gamma$ and conversely, given $G$, letting $\Gamma$ be the inverse image of $G$ under the natural projection map $\Gamma(1) \rightarrow \Gamma(N)\backslash \Gamma(1) \cong \mathrm{SL}_2(\Z/N\Z)$. Recall that we set $h_\Gamma = [\Gamma:\Gamma(N)] = |G|$.  Given an isomorphism class of an abstract finite group $H$, we call $H$ the \emph{type} of the congruence group $\Gamma$ if $G \cong H$. Notice that if $\Gamma$ is a congruence group of level $N$ and type $H$, then the same holds for its $\Gamma(1)$-conjugate subgroups $\gamma^{-1} \Gamma \gamma$ ($\gamma \in \Gamma(1)$). Distinct such conjugate groups are parameterized by $N_\Gamma \defeq \Gamma(1)/\mathrm{Nor}_{\Gamma(1)}(\Gamma)$ (using the orbit-stabilizer theorem for the action of conjugation on the set of subgroups). 

\subsection{The trace modular forms $E_{k,\Gamma}^{a,b}$} 
The space of Eisenstein series $\E_k(\Gamma)$ of weight $k$ for $\Gamma$ is a subspace of $\E_k(\Gamma(N))$. We consider the trace modular forms as in \eqref{traces}; the use of traces to construct modular forms occurs explicitly already in \cite{PeterssonSpur}, cf.\ \cite[VII, \S1.3]{Schoeneberg}. 

\begin{proposition} \label{Tab} let $\Gamma$ denote a congruence group of level $N$ and $k \geq 4$ an even integer.  
\begin{enumerate}
\item \label{Tab1} For any $\gamma \in \Gamma(1)$, $E_{k,\Gamma}^{a,b}|_{\gamma} = E_{k,\gamma^{-1}\Gamma\gamma}^{(a,b)\gamma}$. 
\item \label{Tab2} The trace forms $E_{k,\Gamma}^{a,b}$ from \eqref{traces} generated $\E_k(\Gamma)$, where $(a,b)$ runs over (a set of representatives of the orbit of $\Gamma$ acting on) $(N^{-1} \Z/\Z)^2$. 
\item \label{Tab2a} The dimension of $\E_k(\Gamma)$ is the number of cusps of $\Gamma$. 
\item \label{Tab4} $\widetilde E^{a,b}_{k,\Gamma} \defeq \frac{-k!}{(2 \pi \ii)^k} \,  E_{k,\Gamma}^{a,b} \in \E_k(\Gamma)\cap \Q(\zeta_N)[\![q_N]\!]$ \textup{(}as Fourier series in $q_N=\e(z/N)$\textup{)}. 
\item \label{Tab5} We have the identities $E_{k,\Gamma}^{-a,-b}(z) = (-1)^k E_{k,\Gamma}^{a,b}(z)$ and $\overline{E_{k,\Gamma}^{a,b}(z)} = E^{a,-b}_{k,\Gamma} (- \overline{z} )$. 
\end{enumerate}
\end{proposition} 
\begin{proof} 
\eqref{Tab1} From the definition and \eqref{transfo}, we find (using that $\Gamma(N)$ is normal in $\Gamma$) 
$$ E_{k,\Gamma}^{a,b}|_{\gamma} = \sum_{g \in \Gamma(N)\backslash \Gamma} E_k^{(a,b)g\gamma} = \sum_{\gamma^{-1} g \gamma \in \Gamma(N)\backslash \gamma^{-1}\Gamma\gamma} E_k^{(a,b) \gamma \cdot \gamma^{-1} g \gamma} = E_{k,\gamma^{-1} \Gamma \gamma}^{(a,b)\gamma}. $$ 

\eqref{Tab2} A projection onto the linear subspace $\E_k(\Gamma) \subseteq \E_k(\Gamma(N))$ is given the Reynolds operator $$P(f)\defeq\frac{1}{|G|} \sum_{\gamma \in G} f|_{\gamma}$$ for the modular action \eqref{modact} of the finite group $G= \Gamma(N) \backslash \Gamma$. Note that projection send a set of generators to a set of generators (possibly with more relations). Using \eqref{Tab1}, we find that $E_{k,\Gamma}^{a,b}$ are a non-constant multiple of the images of the Reynolds operator of our set of generators $E_k^{a,b}$ of $\E_k(\Gamma(N))$, i.e., 
$$ |G| \cdot P(E_k^{a,b}) = \sum_{\gamma \in G} E_k^{a,b}|_{\gamma} = \sum_{\gamma \in G} E_k^{(a,b)\gamma} = E_{k,\Gamma}^{a,b}.$$ 
Thus, $\{ E_{k,\Gamma}^{a,b} \}$, where $(a,b)$ runs over $(N^{-1} \Z/\Z)^2$ indeed generates $\E_k(\Gamma)$. In fact, to get a set of generators, it suffices to let $(a,b)$ runs over a set of representatives of the orbit of $\Gamma$ acting on $(N^{-1} \Z/\Z)^2$. 

\eqref{Tab2a} is implied by our assumption that $k \geq 4$ is even (in odd weight $k \geq 3$, the dimension would be zero if $-1 \in \Gamma$ and the number of regular cusps of $\Gamma$ if $-1 \notin \Gamma$); see \cite[\S 4.2, (4.3)]{DiamondShurman}.

\eqref{Tab4} follows rom Proposition \ref{eisfour}\eqref{eis1a}. 

\eqref{Tab5} follows directly from the corresponding relations \eqref{eisrel1} and \eqref{eisrel4}.  
\end{proof} 

\begin{definition} \label{gennormforms} Define the \emph{norm modular form} of weight $k$ for $\Gamma$ as 
$$ \mathcal N_{\Gamma,k} \defeq \prod_{\substack{(a,b) \in (N^{-1}\Z/\Z)^2 \\ \gamma \in N_\Gamma}} E_{k,\gamma^{-1} \Gamma \gamma}^{a,b}.$$ This is $\Gamma(1)$-invariant, non-zero, of weight $kN^2\cdot |N_\Gamma|$, and hence we can consider its zeros in $\F$ as a $\Q$-divisor on the modular curve  $X(1)$. 
\end{definition} 

\subsection{Kluyver sums for congruence groups}  We will apply the theory of Kluyver sums to congruence groups using the following definition. 

\begin{definition} \label{defgenram}
When $\Gamma$ is a congruence group of level $N$ and $v,w \in (\Z/N\Z)^2$ are two row vectors, we set 
$$ \rho_\Gamma^{v,w}(j) \defeq \ 2 \hspace*{-3mm} \sum_{g \in \Gamma(N) \backslash \Gamma} \hspace*{-3mm} \cos(2 \pi j/N v g w^{\top}). $$  
For $(a,b) \in (N^{-1}\Z/\Z)^2$, we set $$ \rho_\Gamma^{a,b}(j) \defeq \rho_\Gamma^{(Na,Nb),(1,0)} = \ 2 \hspace*{-7mm} \sum_{\left( \begin{smallmatrix} \alpha & \beta \\ \gamma & \delta \end{smallmatrix} \right) \in \Gamma(N) \backslash \Gamma} \hspace*{-7mm} \cos(2 \pi(\alpha a + \gamma b)j). $$
 We will also write $\kappa_\Gamma^*$ for the corresponding integers $\kappa_G^*$ in Definition \ref{defkappafinite}.
\end{definition}  

\begin{remark} \label{kappaconjonly}
A subgroup of $\Gamma(1)$ conjugate to a given $\Gamma$ corresponds to a subgroup of $\mathrm{SL}_2(\Z/N\Z)$ conjugate to $G$, and we have seen in Lemma \ref{conjvect} that looping over all conjugates in $G$ (or $\Gamma$) is the same as looping over all vectors $v,w$ modulo $N$. Thus, $\kappa_\Gamma$, as defined in the introduction, indeed equals $\kappa_G$.  
\end{remark} 

\begin{remark}
We recall from Lemma \ref{grl} that for any congruence group $\Gamma$ of level $N$, $\kappa_\Gamma$ divides $N$ and is upper bounded by $2h_\Gamma$. In particular, if $\Gamma$ is of prime level $N=p$, then $\kappa_\Gamma \in \{1,p\}$. By Lang's Theorem, conjugacy classes of subgroups of $\mathrm{SL}_2(\mathbf F_p)$ are the same over $\mathbf F_p$ and $\overline{\mathbf  F}_p$, and all conjugacy classes of subgroups and representatives for each class are know, first in the work of Gierster \cite{Gierster}, so in prime level one only needs to compute $\kappa_\Gamma$ for the families of finite groups on this list.
\end{remark}

\begin{example} \label{gammacosets} Recall the following form of the quotient group $G=\Gamma(N)\backslash \Gamma$ for some standard choices of congruence groups $\Gamma$. 
\begin{enumerate} 
\item  \label{gamma0Ncosets} 
If $\Gamma=\Gamma_0(N)$, $G$ is as in Example \ref{gamma0RS}. 
\item \label{gamma1Ncosets}
If $\Gamma=\Gamma_1(N)$, $G$ is as in Example \ref{gamma1RS}. 
\item \label{Lambdagroupcosets} If $\Gamma=\Lambda$, $G$ is as in Example \ref{LambdaRS}.
\end{enumerate} 
Therefore, the determination of $\kappa_G$ for the subgroups given in Section \ref{detkappa} at the same time determines the values of $\kappa_\Gamma$ for congruence groups $\Gamma$ of these types.  
\end{example}

\subsection{Constant Fourier coefficients and non-vanishing of trace Eisenstein series} 

A modular form in $\E_k(\Gamma)$ is identically zero if and only if it is a cusp form, i.e., the constant term of the Fourier expansion at all cusps is zero. Applying this to $\widetilde E_{k,\Gamma}^{a,b}$, the `value at a cusp' (by abuse of terminology; only non-vanishing is relevant) is the value of the constant Fourier coefficient of $\widetilde E_{k,\gamma^{-1} \Gamma \gamma}^{a',b'}$ at the infinite cusp, for a suitable $a',b'$ and $\gamma \in N_\Gamma$. Recall from Proposition \ref{eisfour} that the value of $\widetilde E_k^{a,b}$ at the infinite cusp equals the value of the Bernoulli polynomial $B_k(a)$. We conclude the following. 

\begin{lemma} Given a congruence group $\Gamma$ of level $N$ and an even weight $k\geq 4$, 
the constant Fourier coefficient of the trace modular form $\widetilde E_{k,\Gamma}^{a,b}$ at the cusp $\infty$ equals
\begin{equation} \label{bsum} \widetilde E_{k,\Gamma}^{a,b}(\infty) = \sum_{\left( \begin{smallmatrix} \alpha & \beta \\ \gamma & \delta \end{smallmatrix} \right) \in \Gamma(N) \backslash \Gamma}  B_k(\{\alpha a + \gamma b\}). \end{equation} 
for all $(a,b) \in (N^{-1}\Z/\Z)^2$. \qed
\end{lemma}  

To get a grip on the (non-)vanishing of this coefficient, we plug in the Fourier series \begin{equation} \label{berfour} B_k(x)= c \sum\limits_{j \neq 0} \e(j x)/j^k \ \mbox{ (valid for $k \geq 2$ and $x\in [0,1]$)}, \end{equation} where $c\defeq - k!/(2 \pi \ii)^k$ and rewrite the right hand side in \eqref{bsum} to find (collecting together positive and negative values of $j$)
\begin{equation} \label{vanrho} E_{k,\Gamma}^{a,b}(\infty) =  \sum_{j \neq 0} \frac{\hat\rho^{a,b}_\Gamma(j)}{j^k}  =  \sum_{j> 0} \frac{\rho^{a,b}_\Gamma(j)}{j^k} \end{equation} 
in terms of the Kluyver sums from Definition \ref{defgenram}.

The finiteness of $\kappa_\Gamma^{a,b}$ not only shows the non-vanishing of all Eisenstein series, but even non-vanishing at all cusps. 
\begin{corollary} For a congruence group $\Gamma$, $(a,b) \in (N^{-1}\Z/\Z)^2$ and even weight $k \geq 4$, none of the Eisenstein series $E_{k,\Gamma}^{a,b}$ vanish at any cusp. 
\end{corollary} 

\begin{proof} 
The non-vanishing of $E_{k,\Gamma}^{a,b}$ at every cusp follows from the non-vanishing of $E_{k,\Gamma}^{a,b}|_\gamma$ for every $a,b$ and $\gamma \in \Gamma(1)$, at the infinite cusp. We know that $E_{k,\Gamma}^{a,b}|_\gamma = E_{k,\Gamma'}^{a',b'}$ for some (possibly different) index $(a',b')$ and a (possibly different) congruence group $\Gamma'$.
By \eqref{vanrho}, 
$$ E_{k,\Gamma}^{(a',b')}(\infty) = \frac{\rho_\Gamma^{a',b'}(\kappa_\Gamma^{a',b'})}{(\kappa_\Gamma^{a',b'})^k} + O\left(\frac{1}{(\kappa_\Gamma^{a',b'})^{k+1}}\right),$$ and  
$E_{k,\Gamma}^{(a',b')}(\infty) \not \equiv 0$ is equivalent to $\kappa_\Gamma^{a',b'} < +\infty,$ which follows from Lemma \ref{grl}\eqref{grl2}. 
\end{proof} 

\begin{remark} 
Vandiver \cite[Thm.~1]{Vandiver} has shown a `von Staudt-Claussen' theorem for rational values of Bernoulli polynomials,  and applying this to \eqref{bsum}, we find that for all even $k \geq 4$, we have 
\begin{equation} \widetilde E_{k,\Gamma}^{a,b} (\infty) \in  h_\Gamma \sum_{\stackrel{\ell-1 \mid k}{\ell \nmid N}} \frac{1}{\ell} + \Z,
\end{equation}
which determines the fractional part of the constant Fourier coefficient. 
\end{remark} 

\subsection{Constant Fourier coefficients for specific congruence groups} \label{cfc}

We now calculate the constant Fourier coefficient in case of some standard congruence groups. 

\begin{example} \label{gammaN}
If $\Gamma=\Gamma(N)$, then we already know $\widetilde E_{k,\Gamma(N)}^{a,b}(\infty) = B_k(a).$ 
\end{example} 
 
\begin{example} \label{gamma1N}
If $\Gamma=\Gamma_1(N)$, by Example \ref{gammacosets}\eqref{gamma1Ncosets}, \eqref{bsum} gives $\widetilde E_{k,\Gamma_1(N)}^{a,b}(\infty) = N B_k(a).$
\end{example} 

\begin{example} \label{Lambdagroup} If $\Gamma=\Lambda$, by Example \ref{gammacosets}\eqref{Lambdagroupcosets}, \eqref{bsum} simplifies to $E_{k,\Lambda}^{a,b}(\infty) = 2 B_k$.
\end{example}

If $\Gamma=\Gamma_0(N)$, we have the following closed formula for the constant Fourier coefficient. 
\begin{proposition}
Let $k \geq 2$ be even and $a=A/N \in N^{-1} \Z/\Z$ for some $A \in [1,N-1]$. Set $M\defeq N/(N,A)$. Then 
\begin{equation} \label{explsum} \widetilde E_{k,\Gamma_0(N)}^{a,b}(\infty) = N \sum_{\alpha \in (\Z/N\Z)^*} B_k(\{ \alpha a\}) = \frac{B_k N}{M^{k-1}} \frac{\varphi(N)}{\varphi(M)}\prod_{\substack{p \mid M \\ p \text{\ \textup{prime}}}} (1-p^{k-1}). \end{equation} 
\end{proposition}

\begin{proof} 
In this case \eqref{bsum} simplifies to \begin{equation} \label{bkaa} \widetilde E_{k,\Gamma_0(N)}^{a,b}(\infty) = \sum_{\alpha \in (\Z/N\Z)^*} B_k(\{ \alpha a\}). \end{equation}
We plug in the Fourier series \eqref{berfour} for $B_k(x)$ and rewrite the right hand side as 
\begin{equation} \label{ckr} c\sum_{j \neq 0} \frac{c_N(jA)}{j^k} \end{equation} 
in terms of classical Ramanujan sums \eqref{deframsum}. 
Using \eqref{cmj} and the second identity in \eqref{kluyver}, we rewrite \eqref{ckr} as follows, using the notation $c'\defeq c \varphi(N)/\varphi(M)$: 
\begin{align} \label{simp} 
c' \sum_{j \neq 0} \frac{c_M(j)}{j^k} & = c' \sum_{j \neq 0} \frac{1}{j^k} \sum_{d \mid (M,j)} \mu(M/d) d  
 = c' \sum_{d \mid M} d \mu(M/d) \sum_{\substack{j\neq 0 \\  d \mid j}} \frac{1}{j^k} \nonumber \\ &= 2 \zeta(k) c' \sum_{d \mid M} \mu(M/d) \frac{1}{d^{k-1}}  \stackrel{(*)}{=} \frac{2 \zeta(k)}{M^{k-1}} c' \prod_{\substack{p \mid M \\ p \text{\ prime}}} (1-p^{k-1}). 
\end{align} 
To prove $(*)$, notice that both sides are multiplicative (the left hand side is a Dirichlet convolution of the multiplicative functions $\mu$ and $d \mapsto d^{1-k}$) and the identity is easy to prove for $M$ a prime power. 
\end{proof}

\section{Downward imaginary concentration}\label{section: Downward imaginary concentration}

To prove that zeros of Eisenstein series have bounded imaginary part---in a sense, are being `pushed away' by the infinite cusp---we will analyse them by truncation. This is a standard method, already used in \cite{RSD} in its easiest form. Since we will apply it to arbitrary congruence groups, we set it up in greater generality. 

\subsection{Main term, remainder, constant term} 
We start by defining main terms and remainder terms for the Eisenstein series, determined by a cut-off parameter $t$. 

\begin{definition} \label{cutdef} For $(a,b) \in (N^{-1} \Z/\Z)^2$, $k \geq 4 $ an integer and $t>0$ a positive real number or infinity, we define the following truncated functions for a congruence group $\Gamma(N)$:
\begin{enumerate} \setlength{\itemsep}{6pt}
\item `Main term'\ \  $ \displaystyle{M^{a,b,\leq t}_k(z) \defeq \sum\limits_{m^2+n^2 \leq t} \frac{\e(an-bm)}{(mz+n)^k}}$; 
\item `Remainder term' \ \   $ \displaystyle{R^{a,b,> t}_k(z) \defeq \sum\limits_{m^2+n^2 > t} \frac{\e(an-bm)} {(mz+n)^k}}$; 
\item `Constant term' \ \  $ \displaystyle{C_k^{a,b,\leq t} \defeq \sum\limits_{n^2 \leq t} \frac{\e(an)}{n^k}}$.
\end{enumerate} 
For a general congruence group $\Gamma$ of level $N$, recall that 
$E_{k, \Gamma}^{a,b} = \sum\limits_{\gamma \in G} E_k^{(a,b)\gamma}$ where $G = \Gamma(N)\backslash \Gamma$. Similar to the above, we define 
\[
M^{a,b,\leq t}_{k,\Gamma}(z) \defeq  \sum_{\gamma \in G} E_k^{(a,b)\gamma, \leq t}, \ 
R^{a,b,> t}_{k,\Gamma}(z) \defeq  \sum_{\gamma \in G} E_k^{(a,b)\gamma,> t}, \ 
C_{k,\Gamma}^{a,b,\leq t}  \defeq  \sum_{\gamma \in G} C_{k,\Gamma}^{(a,b)\gamma, \leq t}.  \qedhere
\]

\end{definition} 

\begin{remark} Notice 
$\displaystyle{ C_{k,\Gamma}^{a,b,\leq \infty} = \lim_{t \rightarrow + \infty} C_{k,\Gamma}^{a,b,t} = \lim_{y \rightarrow + \infty} E_{k,\Gamma}^{a,b}(\ii y) = E_{k,\Gamma}^{a,b}(\infty),} $
and this constant Fourier coefficient at the infinite cusp is given by formula \eqref{bsum}.
\end{remark} 

We immediately note the following crucial relation between (non)-vanishing of truncated constant terms and  of Kluyver sums.

\begin{proposition} \label{decaycontrol} 
If $\Gamma$ is a congruence group of level $N$, $k \geq 4$ an even weight and $(a,b) \in (N^{-1}\Z/\Z)^2$ an index, then 
\begin{enumerate} 
\item \label{u1} for any truncation parameter $t>0$, we have
$$C_{k,\Gamma}^{a,b,\leq t} = \sum_{0<|j|<\sqrt{t}} \frac{\hat\rho_\Gamma^{a,b}(j)}{j^k} = \sum_{0<j<\sqrt{t}} \frac{\rho_\Gamma^{a,b}(j)}{j^k}; $$
\item \label{u2} we have the asymptotics as $k \rightarrow + \infty$,  $$ \left| \sum_{j>0} \frac{1}{j^k} \rho_{\Gamma}^{a,b}(j) \right| \sim |\rho_\Gamma^{a,b}(\kappa_\Gamma^{a,b})| \frac{1}{(\kappa_\Gamma^{a,b})^k}; $$
\item \label{usim} for a sufficiently large cut-off parameter $t$, we have 
$$ \lim_{k \rightarrow + \infty} \sqrt[k]{|C_{k,\Gamma}^{a,b,\leq t}|} = 1/\kappa_\Gamma^{a,b}. $$ 
\end{enumerate} 
\end{proposition}

\begin{proof}

\eqref{u1} The truncated constant term of the `trace' Eisenstein series is by definition the sum of the truncation of the Eisenstein series of principal congruence groups that occur in the trace sum over the group $G$. With $C_k^{a,b,\leq t} = \sum_{0<|n|\leq \sqrt{t}} \e(an)/n^k$, we find the result. 

 For ease of reading the rest of the proof, we leave out $a,b,\Gamma$ from the notation. 
 
\eqref{u2} This follows from $$ \lim_{k \rightarrow + \infty} \left| 1 + \sum_{j>\kappa} \frac{\rho(j)}{\rho(\kappa)} \left( \frac{k}{j}\right)^k \right| = 1, $$
since $|\rho(j)|$ is absolutely bounded (e.g., by $h_\Gamma$). 

\eqref{usim} 
Take $t> \kappa$. Then the same reasoning as in \eqref{u2} applies: 
$$ \lim_{k \rightarrow + \infty} \sqrt[k]{|C_{k}^{\leq t}|} = \lim_{k \rightarrow + \infty} \frac{1}{\kappa} \sqrt[k]{\left| \rho(\kappa) + \sum_{j>\kappa}^{\sqrt{t}} \rho(j) \left(\frac{\kappa}{j}\right)^k \right|} =  \frac{1}{\kappa} \lim_{k \rightarrow + \infty} \sqrt[k]{|\rho(\kappa)|} = 1/\kappa,  $$
because $\rho(\kappa) \neq 0$. 
\end{proof}

\subsection{Estimates of remainders} \label{secrem}

In this subsection, we consider the size of the remainder of the Eisenstein series in the fundamental domain, after we delete the first terms up to a cut-off parameter.

We {first} 
deduce an estimate  for the length of vectors in the lattice $\Z \oplus z \Z$ in terms of $\Im(z)$ in great generality, because in later proofs, we will need to `tune' the constant $\alpha_\kappa$ in the statement to be large enough. 

\begin{lemma} \label{im2est} If $z \in \F$ and $|z| \geq \kappa$ \textup{(}in particular, if $\Im(z) \geq \kappa$\textup{)} for some $\kappa \geq 1$, then for all pairs of integers $(m,n)$, $$|mz+n|^2 \geq \alpha_\kappa \cdot (m^2+n^2), $$
where 
{
\begin{equation}
	\label{alphavalue}
	\alpha_\kappa = \frac{1}{2}\left(\kappa^2+1-\sqrt{(\kappa^2-1)^2+1}\right).
\end{equation}
}
\end{lemma} 
{\begin{proof}
		We can assume $(m,n)\neq(0,0)$. Since $|z|\geq\kappa$
		and $|\Re(z)|\leq 1/2$, for any $\alpha\in(0,1)$
		we have
		\begin{align*}
			|mz+n|^2-\alpha(m^2+n^2)
			&= (|z|^2-\alpha)m^2+(1-\alpha)n^2
			+2mn\Re(z)\\
			&\geq (\kappa^2-\alpha)m^2+(1-\alpha)n^2-|mn|\\
			&= \left(|m|\sqrt{\kappa^2-\alpha}
			-|n|\sqrt{1-\alpha}\right)^2\\
			&\quad +|mn|\left(
			2\sqrt{(\kappa^2-\alpha)(1-\alpha)}-1
			\right).
		\end{align*}
		The right-hand side is non-negative whenever
		\[
		(\kappa^2-\alpha)(1-\alpha)\geq\frac14.
		\]
Solving the corresponding quadratic equation for $\alpha$,
we find that the largest admissible value is precisely
$\alpha_\kappa$ as in \eqref{alphavalue}.
	\end{proof}
}

\begin{remark} Note that $\lim\limits_{\kappa \rightarrow \infty} \alpha_\kappa = 1$ and  $ \alpha_1 \geq 0.5$, $\alpha_2 \geq 0.9$, $\alpha_3 \geq {0.96}
$. \end{remark}

\begin{lemma} \label{upperR} Let $\Gamma$ denote a congruence group of level $N$. If $z \in \F$ and $|z| \geq \kappa$ \textup{(}in particular, if $\Im(z) \geq \kappa$\textup{)}, then for $k \geq 4$ and  $(a,b) \in (N^{-1} \Z/\Z)^2$, we have  
$$ |R^{a,b, >t}_{k,\Gamma}(z)| \leq h_\Gamma \frac{4(t+1)(2 \sqrt{t+1}+1)}{((t+1) \alpha_\kappa)^{k/2}}, $$
independently of $a,b$ {for any non-negative integer $t$}. 
\end{lemma} 

\begin{proof} Since $|R^{a,b, >t}_{k,\Gamma}(z)| \leq h_\Gamma \cdot \max\limits_{a,b} |R^{a,b, > t}_{k}(z)|$, it suffices to prove the upper estimate (that is independent of $a,b$) for $\Gamma=\Gamma(N)$.  
We insert the estimate from Lemma \ref{im2est} into the definition of $R^{a,b,>t}_k(z)$ and now fix $a,b$ and leave it out of the notation. We let $r_2(M)$ denote the number of pairs of integers $(m,n) \in \Z^2$ for which $M=m^2+n^2$, and note that {$r_2(M) \leq 4 \sqrt{M}$.} 
This gives
{
\begin{align*}
	|R^{>t}_k(z)|
	&\leq \sum_{m^2+n^2 \geq t+1}
	\frac{1}{(\alpha_\kappa(m^2+n^2))^{k/2}} \leq \frac{1}{\alpha_\kappa^{k/2}}
	\sum_{M\geq t+1}\frac{r_2(M)}{M^{k/2}} \leq \frac{4}{\alpha_\kappa^{k/2}}
	\sum_{M\geq t+1}\frac{1}{M^{(k-1)/2}} \\ & \leq \frac{4}{\alpha_\kappa^{k/2}}
	\left(
	\frac{1}{(t+1)^{(k-1)/2}}
	+\int_{t+1}^\infty
	\frac{\d x}{x^{(k-1)/2}}
	\right) = \frac{4}{\alpha_\kappa^{k/2}(t+1)^{(k-1)/2}}
	\left(1+\frac{2(t+1)}{k-3}\right)\\
	&\leq \frac{4}{\alpha_\kappa^{k/2}(t+1)^{k/2}}
	\left(\sqrt{t+1}+2(t+1)^{3/2}\right)\leq \frac{4(t+1)(2\sqrt{t+1}+1)}
	{((t+1)\alpha_\kappa)^{k/2}}.
\end{align*}
}
\end{proof} 

\subsection{Estimates of main terms, up to constant terms} 

We now estimate the main term minus constant terms in a region of $\F$ that is sufficiently far up along the imaginary axis. 

\begin{lemma} \label{upperMu}  Let $\Gamma$ denote a congruence group of level $N$. If $z \in \F$ and $\Im(z) \geq \kappa+c/k$ for a positive real constant $c \geq 0$, then for any $(a,b) \in (N^{-1} \Z/\Z)^2$, and weight $k$ with $$  k \geq \max \{ 2c^2/\kappa^2,4\}$$ we have  
$$ |M^{a,b, \leq t}_{k,\Gamma}(z)-C^{a,b, \leq t}_{k,\Gamma}| \leq h_\Gamma \frac{8 \sqrt{t} \zeta(k)}{\kappa^k} \me^{-c/\kappa}. $$
\end{lemma} 

\begin{proof} Since $$|M^{a,b, \leq t}_{k,\Gamma}(z)-C^{a,b, \leq t}_{k,\Gamma}| \leq h_\Gamma \cdot \max\limits_{a,b} |M^{a,b, \leq t}_{k}(z)-C^{a,b, \leq t}_{k}|,$$ it suffices to prove the upper estimate (that is independent of $a,b$) for $\Gamma=\Gamma(N)$.  

For any integers $m,n$, we have 
$$|mz+n| \geq |\mathrm{Im}(mz+n)| = m |\mathrm{Im}(z)| \geq m(\kappa+c/k),$$ and so we can estimate
\begin{align*} |M^{a,b,\leq t}_k(z)-C_a^{\leq t}|
&\leq \sum_{\substack{m^2+n^2 \leq t \\ m \neq 0}}  \frac{1}{m^k} (\kappa+c/k)^{-k} \\ 
& \leq \sum_{m \neq 0} \frac{1}{m^k} \cdot \#\{n \in \Z \colon n^2 \leq t \} \cdot \frac{1}{\kappa^k} (1+c\kappa^{-1}/k)^{-k} \\
& \leq \frac{2 \zeta(k) \cdot 2 \sqrt{t} }{\kappa^k} (1+c\kappa^{-1}/k)^{-k}.
\end{align*} 
Now, in view of the identity $
\exp(kx/(k+x))\le (1+x/k)^k
$
valid for all $x\gt 0$ we see that 
\[
(1+c\kappa^{-1}/k)^{-k} \leq \me^{-\frac{c}{\kappa}}\me^{\frac{c^2}{k\kappa^2+c\kappa}}\le 2\me^{-\frac{c}{\kappa}},
\] 
where the last inequality follows by our assumption that $k\ge \frac{2c^2}{\kappa^2}\gt \frac{c^2}{\kappa^2\log{2}}-\frac{c}{\kappa}$. This proves the result.
\end{proof}

We can combine this with the estimate of the remainder, to give the following criterion for downward concentration of the zeros.  

\begin{theorem} \label{condu} Let $\Gamma$ denote a congruence group of level $N$ and let $(a,b) \in (N^{-1}\Z/\Z)^2$ be an index. Suppose that for some $\kappa\geq 1$ and {some integer} $t>\kappa^2/\alpha_\kappa-1$ we have 
\begin{equation} \label{uOmega}  |C_{k,\Gamma}^{a,b,\leq t}| \gg_{\Gamma} \frac{1}{\kappa^{k}}.     
\end{equation}
Then for some constant $c=c^{a,b}_\Gamma$, depending only on $\Gamma$ and $(a,b)$, and for all even weights $k \geq 4$ and $z \in \F$ with $E_{k,\Gamma}^{a,b}(z)=0$ we have $\Im(z)<\kappa+c/k$. 
\end{theorem}

Condition \eqref{uOmega} means that $ |C_{k,\Gamma}^{a,b,\leq t}| \geq \tilde c_\Gamma \kappa^{-k}$ for some constant $\tilde c_\Gamma>0$, so the requirement is that the truncated constant Fourier coefficient doesn't decay faster than exponentially in $k$ (recall $\kappa^{-1}<1$). As we have seen in Proposition \ref{decaycontrol}, the decay rate is controlled by the (non-)vanishing of Kluyver sums. 

\begin{proof} 
Asumme, on the contrary, that $z \in \F$ with $E_{k,\Gamma}^{a,b}(z)=0$, and we have $\Im(z) \geq \kappa+c/k$ for some $c>0$ (to be determined) and sufficiently large $k$. We want $c$ and $k$ to satisfy 
\begin{equation} \label{ck}  k\ge 2c^2/\kappa^2. \end{equation} Then we can use the estimates from Lemma \ref{upperR} and \ref{upperMu} to get 
\begin{align*}  |E_{k,\Gamma}^{a,b}(z) - C_{k,\Gamma}^{a,b,\leq t} | &\leq |R^{a,b,>t}_{k,\Gamma}(z)| + |M^{a,b,\leq t}_{k,\Gamma}(z) - C_{k,\Gamma}^{a,b,\leq t} |  \\
&\leq h_\Gamma \frac{4(t+1)(2 \sqrt{t+1}+1)}{((t+1) \alpha_\kappa)^{k/2}} + h_\Gamma \frac{8 \sqrt{t} \zeta(k)}{\kappa^k} \me^{-c/\kappa}.
\end{align*} 
By plugging in a zero of $E_{k,\Gamma}^{a,b}$ on the left, this leads to a contradiction if the right hand side in this estimate is strictly smaller than $| C_{k,\Gamma}^{a,b,\leq t} |$. We can rearrange that condition into a lower bound for $c$ as follows, 
\begin{equation}\label{eq: LowerBound c}
c > -\kappa  \log \left( \frac{1}{h_\Gamma 8 \sqrt{t} \zeta(k)} \left( | C_{k,\Gamma}^{a,b,\leq t} | \kappa^k  - 4 h_\Gamma (t+1) (2 \sqrt{t+1}+1) \left( \frac{\kappa}{\sqrt{(t+1) \alpha_\kappa}} \right)^k\right) \right).  
\end{equation}
\emph{if} the logarithm in this expression exists, which, for sufficiently large $k$, leads to condition 
\begin{equation} \label{liminfgamma} \liminf_{k \rightarrow + \infty} \kappa^k |C_{k,\Gamma}^{a,b,\leq t}| - 4 h_\Gamma (t+1) (2 \sqrt{t+1}+1) \left( \frac{\kappa}{\sqrt{(t+1) \alpha_\kappa}} \right)^k > 0.
\end{equation}  
By our assumption that $t>\kappa^2/\alpha_\kappa-1$, the second term tends to zero exponentially fast in $k$. By assumption \eqref{uOmega}, the first term stays larger than a strictly positive constant for $k$ tending to infinity. Therefore, condition \eqref{liminfgamma} is satisfied. With $\zeta(k)$ monotonously decreasing to $1$ as $k$ tends to infinity, this allows us to define a lower bound for $c$ only depending on $\Gamma$ and $a,b$ (and a fixed truncation level). 

Since we can make sure that \eqref{ck} holds by choosing $k$ sufficiently large, this establishes the result for $k$ large enough w.r.t.\ the fixed parameters. Since for fixed $\Gamma$, we have missed at most finitely many weights $k$ and hence finitely many zeros of the Eisenstein series in $\F$, we can increase the constant $c$ to have $\Im(z)<\kappa+c/k$ also for those zeros, and the Theorem is proven.  
\end{proof} 

\subsection{Downward imaginary concentration of the zeros}

We can now analyse the conditions in Theorem \ref{condu} for different given groups $\Gamma$, indices $(a,b)$, and for different choices of the cutoff parameter $t$ and $\kappa$.  As noticed before, the upper bounds on the remainder and  the`main term minus constant term' are easy and follow from the ones for principal congruence groups. The interesting challenge is to provide a lower bound for the (truncated) constant terms, as in \eqref{condu}, and this is determined by (non-)vanishing of Kluyver sums. 

\begin{proof}[Proof of Theorem \ref{mainIm} (except optimality)] The bounds on $\kappa_\Gamma$ follow from Lemma \ref{grl}\eqref{grl2}\&\eqref{grl3}. 
The condition in Theorem \ref{condu} hold for $\kappa = \kappa^{a,b}_\Gamma$ by choosing the truncation parameter to satisfy 
$$ t > \max \{ \kappa^{a,b}_\Gamma, (\kappa_\Gamma^{a,b})^2/\alpha_{\kappa_\Gamma^{a,b}} -1\}. $$ 
In fact, this gives at the level of an individual Eisenstein series 
 \begin{equation} \label{imupab} \Im(z) <  \kappa^{a,b}_\Gamma + O(1/k) \mbox{ for all } z \in Z_{k,\Gamma}^{a,b}. \end{equation} 
The exact values of the constants $\kappa_\Gamma$ for the groups under consideration have already been determined in Section \ref{detkappa}, Examples \ref{kappatrivialgroup}--\ref{kappafinitegamma0}. 
\end{proof} 

\begin{example} \label{excondugammaN} To see how this works in practice, fix $N$, consider $\Gamma=\Gamma(N)$ and set $a=A/N$ for some $A \in [0,N-1]$. 

\begin{itemize} 
\item First, assume that $N$ is not divisible by $4$ or $a \neq \pm 1/4$. We know we can choose $\kappa = 1$, and we choose $t=2$. 
Then indeed, $|C_k^{a,b,\leq 2}| = 2| \cos(2 \pi A/N)|$, and observe that the smallest possible value of $2| \cos(2 \pi A/N)|$ occurs when $A=\floor{N/4}$ or $A=\floor{N/4}+1$. Without loss of generality, assume that such minimum occurs at $A=\floor{N/4}$. Then, using the inequality $|\sin(\pi y)|\ge 2|y|$ valid for all $y\in[-1/2,1/2]$, we see that for $N\gt 3$ (so $2\{N/4\}/N\lt 1/2$), we have
    \begin{equation}\label{eq: LowerBound uA}
    |C_k^{a,b,\leq 2}|\ge 2\Big|\sin\Big(\frac{2\pi}{N}\Big\{\frac{N}{4}\Big\}\Big)\Big|\ge 4\Big|\frac{2}{N}\Big\{\frac{N}{4}\Big\}\Big|\ge \frac{2}{N},
    \end{equation}
and it is clear that this bound also extends for $N\ge 1$.
\item The above reasoning cannot work for $N$ divisible by $4$ and $a=\pm 1/4$, because $C_k^{\pm 1/4,\leq 3}=0$. But in this case, we can choose $\kappa=2$ and $t = 4$. With these choices, using $\alpha_2>0.9$, we indeed have $t+1=5>\kappa^2/\alpha_\kappa = 4/0.9$. Also, 
$$|C_k^{\pm \frac{1}{4},b,\leq 4}| = 2| \cos(2 \pi \cdot 1/4) + 1/2^k \cos(4 \pi \cdot 1/4)| = \frac{1}{2^{k-1}} =\frac{2}{\kappa^k} \gg \frac{1}{\kappa^{k}}.$$  (Notice that for $a = \pm 1/4$, no value $\kappa<2$ can satisfy \eqref{uOmega} whatever large cut-off value $t$ one choses.)\qedhere
\end{itemize} 
\end{example} 

\begin{remark}
The proof of optimality of the bound in Theorem \ref{mainIm} will be completed after the proof of Theorem \ref{mainconfig}, cf.\ Example \ref{opt}. 
\end{remark}

In the case of $\Gamma(N)$ with $N$ not divisible by $4$, we can analyse the inequalities from the above reasoning and the proof of Theorem \ref{condu} in detail to give  an explicit value for a  constant $c=c_{\Gamma(N)}$. 

\begin{proof}[Proof of Proposition \ref{explicit} (part 1: explicit bound on $\Im(z)$)] Suppose $N\ge 3$ is not divisible by $4$ and assume that $k\ge 2c^2$ (so that $c,k$ satisfy \eqref{ck}). Following the proof of Theorem \ref{condu} and using \eqref{eq: LowerBound uA}, we see that it suffices to choose $c$ such that
\[
c\gt -\log\Big(\frac{1}{8\sqrt{2}\zeta(4)}\Big(\frac{2}{N}-\frac{54}{(1.5)^{k/2}}\Big)\Big).
\]
Now, since $k\ge 2(\log(13N))^2\gt \frac{2}{\log(1.5)}\log(54 N)$, we see that 
$
{2}/{N}-{54}/{(1.5)^{k/2}}\gt {1}/{N}.
$
Therefore, it is enough to require $c\gt \log(8\sqrt{2}\zeta(4)N)$, so $c=\log(13 N)$ suffices. 
\end{proof}

The proof of the second part of Proposition \ref{explicit} will be given in the next section. 

\section{Transcendence} 

\subsection{General finiteness results} In this section, we will prove Theorem \ref{tran}. 
Recall that a \emph{CM-point} in $\Hh$ is a point belonging to some imaginary quadratic number field. We aim to use T.~Schneider's result that if both $z$ and $j(z)$ are algebraic, then $z$ is a CM-point. Algebraicity of $j(z)$ at a zero $z$ of a given modular form $f$ depends on the Fourier coefficients. For Eisenstein series, these all lie in a cyclotomic field, whose Galois action permutes Eisenstein series. The theory of complex multiplication will allow us to show that if an Eisenstein series has an algebraic zero, then there is also a (possibly different) Eisenstein series with an algebraic zero of small modulus (based on an idea in \cite{Kohnen}).

\begin{proof}[Proof of Theorem \ref{tran}] 
Assume $E_{k,\Gamma}^{a,b}(z)=0$ for some algebraic $z \in \overline \Q$ and some $a,b,k$. By Proposition \ref{eisfour} (ii) we know that $ \widetilde E_{k,\Gamma}^{a,b} \in \Q(\zeta_N)[\![q_N]\!]$. 

For an integer $r$ coprime to $N$, the action of $\sigma_r \in (\Z/N\Z)^* \cong \mathrm{Gal}(\Q(\zeta_N)/\Q)$ given by $\sigma_r(\zeta_N)\defeq\zeta_N^r$, acting componentwise on the Fourier series of $\widetilde E_k^{a,b}$, is given by   
$$ \sigma_r(\widetilde  E_k^{a,b}) = \widetilde E_k^{a,rb} = \widetilde E_k^{(a,b) \gamma_r} \mbox{
where }\gamma_r = \left( \begin{smallmatrix} 1 & 0 \\ 0 & r \end{smallmatrix} \right) \in \mathrm{GL}_2(\Q) \cap \mathrm{M}_2(\Z).$$
{Returning to the general congruence group $\Gamma$ of level $N$, write $G=\Gamma(N) \backslash \Gamma = \pi(\Gamma) \leq \mathrm{SL}_2(\Z/N\Z)$, where $\pi \colon \M_2(\Z) \rightarrow \M_2(\Z/N\Z)$ is the reduction modulo $N$. Using  
the definition of Eisenstein series for general $\Gamma$ of level $N$ as a trace, we find  
$$  \sigma_r(\widetilde  E_{k,\Gamma}^{a,b}) = \sum_{g \in G} \widetilde  E_k^{(a,b)g \gamma_r} = \sum_{g \in G} \widetilde  E_k^{(a,br) \gamma_r^{-1} g \gamma_r}.$$
We claim this is also an Eisenstein series, albeit on a possibly different arithmetic group $\Gamma'$ of the same level $N$. The method we use is inspired by \cite[Prop.~2.1]{BN}. It suffices to notice that, although $\gamma_r$ is not necessarily invertible over $\Z$, since $r$ is invertible modulo $N$, its image modulo $N$ in $M_2(\Z/N\Z)$ is. Hence we can set 
\[ G' \defeq \pi(\gamma_r)^{-1} G \pi(\gamma_r) \mbox{ and } \Gamma' \defeq \pi^{-1} (G') \cap \Gamma(1), \]
so that, noting $G' \leq \mathrm{SL}_2(\Z/N\Z)$, we indeed find $\Gamma(N) \leq \Gamma' \leq \Gamma(1)$ with 
\[ \sigma_r(\widetilde  E_{k,\Gamma}^{a,b}) = \widetilde E_{k,\Gamma'}^{a,br}.  \]
Now 
\begin{equation} \label{nfortab} F \defeq \prod_{\gamma \in \Gamma \backslash \Gamma(1)} \prod_{\sigma \in \mathrm{Gal}(\Q(\zeta_N)/\Q) }  \sigma(\widetilde E_{k,  \Gamma}^{(a,b)}|_\gamma)\end{equation} }
is a product of Eisenstein series over various congruence groups of level $N$, and is itself a modular form for $\Gamma(1)$ with rational Fourier coefficients. (The need to also take the product over the Galois conjugates, not just the $\Gamma(1)$-cosets, arises from the later application of complex multiplication theory and is not necessary to prove that algebraic zeros are CM.) 

Choose suitable positive integers $u,v,w$ such that $G\defeq F E_4^u E_6^v / \Delta^w$ is of weight zero, where $E_4,E_6, \Delta$ are the classical Eisenstein series and discriminant function with rational Fourier coefficients. Since $G$ is a modular function for $\Gamma(1)$ with at most a pole at the cusp, it can be written as $P(j(z))$ for some polynomial $P \in \Cc[j]$, where $j$ is the $j$-invariant with rational Fourier coefficients. Since $G$ and $j$ have Fourier coefficients in $\Q$, in fact, $P \in \Q[j]$. If $F(z)=0$ and $j(z) \notin \{0,1728\}$ (the two excluded values corresponding to the CM-points $\ii$ and $\rho$), then $G(z)=0$, and so $j(z)$ is a zero of $P(j)$, and hence is algebraic. T.~Schneider's theorem \cite[IIc]{Schneider} then implies that $z$ is a CM-point. If $j(z) \in \{0,1728\}$, $z \in \{\ii, \rho\}$ is also CM.

We now appeal to the theory of complex multiplication \cite[\S 7]{Coxprimes}. 
The set $\{z_1,\dots,z_h\}$ (including $z$) of representatives modulo the action of $\Gamma(1)$ of points in $\Hh$ that satisfy a quadratic equation of discriminant $D$ (determined by $z$), coprime integer coefficients and positive leading coefficient is in bijection with the class group of the order $\mathcal O_D$ in the number field $K \defeq \Q(\sqrt{D})$ (since $\mathcal O_D$ is not necessarily maximal, this means the group of \emph{proper} ideal classes) \cite[Thm.\ 7.7]{Coxprimes}. Under this bijection, the equation $az^2+bz+c$ corresponds to the class of the fractional ideal $\langle 2a, b+\sqrt{D} \rangle$. The trivial ideal class corresponds to the point $\alpha_D$ given by 
\begin{equation} \alpha_D\defeq  \left\{ \begin{array}{ll} \sqrt{D/4} & (D = 0 \mbox{ mod } 4), \\ -1/2+\sqrt{D}/2 & (D = 1 \mbox{ mod } 4). \end{array} \right. \end{equation} 
We notice that $\alpha_D \in \mathcal F$, and the imaginary part of $\alpha_D$ is 
$ \Im\, \alpha_D = \sqrt{|D}|/2.$
At the same time, $\{j(z_i)\}$ is the full set of Galois conjugates of $j(z)$ over $K$. 
Since $P(j(z))=0$ and $P(j) \in \Q[j]$, all $j(z_i)$ are also zeros of $P$, so all $z_i$ are also zeros of $F$. In particular, $P(j(\alpha_D))=0$ and, given the form of $F$ in \eqref{nfortab}, there is an $a',b'$ (maybe different from $a,b$) and a congruence group $\Gamma'$ of level $N$ (maybe different from $\Gamma$) such that $\widetilde E_{k,\Gamma'}^{a',b'}(\alpha_D)=0.$ 

{We claim that $\kappa_{\Gamma'} = \kappa_\Gamma$. This follows from the construction of $\Gamma'$, whose image modulo $N$ is conjugate to $G$, and $\kappa_{\Gamma}$ only depending on $G$ and being the same for groups conjugate to $G$, see, e.g., Remark \ref{kappaconjonly}.} 

We conclude from Theorem \ref{mainIm} that $\Im(\alpha_D) = \sqrt{|D}|/2 \leq \kappa_\Gamma + c_\Gamma/k$ and thus 
\begin{equation} \label{nid} |D| \leq  \lfloor (2 \kappa_\Gamma + 2 c_\Gamma/k)^2 \rfloor. \end{equation} for an appropriate constant $c_\Gamma$. Choosing $k$ so large (depending on the constant, that is, on the group $\Gamma$) so that this upper bound equals $4\kappa_\Gamma^2$ then implies the result; the number of excluded $k$ is finite, and hence so is the number of excluded Eisenstein series, and hence, by the finiteness of $Z_{k,\Gamma}^{a,b}$, so is the number of excluded potentially algebraic zeros. 

It remains to establish the explicit values of potential algebraic zeros for small values of $\kappa_\Gamma$ in sufficiently large weight.  

For $\kappa_\Gamma=1$, we find an upper bound $|D| \leq 4$. This implies $|D|=3$ or $|D|=4$, corresponding to $z=\rho$ or $z=\ii$ (note that these discriminants have class number one, so that these two points are the only ones in the $\Gamma(1)$-orbit). 

For $\kappa_\Gamma=2$, we find an upper bound $|D| \leq 16$. We consider all possible way of writing $D=D_0 f^2$ with a fundamental discriminant $D_0$ and conductor $f$, and we find the possibilities in the top two lines of Table \ref{orders}. 
\begin{table}[t] 
\begin{tabular}{l|cccccccc}
$D$ & $-3$ & $-4$ & $-7$ & $-8$ & $-11$ & $-12$ & $-15$ & $-16$ \\
\hline
$f$ & $1$ &  $1$ &  $1$ &  $1$ &  $1$ &  $2$ & $1$ & $2$ \\
$h$ & $1$ &  $1$ &  $1$ &  $1$ &  $1$ &  $1$ & $2$ & $1$ \\
 $\{z_i\}_{i=1}^h$ & $\rho$ & $\ii$ & $\frac{1+\sqrt{-7}}{2}$ & $\sqrt{-2}$ & $\frac{-1+\sqrt{-11}}{2}$ &  $\sqrt{-3}$ &  $\left\{ \begin{array}{l} \frac{-1+\sqrt{-15}}{2} \\ \frac{-1+\sqrt{-15}}{4} \end{array}\right.$ & ${2 \ii}$ \\
\end{tabular}
\caption{All imaginary orders of discriminant $D$ with $|D| \leq 16$, with conductor $f$, class number $h$, and representatives for the ideal classes in $\F$.} \label{orders} 
\end{table} To compute the class number of a non-maximal order, one may use the relation between that class number and the one of the corresponding maximal order, see, e.g., \cite[Thm.\ 7.28]{Coxprimes}. All of this has been automated: for such orders of discriminant $\texttt D$, e.g., in SageMath \cite{Sage} one may use \texttt{D.class\_number()} from \texttt{sage.schemes.elliptic\_curves.cm}.   
The other lines of the table indicate the class number, a representative $z_1=\alpha_D$ for the trivial ideal class, and, in the unique case of class number two, also the unique representative in $\F$ of the non-trivial ideal class, which may be computed instantaneously using {\tt QuadraticField(-15).class\_group().gens()} in SageMath.  \qedhere
\end{proof}

\subsection{Explicit results}

In principle, one can compute an explicit lower bound on the weight $k$ so that only algebraic zeros on a given list can occur (and let this list be determined by $\kappa_\Gamma$ as in the above proof).  This is stated for $\Gamma=\Gamma(N)$ with $N$ not divisible by $4$ in Proposition \ref{explicit}, that we now prove. 

\begin{proof}[Proof of Proposition \ref{explicit} (part 2: algebraic zeros)] The only place in the proof where we needed to choose $k$ sufficiently large was to guarantee that the upper bound in \eqref{nid} equals $4 \kappa_\Gamma^2$. For $\Gamma(N)$ with $N$ not divisible by $4$, for $k>2(\log(13N))^2$ we can use the value $c_{\Gamma(N)} = \log(13N)$ as given in the proof of the first part of the proposition. 
\end{proof}  

\begin{corollary} \label{absupper} Fix a level $N$. The total number of algebraic zeros of all Eisenstein series of sufficiently large (w.r.t.\ $N$) even weight $\geq 4$ and all congruence groups of level $N$ does not exceed $5N^4.$
\end{corollary} 

\begin{proof} Looking at the above proof, for sufficiently large weight $k$, the number of such algebraic zeros is upper bounded by $\sum h(\mathcal O)$ where the sum runs over all imaginary quadratic orders $\mathcal O$ of discriminant upper bounded by $4 \kappa_\Gamma^2$ in absolute value.  A very crude way of estimating this goes as follows: the class number of an order of discriminant $-D_0 f^2$ ($-D_0$ a fundamental discriminant) is, by \cite[Thm.\ 7.28]{Coxprimes}, upper bounded by $h(-D_0) \Psi(f)$ with $\Psi(f) = f \prod_{p \mid f} (1+1/p)$ the Dedekind $\Psi$-function. Now $h(-D_0) \leq |D_0|/2$ and $\Psi(f) \leq f^2$. If $\Gamma$ is of conductor $N$, then $\kappa_\Gamma \leq N$, so we find the indicated bound by summing \[\sum h(\mathcal O) \leq  \sum_{|D_0|f^2 \leq 4\kappa_\Gamma^2} \frac12 |D_0| f^2 \leq \sum_{x \leq 4N^2} x/2 \leq 5N^4. \qedhere \]  
\end{proof} 
\begin{remark} The dependence on $N$ can be improved by using $h(-D_0) \leq \sqrt{|D_0|}/\pi (\log |D_0|)/2$ (cf.\ \cite{Ramare}) and $\Psi(f) \leq \me^\gamma f \log \log f$ for $f \geq 31$ (cf.\ \cite{Sole}), leading to an upper bound in Corollary \ref{absupper} of the form $\ll N^{3 + \epsilon}$ with an implicit constant. 
\end{remark}

\subsection{Algebraic zeros for principle congruence groups of level divisible by 4}
Our objective is to prove Proposition \ref{evenN}, i.e., to show that $E^{a,b}_k(z)\neq 0$ for the CM points in Theorem \ref{tran} (ii) with $z\neq \ii,\rho$. The argument is based on studying the irrationality measure of $y_n \defeq \arctan(\sqrt{n})/\pi$ for $n\in\bb{Z}_{>0}$. We begin by studying when $\arctan(\sqrt{n})$ is a rational multiple of $\pi$. We set 
\[
    \alpha_n\defeq    \me^{\ii \pi y_n} = \frac{1+\ii \sqrt{n}}{\sqrt{1+n}},
    \]
\begin{lemma}\label{lem: arctan(sqrt(n))/pi irrational}
    Let $n$ be a non-negative integer. Then $y_n \in\bb{Q}$ if and only if $n\in \{0,1,3\}$.
\end{lemma}
\begin{proof}
    Assume that $y_n \in\bb{Q}$. Then, 
    $
   \alpha_n
    $    is a root of unity, and $\bb{Q}(\me^{\ii \pi y})/\bb{Q}$ is an extension of degree $\le 4$. This shows that $(1+\ii \sqrt{n})^m$ is real for some positive integer $m$ such that $\vph(m)\le 4$, and this gives us a finite list for $m$. A finite computation reveals that the only possibilities are $n=0,1,3$. 
\end{proof}
\begin{lemma}
    For any non-negative integer $n$,  $y_n$ has finite irrationality measure.
\end{lemma}
\begin{proof}
    We take our inspiration from \cite{IrrationalityMeasure}. If $n=0,1,3$, then $y_n$ is a rational number, so it has irrationality measure $1$. For other $n$, $\log{\al_n}/\log(-1)=y_n$. Since, by the previous lemma, $y_n$ is irrational, the logarithmic form
    $
    L=\be_1 \log{\al_n} - \be_2 \log(-1)
    $
    is non-vanishing for any pairs of integers $(\be_1,\be_2)\neq (0,0)$. The Baker--W\"ustholz theorem \cite{BW} gives a lower bound of the form 
    $
    \log{|L|}\gt -C_1
    $
    for some $C_1=C_1(n)\gt 0$ only depending on $\al_n$, so only on $n$. This shows that for any positive integers $p,q$, 
    \[
    \Big|y_n-\frac{p}{q}\Big| = \Big|\frac{\log{\al_n}}{\log(-1)}-\frac{p}{q}\Big| \gt \frac{\me^{-C_1}}{q}\gt \frac{1}{q^{C_2}}
    \]
    for some $C_2=C_2(n)$. 
\end{proof}
Observe that by Lemma \ref{lem: arctan(sqrt(n))/pi irrational}, $|\sin(k\arctan(\sqrt{n}))|$ is non-vanishing when $n\not\in \{0,1,3\}$. We show more, namely, that it cannot decrease exponentially in $k$. 
\begin{lemma}\label{lem: LowerBoundFor sin(k arctan)}
     Let $n$ be a non-negative integer different from $0,1,3$, and let $\be$ denote the irrationality measure of $y_n$. Then, for every sufficiently large integer $k$, we have the following inequalities:
     \begin{enumerate}
         \item $|\sin(k\arctan(\sqrt{n}))|\gt \frac{1}{k^{\be}}.$
         \item $1-\cos(k\arctan(\sqrt{n}))\gt \frac{1}{k^{2\be}}.$
         \item $1+\cos(k\arctan(\sqrt{n}))\gt \frac{1}{k^{2\be}}.$
         \item $|\cos(k\arctan(\sqrt{n}))|\gt \frac{1}{k^{\be}}.$
     \end{enumerate}
\end{lemma}
\begin{proof}
    Let $x\defeq k\arctan(\sqrt{n})$, and observe that for any integer $m$, 
    $
    |\sin(x)|=|\sin(x-\pi m)|.
    $
    Let $m\defeq \floor{x+\f12}$. Then, combining the fact that $|y-\floor{y+\f12}|\le \frac{1}{2}$ for any real $y$ with the inequality $|\sin(\pi y)|\ge 2|y|$ valid for $-\f12 \le y\le \f12$, we have
    \[
    |\sin(x-\pi m)|=\Big|\sin\Big(\pi\Big(\frac{x}{\pi}-m\Big)\Big)\Big|\ge 2\Big|\frac{x}{\pi}-m\Big|=2k\Big|\frac{\arctan\sqrt{n}}{\pi}- \frac{m}{k}\Big|\gt \frac{1}{k^{\be}},
    \]
    where the last inequality holds for sufficiently large $k$, proving (i). Now (ii) follows immediately from (i) via the identity $1-\cos{x}=2\sin^2(x/2)$. To prove (iii), note that for all $y\in [0,2\pi]$, $1+\cos{y}\ge (y-\pi)^2/6$. Let $m'\defeq \floor{x}$, so that $x-2\pi m'\in [0,2\pi]$. Then 
    \[
    1+\cos{x}=1+\cos(x-2\pi m')\ge \frac{(x-2\pi m'-\pi)^2}{6}=\frac{k^2\pi^2}{6}\Big| \frac{\arctan(\sqrt{n})}{\pi}-\frac{2m'+1}{k}\Big|^2\gt \frac{1}{k^{2\be}}
    \]
    for sufficiently large $k$. Finally, (iv) follows by writing $\cos(x)=\sin(x+\frac{\pi}{2})$, and then doing a completely analogous procedure as in the proof of (i).
\end{proof}
We are ready to prove the following result, which implies Proposition \ref{evenN}. 
\begin{theorem}
    For all sufficiently large even integer $k$, and any level $N$, $E_k^{a,b}(z)\neq 0$ for $z\in \{\sqrt{-2}, \sqrt{-3}, \frac{-1+\sqrt{-7}}{2}, \frac{-1+\sqrt{-11}}{2},\frac{-1+\sqrt{-15}}{2}, \frac{-1+\sqrt{-15}}{4} \}.$
\end{theorem}
\begin{proof} We will use the shorthand notation $u_j= \cos(2\pi j/N)$.
First assume that $u_A\neq 0$. Observe that if $|z| \neq 1$, then by Lemma \ref{upperR},
\[
|E_k^{a,b}(z)-2u_A|=\Big|\frac{2u_B}{z^k}+R_k^{a,b,\gt 1}(z)\Big|
\]
decreases exponentially in $k$. This shows that for sufficiently large $k$, $|E_k^{a,b}(z)-2u_A|\lt |2u_A|$, i.e., $E_k^{a,b}(z)\neq 0$
for $z\in\{\sqrt{-2}, \sqrt{-3}, \frac{-1+\sqrt{-7}}{2}, \frac{-1+\sqrt{-11}}{2},\frac{-1+\sqrt{-15}}{2}\}$. Next, for the point $z_0\defeq \frac{-1+\sqrt{-15}}{4}$ (lying on the unit circle),  we let
\[
M_0\defeq 2\Big(u_A+\frac{u_B}{z_0^k}\Big)=2(u_A+u_B \me^{\ii k\arctan{\sqrt{15}}}),
\]
and this can only be zero if $\arctan{\sqrt{15}}$ is a rational multiple of $\pi$, but that is false by Lemma \ref{lem: arctan(sqrt(n))/pi irrational}. Now, let $R_0\defeq E_k^{a,b}(z_0)-M_0=R_k^{a,b,\gt 1}(z_0)$,
{

Since $|z_0|=1$ and $\Re(z_0)=-1/4$, we have
\[
|mz_0+n|^2=m^2-\frac{mn}{2}+n^2
\geq \frac34(m^2+n^2).
\]
Proceeding as in the proof of Lemma \ref{upperR}
with $t=1$, we obtain
\[
|R_0|
\leq \frac{8(2\sqrt{2}+1)}{(3/2)^{k/2}}
=O\Big(\Big(\frac23\Big)^{k/2}\Big).
\]
}
Thus, $R_0$ decreases exponentially in $k$. Therefore, if $|u_A|\neq |u_B|$, we see that 
\[
|M_0|\ge 2\big| |u_A|-|u_B|\big|\gt |R_0|
\]
for sufficiently large $k$, showing that $E_k^{a,b}(z_0)\neq 0$. If $|u_A|=|u_B|$, then by Lemma \ref{lem: LowerBoundFor sin(k arctan)}, we have
\[
|M_0|=2|u_A|\sqrt{2(1\pm\cos(k\arctan\sqrt{15}))}\gt \frac{|u_A|}{k^{\be_0}},
\]
where $\be_0$ is the irrationality measure of $\arctan(\sqrt{15})$. Since $R_0$ decreases exponentially in $k$, $|R_0|\lt |M_0|$ for sufficiently large $k$, showing again that $E_k^{a,b}(z_0)\neq 0$.

Now consider the case where $u_A=0$, and without loss of generality assume that $A=N/4$. Since $u_{A-B}=v_B$ and $u_{A+B}=- v_B$ with $v_j\defeq \sin(2\pi j/N)$, we have 
\[
E_k^{a,b}(z)=2\Big(\frac{u_B}{z^k}+v_B\Big(\frac{1}{(z+1)^k}-\frac{1}{(z-1)^k}\Big)-\frac{1}{2^k}+\frac{u_{2B}}{(2z)^k}\Big)+R_k^{a,b,\gt 2}(z).
\]
We again consider two separate cases.

\textbf{Case 1: $\boldsymbol{u_B\neq 0}$. }
When $z\in\{\sqrt{-2},\sqrt{-3},\frac{-1+\sqrt{-15}}{4}\}$, a routine computation via Lemma \ref{upperR} shows that
$|E_k^{a,b}(z)-{2u_B}/{z^k}| \lt |{2u_B}/{z^k}|$
for sufficiently large $k$. Assume that $z_1\defeq \frac{-1+\sqrt{-7}}{2}$, and let  
$
M_1\defeq 2({u_B}/{z_1^k}+{v_B}/{(z_1+1)^k}).
$
Then, by Lemma \ref{upperR}, $|E_k^{a,b}(z_1)-M_1|\le C\cdot 2^{-k}$ for some $C\gt 0$, showing that if $|u_B|\neq |v_B|$, then
$|M_1|\ge {| |u_B|-|v_B||}/{2^{k/2}} \gt {C}/{2^k}\ge |E_k^{a,b}(z_1)-M_1|.
$
If $|u_B|=|v_B|$, then again by Lemma \ref{lem: LowerBoundFor sin(k arctan)},
\[
\Big| \frac{u_B}{(-\f12+\f{\sqrt{7}}{2}\ii)^k}+\frac{v_B}{(\f12+\f{\sqrt{7}}{2}\ii)^k}\Big|\ge 2^{-k/2} |u_B| \cdot |\me^{ik\arctan\sqrt{7}}\pm\me^{-ik\arctan\sqrt{7}}|\gt \frac{2^{-k/2}|u_B|}{k^{\be_1}},
\]
where $\be_1$ is the irrationality measure of $\arctan(\sqrt{7})$. This shows that $|E_k^{a,b}(z_1)-M_1|\le |M_1|$, and we conclude that $E_k^{a,b}(z_1)\neq 0$. The case $z\in\{\frac{-1+\sqrt{-11}}{2},\frac{-1+\sqrt{-15}}{2}\}$ is completely analogous. 

\textbf{Case 2: $\boldsymbol{u_B= 0}$. }
We start with $z_2\defeq \sqrt{-2}$, and we claim that in this case the main contribution is given by 
\[
M_2\defeq v_B\Big(\frac{1}{(z_2+1)^k}-\frac{1}{(z_2-1)^k}\Big)=4v_B 3^{-k/2}\sin(k\arctan\sqrt{2} )\ii.
\]
By Lemma \ref{lem: LowerBoundFor sin(k arctan)}, $|M_2|\gt 3^{-k/2}|v_B|/k^{\be_2}$, where $\be_2$ is the irrationality measure of $\arctan\sqrt{2}$. Since $R_2\defeq E_k^{a,b}(\sqrt{-2})-M_2 \le C\cdot 2^{-k}$ for some $C\gt 0$ and $v_B\neq 0$, for sufficiently large $k$, $|R_2|\lt |M_2|$, showing that $E_k^{a,b}(z_2)\neq 0$.

If $z_3\defeq \sqrt{-3}$, then it may be the case that $\sin(\arctan\sqrt{3} k)=\sin(\pi k/3)=0$ when $k\equiv 0\mod{3}$, but in that case, we just pick the main term to be $-2^{1-k}$, and it is easy to show that for sufficiently large $k$, the remainder is smaller than this:
\[
|E_k^{a,b}(\sqrt{-3})+2^{1-k}|\le \frac{2}{(2\sqrt{3})^k}+\sum_{n^2+m^2 \ge 5}\frac{1}{|m\sqrt{-3}+n|^k} \lt 2^{1-k}.
\]
When $k\not\equiv 0\mod{3}$, then $|\sin(\pi k/3)|=\sqrt{3}/2$, so that the main term is 
\[
M_3\defeq 2v_B\Big(\frac{1}{(1+\sqrt{3}\ii)^k}-\frac{1}{(1-\sqrt{3}\ii)^k}\Big)=\pm 2v_B 3^{-k/2}\cdot \sqrt{3}\ii ,
\]
and again $|E_k^{a,b}(z_3)-M_3|\lt |M_3|$ for sufficiently large $k$. In any case, this shows that $E_k^{a,b}(z_3)\neq 0$.

Now, if $z=z_1=\frac{-1+\sqrt{-7}}{2}$, then since $u_B=0$, we have $v_B=\pm 1$, so the main contribution in absolute value is $|\frac{2v_B}{(z_1+1)^k}|=2^{1-k/2}$, and since all the other terms decrease as $O(2^{-k})$, then for sufficiently large $k$, 
$|E_k^{a,b}(z_1)-{2v_B}/{(z_1+1)^k}|\lt|{2v_B}/{(z_1+1)^k}|,
$
i.e., $E_k^{a,b}(z_1)\neq 0$. When $z\in\{\frac{-1+\sqrt{-11}}{2},\frac{-1+\sqrt{-15}}{2},\frac{-1+\sqrt{-15}}{4}\}$, the argument is completely analogous as for $z_1$.
\end{proof}

{

\begin{remark}
The statement of Proposition \ref{evenN} is optimal in the sense that all three points can occur as zeros for arbitrarily large even weights. Indeed, the classical Eisenstein series $E_k$ vanishes at $\ii$ when $k\equiv 2\pmod 4$, and at $\rho$ when $k\equiv 2\pmod 6$ (cf.\ also Proposition \ref{zerosat} infra). These are the only possibilities if $N$ is not divisible by $4$ and $k$ is sufficiently large. 
	
To see that $2\ii$ can also occur for infinitely many weights if $N$ divisible by $4$, consider $a=1/4$ and $b=1/2$. For every even $k\geq 4$, we have
$E_k^{1/4,1/2}(2z)=2^{-k}E_k^{1/2,1/2}(z).$
By the transformation formula \eqref{eisrel3} and \eqref{eisrel1}, we have
$E_k^{1/2,1/2}(-1/z)=z^k E_k^{1/2,1/2}(z).$ Setting $z=\ii$, we obtain
$(1-\ii^k)E_k^{1/2,1/2}(\ii)=0.$ Thus, whenever $k\equiv 2\pmod 4$, we have 
\[E_k^{1/4,1/2}(2\ii)=2^{-k}E_k^{1/2,1/2}(\ii)=0.\]
This is the only further possibility if $4 \mid N$ and $k$ is sufficiently large. 
\end{remark}
}

\section{Convergence to a configuration of geodesic segments}  

\subsection{Comparison of Euclidean and hyperbolic metric in compact subsets of $\Hh$} In this section, we will be talking about convergence of sequences and sets in $\Hh$, and we first briefly discuss why our results are the same in the Euclidean and hyperbolic metric; we will use the Euclidean metric in all the proofs. 
Recall that on $\Hh$ the Euclidean metric is $d_E(z,w)\defeq |z-w|$ and the hyperbolic metric $d_H(z,w) \defeq 2 \mathrm{arctanh}(|z-w|/|z-\overline{w}|).$
Although the Euclidean and hyperbolic metric can be quite different in all of $\Hh$ (think of the distance from $\ii$ to a point $y \ii$ with $y \rightarrow 0$: the Euclidean distance is bounded, but the hyperbolic distance diverges), in bounded compact regions, such as $$\Omega_\kappa \defeq \overline \F \cap \{ \Im(z) \leq \kappa \},$$ there are two-point comparison theorems for $d_E$ and $d_H$, that we now recall from \cite[Thm.\ 2(i)]{MaMaMi}; we will apply this to $\overline \F_\Gamma = \Omega_{\kappa_\Gamma}$. As ingredients, we need to know that the hyperbolic distance is defined via the density $\lambda(z) |\mathrm d z|$ with $\lambda(z)=1/\Im(z)^2 = 4/(z - \overline z)^2$.  Also recall that in $\Omega_\kappa$, 
$ \sqrt{3}/2 \leq \Im(z) \leq \kappa$, so for any two points $z,w \in \Omega_\kappa$, $\min\{\lambda(z),\lambda(w)\} \leq 4/3$ and $\max\{\lambda(z),\lambda(w)\} \geq 1/\kappa^2$. We also need to compute the domain constant 
$$ \eta(\Omega_\kappa) = \sup_{z \in \Omega_\kappa} |\partial_z \log \lambda(z)|/\lambda(z) =  \frac12 \sup_{z \in \Omega_\kappa} \Im(z) = \kappa/2. $$ Since this is bounded, $\Omega_\kappa$ is universality perfect, and \cite[Thm.\ 1]{MaMaMi} implies that 
\begin{equation} \label{twopointcomp} 1 - \me^{-\kappa d_H(z,w)} \leq \frac34 \kappa d_E(z,w) \ \mbox{ and }\ \frac{1}{\kappa} d_E(z,w) \leq \me^{\kappa d_H} -1. 
\end{equation}
As an application, we get, for example, the following. 
\begin{lemma}  If $d_E(z,w) \asymp 1/k$, then also $d_H(z,w) \asymp 1/k \ (z,w \in \Omega_\kappa).$
\end{lemma} 
\begin{proof}
Assume $c_1/k \leq d_E(z,w) \leq c_2/k$ for some constants $c_1,c_2>0$ From \eqref{twopointcomp}, we conclude that 
$ \log(1+c_1'/k)/\kappa \leq d_H(z,w) \leq -\log (1-c_2'/k)/\kappa$
for some constants $c_1',c_2'>0$.
The upper bound is further upper bounded by $c_1''/k$, and the lower bound is further lower bounded (e.g., using $\log(1+x) \geq x/\sqrt{1+x}$ for $-1<x<0$) by $c_2''/k$, and this proves the result.
\end{proof}

\subsection{Limits of sets} 

We first clarify different notions of set convergence that we will use, see, e.g., \cite[Ch.\ 4]{RockWet}. 
Recall that the \emph{Hausdorff} distance between two closed non-emtpy sets $X$ and $Y$ in $\Cc$ is 
$\delta(X,Y) \defeq \sup_{z \in \Cc} |\inf_{x \in X} |x-z| - \inf_{y \in Y} |y-z||.$
Suppose we are given a collection of subsets $\{X_k\}_{k \geq 1} \subseteq \Cc$. We study convergence of closed non-empty sets with respect to the Hausdorff distance using the following related concepts. 
The \emph{outer limit} (or `Kuratowski limit superior') of the collection is 
$$ \LS X_k \defeq \{ z \in \Cc \colon \exists k_j \rightarrow + \infty, X_{k_j} \ni z_j \rightarrow z \}.$$
The \emph{inner limit} (or `Kuratowski limit inferior') of the collection is $$ \LI X_k \defeq \{ z \in \Cc \colon \exists k_{\min}, \forall k>k_{\min}, \exists X_{k} \ni z_k \rightarrow z \}.$$
If $\LS X_k = \LI X_k$, we denote this by $\LH X_k$ and call it the \emph{Kuratowski limit}. We will use that 
if $X_k$ is a collection of sets contained in a \emph{bounded} domain of $\Cc$, then convergence in Hausdorff distance of $X_k$ to a set $X$ in the domain is equivalent to $\LH X_k = X,$ see \cite[4.13]{RockWet}. 

 \subsection{Hausdorff limit of the zeros}

Fix a congruence group $\Gamma$ and an index $(a,b) \in (N^{-1}\Z/\Z)^2$. In this subsection, only the weight $k$ will vary, so we use the shorthand notation 
$$ F_k(z) \defeq E_{k,\Gamma}^{a,b}(z). $$  Also, we will assume throughout that $k \in 2\Z_{\geq 2},$
in particular when writing limits as $k$ increases, without stating this explicitly. Furthermore, let $D$ denote a (sufficiently small) connected open domain in $\Hh$ containing $\overline \F_\Gamma$. {
	We choose $D$ sufficiently small such that $
	|mz+n|^2\geq \frac14(m^2+n^2)
	$
	for all $z\in D$ and $(m,n)\in\Z^2$, which is possible by Lemma \ref{im2est}.
}

To study the limit of the zero sets of the Eisenstein series $F_k$ as the (even) weight $k$ tends to infinity, we adapt  a method to determine the limit of zero sets of finite linear combinations of $k$-th powers of analytic functions as $k$ tends to infinity used by Sokal \cite{Sokal} in combinatorics and statistical physics. We are in a slightly different situation: Eisenstein series are infinite sums (over a lattice), but linear combinations of powers of particularly easy functions (of the form $1/(mz+n)$ for integers $m,n$), and we need only study the zeros in the compact set $\overline \F_\Gamma = \{ z \in \overline \F \colon \Im(z) \leq \kappa_\Gamma\}$; this mixture of generalisation and specialisation allows us to tweak the method of proof from \cite[Theorem 1.5]{Sokal}. 

\begin{proposition} \label{sokal} Rewrite the functions $F_k \colon \Hh \rightarrow \Cc$ \textup{(}$k \in 2\Z_{\geq 2}$\textup{)} in the form 
\[
F_k(z) = \sum_{\mathbf{v} \in \mathbb I} \frac{c_{\mathbf v}}{f_{\mathbf v}(z)^k},
\]
where $f_{\mathbf v}(z) = mz + n$ is a family of linear polynomials for $\mathbf v = (m,n) \in \hat{\mathbb{I}} \subset \mathbb I = (\Z_{> 0} \times \Z) \cup (\{ 0 \} \times \Z_{>0})$, where $\hat{\mathbb{I}}$ is a nonempty index set chosen such that $c_{\mathbf v} \neq 0$ for $\mathbf v \in \hat{\mathbb{I}}$. Write $\mathcal Z_k$ for the zero set of $F_k$ in $D$. Then  
$$ \LH_{\substack{k \rightarrow + \infty\\ k \in 2\Z_{\geq 2}}} \mathcal Z_k =  Y \defeq \bigcup_{{\mathbf{v} \neq \mathbf{w}} \in \hat{\mathbb{I}}} Y_{\mathbf{v}, \mathbf{w}} \mbox{ with } Y_{\mathbf{v}, \mathbf{w}} = \{ z \in D \colon |f_{\mathbf v}(z)| =|f_{\mathbf w}(z)| \leq  |f_{\mathbf v'}(z)| \mbox{ for all }\mathbf v' \in \hat{\mathbb{I}}    \}.$$
\end{proposition}

\begin{remark} We see that the zeros sets, that are (zero-dimensional) algebraic subsets of $\Cc$, have as limit sets semi-algebraic sets defined by inequalities. This demonstrates how the collection of algebraic subsets is not closed under Hausdorff limits, but the collection of semi-algebraic sets in this case is, cf.\ also \cite{Cesar}. 
\end{remark} 

\begin{remark} \label{geodef} Notice that solutions to equations of the form $|f_{\mathbf v}(z)|=|f_{\mathbf w}(z)|$ in $\Hh$ for $\mathbf v=(m_1,n_1) \neq \mathbf \pm w =  (m_2,n_2)$ in $\hat{\mathbb{I}}$ are geodesics. More precisely, we can rewrite the defining equation as $|\gamma z|=1$ for $\gamma = \left(\begin{smallmatrix} m_1 & n_1 \\ m_2 & n_2 \end{smallmatrix} \right) \in \mathrm{M}_2(\Z)$, so the solution set is the inverse image under $\gamma$ of the unit circle. If $m_1=0$, they form half circles centered at rational points. If $m_1 \neq 0$ and the matrix $\gamma$ were not invertible over $\Q$, then $m_1n_2=m_2n_1$; if we then set $q=n_1/m_1=n_2/m_2$, the equation becomes $|m_1(z+q)|=m_2|z+q|$, implying that $m_2 = \pm m_1$ and then also $n_2= \pm n_1$, a situation that is excluded. Thus, $\gamma \in \mathrm{M}_2(\Z) \cap \mathrm{GL}_2(\Q).$ 
The solution set is a vertical line precisely if $m_1=m_2\eqqcolon m (\neq 0)$ and then given by $\Re(z) = ({n_2-n_1})/({2m})$. 
\end{remark} 

Before presenting the (somewhat lengthly) proof of Proposition \ref{sokal}, let us show how it implies Theorem \ref{mainconfig}. For this, we use the following crucial lemma (that will also be used in the proof of Proposition \ref{sokal}). The gist of the argument is that, whilst the infinitely many terms in $F_k$ might contribute geodesic segments to the limit set of the zeros, in fact, by downward imaginary concentration, only finitely many genuinely play a role to find the limit inside $\overline \F$ (even inside the compact $\overline \F_\Gamma$). 

\begin{lemma}
\label{LemmaTopoY}
Using the notation of Proposition \ref{sokal},
there exist finite subsets $I, J \subset \hat{\mathbb{I}}$ such that
\[Y = \bigcup_{{\mathbf{v} \neq \mathbf{w}} \in I}  \{ z \in D \colon |f_{\mathbf v}(z)| =|f_{\mathbf w}(z)| \leq  |f_{\mathbf v'}(z)| \mbox{ for all }\mathbf v' \in J   \};\]
in fact, we can choose $I = \{ \mathbf v \in \hat{\mathbb{I}} \colon ||\mathbf v|| \leq \sqrt{2} \, \kappa^{a,b}_\Gamma \}$ and $J = \{ \mathbf v \in \hat{\mathbb{I}} \colon ||\mathbf v|| \leq \sqrt{7}\, (\kappa^{a,b}_\Gamma)^2 \}.$ In particular, $Y$ is contained in a finite number of geodesic segments.
\end{lemma}

\begin{proof}
The result follows from the following claim. 
$$\mbox{ For all }\mathbf v_1, \mathbf v_2 \in \hat{\mathbb{I}},\mbox{ if }||\mathbf v_1|| > \sqrt{2} \sup\limits_{z \in D} |f_{\mathbf v_2}(z)|,\mbox{ then }|f_{\mathbf v_1}(z)| > |f_{\mathbf v_2}(z)|\mbox{ for all }z \in D.$$
This claim follows from  Lemma \ref{im2est} with $\kappa=1$ (so $\alpha_\kappa=1/2$), since it implies that $|f_{\mathbf v_1}(z)| \geq ||\mathbf v_1||/\sqrt{2}$ for all $z \in D$ (the lemma is about $z \in \F$ instead of $z \in D$, but since we choose $D$ arbitrary close to $\F$, this makes no difference). 

The first consequence is that $$Y = \bigcup\limits_{{\mathbf{v} \neq \mathbf{w}} \in I} Y_{\mathbf{v}, \mathbf{w}}.$$ Indeed, $\mathbf v' = (0,\kappa^{a,b}_\Gamma) \in \hat{\mathbb{I}}$ since $c_{\mathbf v'} \neq 0$ (see Proposition \ref{decaycontrol}). Noticing that $\sup_{z \in D} |f_{\mathbf v'}(z)| = \kappa^{a,b}_\Gamma,$ the claim implies that the inequality in the definition of $Y_{\mathbf{v}, \mathbf{v}'}$ is never satisfied if $||\mathbf v||>\sqrt{2} \kappa^{a,b}_\Gamma$; and symmetrically so if  $||\mathbf w||>\sqrt{2} \kappa^{a,b}_\Gamma$. 

The second consequence is that in the occurring $Y_{\mathbf{v}, \mathbf{w}}$, the inequalities $|f_{\mathbf v}(z)| \leq  |f_{\mathbf v'}(z)|$ can be restricted to $\mathbf v' \in J$. Indeed, by the choice of $D$, there exists arbitrarily small $\epsilon_1, \epsilon_2 >0$ such that $\Im (z) \leq \kappa^{a,b}_\Gamma + \epsilon_1$ and $|\Re (z)| \leq \frac{1}{2 } + \epsilon_2$. We use the shorthand $\kappa = \kappa_\Gamma^{a,b}$. Now for $\mathbf v = (m,n)$, we have $|f_{\mathbf v}(z)|^2 = m^2 |z|^2 + 2 \Re(z) mn + n^2 \leq m^2 (1/4+\kappa^2) + |mn| + n^2 +\epsilon$ for small $\epsilon$, and the right hand side of the inequality attains $g(\kappa)$ as its maximum over $(m,n) \in I$ (i.e., $m^2+n^2 \leq 2 \kappa^2$), where \begin{equation} \label{gkappa} g(\kappa) = \kappa^4 +\kappa^2 \sqrt{\kappa^4-3/2 \kappa^2 +25/16} + 5/4 \kappa^2,\end{equation} 
where we have let $\epsilon_1,\epsilon_2,\epsilon$ tend to zero.   
For $\kappa \geq 1$, we have $2 g(\kappa) \leq (9+\sqrt{17})/2 \cdot \kappa^4 \leq 7 \kappa^4.$   
Hence, by the claim, for any $\mathbf v' \in \hat{\mathbb{I}} $ with $||\mathbf v'|| > 2 g(\kappa)$, the inequality in the definition of $Y_{\mathbf{v}, \mathbf{v}'}$ is always satisfied and may be left out.

Now $Y$ is contained in 
\begin{equation} \label{defcg} \bigcup_{{\mathbf{v} \neq \mathbf{w}} \in I}  \{ z \in \Hh \colon |f_{\mathbf v}(z)| =|f_{\mathbf w}(z)| \}, \end{equation} which, by Remark \ref{geodef}, is a finite union of geodesics in $\Hh$. The (finitely many) inequalities  $|f_{\mathbf v}(z)| \leq  |f_{\mathbf v'}(z)|$ for  $\mathbf v \in I, \mathbf v' \in J$ cut out finitely many segments of these geodesics inside $D$. 
\end{proof}

\begin{remark} 
If a segment of a vertical line occurs in $Y$, say, of the line $|mz+n_1| = |mz+n_2|$, then since $(0,\kappa_\Gamma) \in J$, we in fact need to have $|mz+n_1| = |mz+n_2| \leq \kappa_
\Gamma$ and the points belong to the bounded vertical segment 
\[
\Bigg\{ \Re(z) = \frac{n_2-n_1}{2m} \mbox{ and } \Im(z) \leq \frac{1}{m} \sqrt{ \kappa_\Gamma^2 - \left(\frac{n_1+n_2}2\right)^2} \Bigg\}.\qedhere 
\]
\end{remark} 

We are now interested in 
$$ \overline Z_{\infty,\Gamma} \defeq \LH_{\substack{k \rightarrow +\infty \\ k \in 2 \Z_{\geq 2}}} \overline Z_{k,\Gamma} = \hspace*{-5mm} \bigcup_{\substack{(a,b) \in (N^{-1}\Z/\Z)^2 \\ \gamma \in N_\Gamma}} \hspace*{-5mm}  \overline Z^{a,b}_{\infty,\gamma^{-1}\Gamma\gamma}
\, \mbox{ with } \,\overline Z^{a,b}_{\infty,\Gamma} \defeq \LH_{\substack{k \rightarrow +\infty \\ k \in 2 \Z_{\geq 2}}} \overline Z^{a,b}_{k,\Gamma}.$$

\begin{proof}[Proof of Theorem \ref{mainconfig}]
By Theorem \ref{mainIm}, it follows that for sufficiently {large}
 $k$, $\overline Z^{a,b}_{k, \Gamma} \subset D$. Hence we can directly apply Proposition \ref{sokal} by noting that $\mathcal{Z}_k \cap \overline \F = \overline Z^{a,b}_{k, \Gamma}$. 
As before, set $G =\Gamma(N)\backslash \Gamma$. Looking at the definition of $E_{k,\Gamma}^{a,b}$, we see that the term $1/(mz+n)^k$ with $k$ even occurs with coefficient 
$$ c_{m,n} =  \sum_{g \in G} \exp(2 \pi i (na_g - mb_g)) + \exp(2 \pi i (-na_g + mb_g))$$
where $(a_g, b_g) = (a,b)g$. We conclude that $c_{(m,n)}$ is a Kluyver sum 
$ c_{(m,n)} = \rho_{\Gamma}^{(Na,Nb),(n,-m)}(1).$
The result then follows from Proposition \ref{sokal} in combination with Lemma \ref{LemmaTopoY} by taking into account all conjugate groups $\gamma^{-1} \Gamma \gamma$. 
\end{proof} 

It remains to prove the Proposition \ref{sokal}. 

\begin{proof}[Proof of Proposition \ref{sokal}]  Notice that all $F_k$ are analytic on $\Hh$, since their poles in $\Cc$ are located at rational points on the real line, and we choose $D$ small enough to avoid these poles. Since $D$ is contained in a bounded domain, Hausdorff convergence to $Y$ is equivalent to \begin{equation} \label{tb} \LI \mathcal Z_k = \LS \mathcal Z_k = Y. \end{equation} 
Call an index $\mathbf v \in \hat{\mathbb{I}}$ \emph{dominant at $z \in D$} if $|f_{\mathbf v}(z)| \leq |f_{\mathbf v'}(z)|$ for all $\mathbf v' \in \hat{\mathbb{I}}$ (`dominant' because this means that the term $1/f_{\mathbf v}(z)^k$ in $F_k$ is the largest at $z$). The set $Y$ describes points where at least two distinct indices are dominant. 
For $\mathbf{v} \in \hat{\mathbb{I}} $, define  
\begin{align*}
V_{\mathbf{v}} &= \{ z \in D \colon|f_{\mathbf v}(z)| < |f_{\mathbf v'}(z)| \mbox{ for all $\mathbf v' \neq \mathbf v \in \hat{\mathbb{I}}$}   \} \mbox{ and } V = \bigcup_{\mathbf{v} \in \hat{\mathbb{I}}} V_{\mathbf{v}}, 
\end{align*}
the set of points where precisely one index is dominant. Note that $D = V \sqcup Y$. Since we always have $\LI \mathcal Z_k \subseteq \LS \mathcal Z_k$, the main result \eqref{tb} follows if we show that $V$ is a subset of the complement of $\LS \mathcal Z_k$ and that $Y$ is contained in $\LI \mathcal Z_k$. The first of these statements is proven as follows. 

\begin{lemma} $V \subset (\LS \mathcal Z_k)^{c}$. \end{lemma} 

\begin{proof} Suppose $z_0 \in V_{\mathbf{w}}$. Choose $\epsilon > 0$ such that $\overline{B}_\epsilon \defeq \{ z \colon |z - z_0| \leq \epsilon \} \subset V_{\mathbf{w}}$. Let $C$ be the maximum of $|f_{\mathbf{w}}(z)|$ on  $\overline{B}_\epsilon$. It follows there exists a $\delta > 0$ such that for all $z \in \overline{B}_\epsilon$
\begin{align*}
|f_{\mathbf{w}}(z)| \leq \delta^{-1}\mbox{ and }
\frac{|f_{\mathbf{w}}(z)|}{|f_{\mathbf{v}}(z)|} \leq 1 - \delta \mbox{ for all $\mathbf{v} \neq \mathbf{w}$ with $||\mathbf v||\leq 2 C$}.
\end{align*}
{

Put $T=\lfloor 4C^2\rfloor+1$.
Proceeding as in the proof of Lemma \ref{upperR},
we obtain
\[
\Bigg|
\sum_{\substack{\mathbf v\neq\mathbf w\\
		\|\mathbf v\|>2C}}
\frac{c_{\mathbf v}}{f_{\mathbf v}(z)^k}
\Bigg|
\leq h_\Gamma
\frac{4T(2\sqrt T+1)}{(T/4)^{k/2}}.
\]
}
By the reverse triangle inequality it follows that
{
\begin{align*}
	|F_k(z)|
	&\geq \frac{|c_{\mathbf{w}}|}{|f_{\mathbf{w}}(z)|^k}
	-\sum_{\substack{\mathbf{v}\neq\mathbf{w}\\
			\|\mathbf v\|\leq 2C}}
	\frac{|c_{\mathbf{v}}|}{|f_{\mathbf{v}}(z)|^k}
	-\Bigg|
	\sum_{\substack{\mathbf{v}\neq\mathbf{w}\\
			\|\mathbf v\|>2C}}
	\frac{c_{\mathbf{v}}}{f_{\mathbf{v}}(z)^k}
	\Bigg|\\
	&\geq \frac{|c_{\mathbf{w}}|}{|f_{\mathbf{w}}(z)|^k}
	-\sum_{\substack{\mathbf{v}\neq\mathbf{w}\\
			\|\mathbf v\|\leq 2C}}
	\frac{|c_{\mathbf{v}}|}{|f_{\mathbf{v}}(z)|^k}
	-h_\Gamma\frac{4T(2\sqrt T+1)}{(T/4)^{k/2}}\\
	&\geq \frac{|c_{\mathbf{w}}|}{|f_{\mathbf{w}}(z)|^k}
	\left(
	1-\frac{2h_\Gamma}{|c_{\mathbf{w}}|}
	\#\{\mathbf v:\|\mathbf v\|\leq 2C\}(1-\delta)^k
	-\frac{4h_\Gamma T(2\sqrt T+1)}{|c_{\mathbf{w}}|}
	\left(\frac{4C^2}{T}\right)^{k/2}
	\right)
\end{align*}}
for all $z \in \overline{B}_\epsilon$, where, in the last inequality, we have used that $|c_{\mathbf{v}}| \leq 2 h_\Gamma$ (the definition of the Eisenstein series as trace over $\Gamma/\Gamma(N)$ shows that in the rewritten form $F_k(z)$, the term $1/(mz+n)^k$ has coefficient at most twice---since exactly $\pm(m,n)$ give the same term---a sum of $h_\Gamma$ roots of unity). 
{

The expression inside the brackets tends to $1$ as
$k\to\infty$, since $T>4C^2$, and hence is positive
for sufficiently large $k$. Therefore, for these $k$,
\[
|F_k(z)| \geq |c_{\mathbf{w}}|\delta^k
\left(
1-\frac{2h_\Gamma}{|c_{\mathbf{w}}|}
\#\{\mathbf v:\|\mathbf v\|\leq 2C\}(1-\delta)^k
-\frac{4h_\Gamma T(2\sqrt T+1)}{|c_{\mathbf{w}}|}
\left(\frac{4C^2}{T}\right)^{k/2}
\right).
\]
It follows that $F_k(z)$ is non-vanishing on
$\overline{B}_\epsilon$ for all sufficiently large $k$.
We conclude that $z_0\notin\LS\mathcal Z_k$.
}
\end{proof} 

To prove that $Y \subseteq \LI \mathcal Z_k$ requires more work and is based on the following special case of a  result of Sokal (the proof of which is based on properties of normal families of harmonic functions). Set 
$$ H_k(z) \defeq \frac{1}{k} \log |F_k(z)|.$$ 
This is a well-defined function $H_k \colon D \rightarrow \R \cup \{-\infty\}$, harmonic on $D \setminus \mathcal Z_k$ and taking the `value' $- \infty$ on $\mathcal Z_k$.
\begin{namedtheorem*}{Sokal's Theorem}[Special case of {\cite[Theorem 3.2]{Sokal}}]
Suppose the functions $(|F_k|^{1/k})$ are uniformly bounded on compact subsets of $D$.
Suppose that there does not exist a neighborhood $U \ni z_0$ on which $\lim\limits_{k \rightarrow + \infty} H_k(z)$ extends to a harmonic function, or takes the value $-\infty$. 
Then $z_0 \in \LI \mathcal Z_k$.
\end{namedtheorem*}

We show that the first condition in Sokal's Theorem (that implies that $\{H_k\}$ is a normal family) is satisfied. 
 \begin{lemma}
\label{uniformlyboundedness}
    The sequence $(|F_k|^{1/k})$ is uniformly bounded on $D$.
\end{lemma}
\begin{proof}
    Write $a_{\mathbf{w}}$ for the supremum of $|f_{\mathbf{w}}(z)|^{-1}$ on $D$. 
    {
    
    Proceeding as in the proof of Lemma \ref{upperR}
    with $t=3$ and using the bound on $D$, we obtain
    \[
    |F_k(z)|
    \leq
    \sum_{\substack{\mathbf v\in\hat{\mathbb I}\\
    		\|\mathbf v\|\leq\sqrt 3}}
    |c_{\mathbf v}|a_{\mathbf v}^k
    +80h_\Gamma.
    \]
    This implies 
    \[
    |F_k(z)|^{1/k} \leq \Big(80 h_\Gamma + \sum_{\substack{\mathbf v \in \hat{\mathbb{I}}\\ ||\mathbf v||^2 \leq 3}}  |c_{\mathbf v}| \Big)^{1/k} \cdot \max\left( \{1 \} \cup \{a_{\mathbf v} \colon \mathbf v \in \hat{\mathbb{I}}  \hbox{ and } ||\mathbf v|| \leq \sqrt 3\} \right), 
    \]
    which is clearly uniformly bounded.
    }
\end{proof}

Next, we notice the following limits. 
\begin{lemma}
\label{LimietLog}
For $z \in V_{\mathbf{w}}$,  $H_k(z)$ converges to $- \log |f_\mathbf{w}(z)|$ uniformly in compact subsets of $V_{\mathbf{w}}$. 
\end{lemma}
\begin{proof}
    Let $V' \subset V_{\mathbf{w}}$ be a compact set and let $C$ be the maximum of $|f_{\mathbf{w}}(z)|$ on this set. Let $z_0 \in V'$ and write
    \begin{equation} \label{middle} 
    \frac{f_{\mathbf{w}}(z_0)^k}{c_{\mathbf{w}}} \cdot F_k(z_0) = 1 + \sum_{\substack{\mathbf{v} \neq \mathbf{w} \\ ||\mathbf v|| \leq 2 C}} \frac{c_{\mathbf{v}}}{c_{\mathbf{w}}} \Big(\frac{f_{\mathbf{w}}}{f_{\mathbf{v}}}(z_0)\Big)^k + \sum_{\substack{\mathbf{v} \neq \mathbf{w} \\ ||\mathbf v|| > 2 C}}\frac{c_{\mathbf{v}}}{c_{\mathbf{w}}} \Big(\frac{f_{\mathbf{w}}}{f_{\mathbf{v}}}(z_0)\Big)^k.
    \end{equation} 
    
 {
 
 Put $T=\lfloor 4C^2\rfloor+1$.
 Proceeding as in the proof of Lemma \ref{upperR},
 we obtain
 \[
 \Bigg|
 \sum_{\substack{\mathbf v\neq\mathbf w\\
 		\|\mathbf v\|>2C}}
 \frac{c_{\mathbf v}}{c_{\mathbf w}}
 \left(\frac{f_{\mathbf w}(z_0)}
 {f_{\mathbf v}(z_0)}\right)^k
 \Bigg|
 \leq
 \frac{4h_\Gamma T(2\sqrt T+1)}{|c_{\mathbf w}|}
 \left(\frac{4C^2}{T}\right)^{k/2}.
 \]
  }
The latter expression clearly converges to $0$ as $k \to \infty$. Finally, since $z_0 \in V'$ the (finite) middle sum on the right hand side of Equation \ref{middle} converges uniformly to $0$ as $k \to \infty$: 
 for each $\mathbf{v}$, the ratios  $f_{\mathbf{w}}/f_{\mathbf{v}}$ are uniformly bounded by a quantity smaller than $1$ on $V'$.  \end{proof}

We can now finish the proof that $Y \subset \LI \mathcal Z_k$. Let $z_0 \in Y$ and let $U$ be an arbitrary connected open neighbourhood of $z_0$. Since $Y$ is a finite union of geodesic segments, $V$ is dense in $D$, so $U \cap V$ is dense (open) in $U$. Hence there is at least one index $\mathbf v$ with $U \cap V_{\mathbf v} \neq \emptyset$. If this were the unique such index, then $U \cap V_{\mathbf v} = U \cap V$ is dense in $U$. Since $z_0 \in Y$, there is an index $\mathbf w \neq \mathbf v$ with $z_0 \in Y_{\mathbf v, \mathbf w}$.  Since $U \cap V_{\mathbf v} = U \cap V$ is dense in $U$, $|f_{\mathbf w}(z)| \leq |f_{\mathbf v}(z)|$ for all $z \in U$, while $|f_{\mathbf w}(z_0)| = |f_{\mathbf v}(z_0)|$ at $z_0 \in U$. The function $g(z) \defeq f_{\mathbf w}(z)/f_{\mathbf v}(z)$ is holomorphic on $U$ and maximal at the interior point $z_0$; hence, by the maximum principle, it is constant, contradicting the fact that $\mathbf v \neq \mathbf w \in \hat{\mathbb{I}}$. We conclude that there are at least two indices $\mathbf v \neq \mathbf w$ such that $U_{\mathbf{v}} \defeq U \cap V_{\mathbf v} \neq \emptyset$ and $U_{\mathbf{w}} \defeq U \cap V_{\mathbf w} \neq \emptyset$. Lemma \ref{LimietLog} implies that
\[
\lim_{k \to \infty} H_k(z) =
\begin{cases}
- \log |f_\mathbf{v}(z)| \quad \hbox{ for $z \in U_{\mathbf{v}}$ },\\
- \log |f_\mathbf{w}(z)| \quad  \hbox{ for $z \in U_{\mathbf{w}}$ }.
\end{cases}
\]
Clearly $\lim\limits_{k \to \infty} H_k(z)$ is harmonic on $U_{\mathbf{v}} \cup U_{\mathbf{w}}$. Suppose it extends to a harmonic function $s(z)$ on $U$. This would imply $|f_{\mathbf{v}}| = |f_{\mathbf{w}}| (= \me^{-s(z)})$, hence $f_{\mathbf{v}} \equiv \omega f_{\mathbf{w}}$ for some $\omega$ with $|\omega| = 1$. This is in contradiction with our definition of $I$. Therefore we conclude from Sokal's Theorem that $z_0 \in \LI \mathcal Z_k$.
\end{proof}

\subsection{Examples and applications}

We start by using the general result to explicitly describe the limit set for principal congruence groups, as in Proposition \ref{introcinfty}. 

\begin{proof}[Proof of Proposition \ref{introcinfty}] We will fix the level $N$ and notice that $N_{\Gamma(N)}=1$; hence we leave out $\Gamma=\Gamma(N)$ and `$\gamma$' from all notations. Notice that $c_{m,n} = 2 \cos(2 \pi(an-bm))$. We need to distinguish two cases. 

\textbf{Case 1: $N$ is not divisible by $4$, or $a \neq \pm 1/4$.} Then $c_{0,1} \neq 0$, hence $\kappa^{a,b} = 1$. This means all index sets satisfy $I^{a,b} = \{(1,0),(0,1),(1,1),(1,-1)\}$. All possible combinations of equalities of absolute values of pairs of these define the following six sets: 
\begin{equation} \label{eqN4} 
\{ |z|=1 \}; \{ |z\pm 1|=1\}, \{|z|=|z \pm 1|\}, \{|z+1|=|z-1|\}. \end{equation} 
The first three sets are unit circles centered at $0,\mp 1$, the next two define lines $\Re(z)=\mp 1/2$ and the final one defines the imaginary axis. In the entire complex plane, this defines a configuration $\mathcal C$ of generalized circles, displayed in Figure \ref{Miquel} (solid and dashed lines). 
\begin{figure}[!h]
\includegraphics[width=6cm]{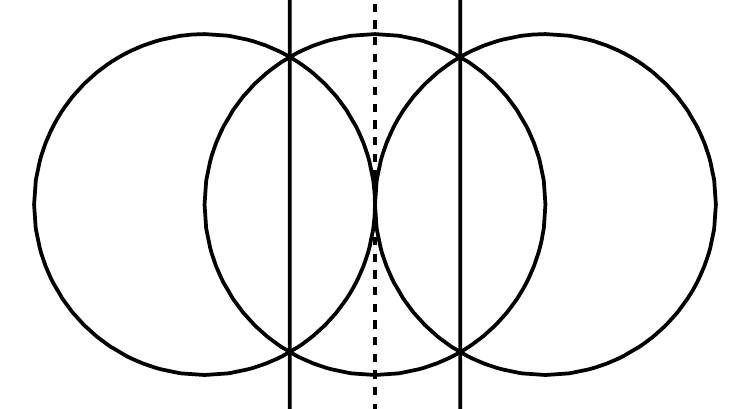}
\caption{The six circles related to the limit configuration of zeros for $\Gamma(N)$ with $N$ not divisible by $4$.} \label{Miquel} 
\end{figure}
 Recall that we only need to determine the part intersecting $\F_{\Gamma(N)}$. From the point of view of inequalities, since the three final equalities  in \eqref{eqN4} do not involve the term `$1$', their contribution to the limit set is contained in $\{|z|=|z \pm 1| \leq 1\}$ and  $\{|z+1|=|z-1|\leq 1\}$. This shows that the limit set is contained in the lower arc $\A \cup \A'$. 
To show that it is, in fact, equal to it, we only need to verify that all further inequalities are automatic for elements on this arc, i.e., that for $z \in \F \cap \{|z|=1\}$, we always have $|mz+n| \geq 1$ for all $(m,n)$. This is clear for $m=0$ or $n=0$, and for $mn \neq 0$, we have $|mz+n|^2 = m^2 |z|^2 + 2mn \Re(z)+n^2 \geq m^2-|mn|+n^2 \geq 1$. 

\textbf{Case 2: $N$ is divisible by $4$ and $a=\pm 1/4$.} In this case, $\kappa^{a,b}=2$ since $c_{0,1}=0$ but $c_{0,2} \neq 0$. We have an upper bound $m^2+n^2 \leq 8$ for $(m,n) \in I^{\pm 1/4,b}$, and a value occurs if and only if $c_{m,n} \neq 0$, where this coefficient is determined in terms of $b$ as in Table \ref{cb}. It follows that all these coefficients are non-zero, except when $b = \pm 1/4$ or $b\in\{0,1/2\}$, when we need to leave out the terms $\{|z|,|2z\pm 1|,|z\pm 2|\}$ and $\{|z\pm 1|,|2z\pm 1|\}$ from the (in)equalities defining the limit set, respectively.  

\begin{table}[t] 
\begin{tabular}{c|cccccccccc}
$(m,n)$ & $(1,0)$ & $(1,1)$ & $(1,-1)$ & $(0,2)$ & $(2,0)$ & $(2,1)$ & $(2,-1)$ & $(1,2)$ & $(1,-2)$ & $(2,2)$ \\ \hline 
$c_{m,n}$ & $c_b$ & $s_b$ & $-s_b$ & $-2$ & $c_{2b}$ & $s_{2b}$ & $-s_{2b}$ & $-c_{b}$ & $-c_{b}$ & $-c_{2b}$ \\
\end{tabular}
\caption{Possible $(m,n) \in I^{\pm 1/4,b}$ with $c_x\defeq 2 \cos(2\pi x)$ and $s_x\defeq 2 \sin(2\pi x)$} \label{cb} 
\end{table}

Concerning the inequalities, 
rather than using the exact value of $g(2)$ from \eqref{gkappa}, it is better to reason as follows: whatever $b$ is, $c_{0,2}=-2 \neq 0$, so the term $|f_{(0,2)}|=2$ occurs in all sets of (in)equalities defining the limit, either on the left hand side or the right hand side. Because of this, we first study whether it is possible to have $z \in \F_{\Gamma(N)}$ for which $|mz+n|<2$ for some integers $m,n$ at all. With $|z| \geq 1$ and $\Re(z) \geq -1/2$, we have $|mz+n|^2 \geq m^2 -|mn|+n^2$. Hence suppose $m^2 -|mn|+n^2 <4$. If $m=0$, then $n = \pm 1$, but the corresponding term is zero. If $n=0$, then $m=\pm 1$, and $(m,n)=\pm(1,0)$ occurs (unless $b = \pm 1/4$). If $mn>0$, we rewrite the requirement as $(m-n)^2+mn<4$, which has only the solutions $(m,n)=\pm(1,1)$ (that occurs unless $b\in\{0,1/2\}$), $(m,n)=\pm(1,2)$ (that occurs unless $b=\pm 1/4$) and $(m,n)=\pm(2,1)$ (that occurs unless $b \in \{0,1/2,\pm 1/4\}$). 
Similarly, if $mn<0$, we find the solutions $(m,n)=\pm(1,-1),\pm(1,-2),\pm(2,-1)$ (for $b\notin\{0,1/2\}$, $b \neq \pm 1/4$ and $b \notin\{0,1/2,\pm 1/4\}$, respectively). Now observe that, since in $\F$, $|z \pm 2| \geq |z|$ and $|2z \pm1| \geq |z|$ always hold, and if $|z \pm 2|$ or $|2z \pm 1|$ occur on the right hand side of the inequalities, then so does $|z|$, we can leave out from the inequalities the terms $|z \pm 2|$ and $|2z \pm 1|$.  
It follows that for any set of (in-)equalities in the definition of the limit set, we only need to consider inequality at most with $2, |z|, |z+1|$ or $|z-1|$ on the right. 

One can now loop over all possibilities; we have done this in Mathematica. For $b \notin \{\pm 1/4, 0, 1/2\}$ we find as limit set in $\F_{\Gamma(N)}$ the solutions to $|z|=\min\{|z+1|,|z-1|,2\}$, i.e.,  the segment of the circle of radius two centered at the origin, as well as the segments of the vertical lines $\Re(z) = \pm 1/2$ below that circle. For $b=\pm 1/4$, we find as limit set in $\F_{\Gamma(N)}$ the intersection with the sets $\{ |z\pm 1|=2 \leq |z \mp1| \}$, i.e., the part of the circle or radius $2$ centered at $\mp 1$ outside the circle of radius $2$ centered at $\pm 1$; as well as the intersection with $|z+ 1|=|z-1| \leq 2$, i.e., $\Re(z)=0 \mbox{ and } |\Im(z)| \leq \sqrt{3}$, the part of the imaginary axis up until $\ii \sqrt{3}$. Finally, the values $b=0,1/2$ contribute nothing new to the limit set. 

\textbf{Conclusion.} For $N$ not divisible by $4$, $\kappa^{a,b}=1$ for all $(a,b)$, so only the first case occurs, and we conclude that $\overline Z_{\infty, 
\Gamma(N)} = \A\cup \A'$ (as drawn using solid lines on the left of Figure \ref{anharmpic}). On the other hand, for $\Gamma(N)$ with $N$ divisible by $4$, there exist values of $(a,b)$ that give $\kappa^{a,b}=1$ and $\kappa^{a,b}=2$, so $\overline Z_{\infty, 
\Gamma(N)}$ is the union of the occurring limit sets from both cases, as indicated in the Proposition (and as drawn using solid lines on the right of Figure \ref{anharmpic}). 
\end{proof}

\begin{remark} \label{curproj} The configuration $\mathcal C$ of generalized circles in Figure \ref{Miquel} is well-known in classical projective geometry; it satisfies Miquel's---sometimes named after Jakob Steiner--- six circle theorem (`if five circles have four triple points of intersection, then the remaining intersection points are on a sixth circle',  see \cite[\S 2.8, Ex.~20]{Hist}); the result applies in our situation, where the sixth `circle' is  the imaginary axis, dashed in Figure \ref{Miquel},  containing the remaining (double) intersection points of the other lines (solid in Figure \ref{Miquel}): the origin and the point at infinity. 

 The configuration $\mathcal C'$ of solid lines in Figure \ref{Miquel} is the union of orbits of the unit circle under two groups of isometries isomorphic to $S_3$ contained in the group of (not necessarily orientation preserving) isometries $\mathrm{PGL}(2,\R)$ of $\Hh$. A general so-called \emph{anharmonic group} is of the form $\mathrm{Aut}(\widehat \Cc \setminus \{p,q,r\})$ for three distinct points $p,q,r$ on the Riemann sphere. It permutes the values of the harmonic/cross ratio of a fourth point with the three given points taken in any possible order. A different description is found in \cite{Hathaway}. All such groups are conjugate in $\mathrm{PGL}(2,\R)$ by the sharp $3$-transitivity of the action of this group on $\widehat \Cc$. and the two relevant choices for us are 
$$A_+=\mathrm{Aut}(\widehat \Cc \setminus \{0,1,\infty\}) = \langle \frac{1}{z}, 1-z \rangle  \mbox{ and } A_-=\mathrm{Aut}(\widehat \Cc \setminus \{0,-1,\infty\}) = \langle \frac{1}{z},-1-z \rangle. $$ 
The configuration $\mathcal C'$ is the union of the orbits of the unit circle under $A_+$ and $A_-$. 
\end{remark} 

\begin{remark} 
In Figure \ref{monroepic}, we superimpose the set of zeros of a finite truncation in finite weight on the limit configuration. As the weight increases, the zeros (dots) converge to the solid arcs; if we would add terms to the truncation, more and more arcs would appear, but outside $\F$. The corresponding phenomenon in condensed matter physics is the behaviour of the so-called Lee--Yang and Fisher zeros related to the (truncation at finite lattice volume of the) partition function of certain lattice models in which an external magnetic field is faded out. For a strikingly similar picture to Figure \ref{monroepic}, see the top part of Figure 3  in \cite{Monroe}, which is the case of the Baxter-Wu model (a three-spin Ising model on a triangular lattice) in the presence of a transverse magnetic field.
\begin{figure}[!h]
\includegraphics[trim={0 0 0 0},clip,width=6cm]{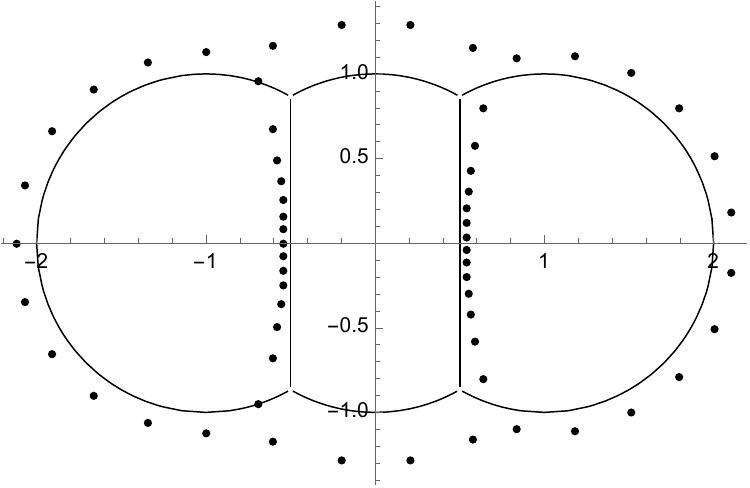} 
\caption{Solutions of the `truncation' $1 + 200/z^{20} - 10/(z + 1)^{20} + 10/(z - 1)^{20} = 0$ (dots), superimposed on the limit configuration (solid lines)} \label{monroepic} 
\end{figure}
\end{remark}

\begin{example}[\emph{Proof of Theorem \ref{mainIm} (optimality)}] \label{opt} We can now also prove the optimality of the constants $\kappa_\Gamma$ in Theorem \ref{mainIm}. 
Since  both $\mathbf v = (1,0)$ and $\mathbf w = (0,\kappa_\Gamma)$ occurs as index in $I^{a,b}_\gamma$ for some choice of Eisenstein series, we have $$Y_{\mathbf v, \mathbf w} = \{ |z|=\kappa_\Gamma \leq \min_{\mathbf v' \in J^{a,b}_\gamma} |f_{\mathbf v'}(z)| \}.$$ We claim that $\kappa_\Gamma \cdot \ii \in Y_{\mathbf v, \mathbf w}$. Indeed, whatever $\mathbf v'=(m,n)$ is, $|f_{\mathbf v'}(\kappa_\Gamma \cdot \ii)| | = \sqrt{m^2 \kappa_\Gamma^2 + n^2} \geq \kappa_\Gamma$.  Hence $\kappa_\Gamma \cdot \ii \in \lim Z_{k,\Gamma}$, which means that there is a sequence of zeros $z_i$ of Eisenstein series of weight $k$ such that $|z_i - \kappa_\Gamma \cdot \ii| \rightarrow 0$ as $k \rightarrow + \infty$. It follows that also $\Im(z_i) \rightarrow \kappa_\Gamma$. 
\end{example}

In principle, the method can be executed to compute the limit configuration for any congruence group. We illustrate this with an application to  $$\Gamma^0(N) \defeq \left (\begin{smallmatrix} 0 & 1 \\ -1 & 0 \end{smallmatrix} \right) \Gamma_0(N) \left (\begin{smallmatrix} 0 & 1 \\ -1 & 0 \end{smallmatrix} \right) = \{ \gamma \in \Gamma(1) \colon \gamma \equiv \left( \begin{smallmatrix} \ast & 0 \\ \ast & \ast \end{smallmatrix} \right) \mbox{ mod } N \}.$$ 

\begin{proposition} \label{Gamma0opt} The Hausdorff limit as $k \rightarrow + \infty$ of the \textup{(}non-$\Gamma(1)$-invariant\textup{)} set consisting of all zeros in $\overline \F$ of $\{E_{k,\Gamma^0(N)}^{a,0} \colon a \in N^{-1}\Z/\Z\}$
contains the geodesic segment $C_d \defeq \{z \in \overline \F \colon |z|=d\}$ for $d \geq 1$ precisely if $d$ divides $N/\mathrm{rad}(N)$. 
\end{proposition} 

\begin{proof}  
From the description of the limit set in Proposition \ref{sokal}, we see that $C_d$ is part of the limit set precisely in the form of $Y_{(0,\lambda d) , (\lambda,0)}$ for some integer $\lambda$, provided $c_{(0,\lambda d)} \neq 0$, $c_{(\lambda ,0)} \neq 0$ and  $|\lambda d|=|\lambda z| \leq |mz+n|$ for all $c_{(m,n)}\neq 0$. 

We study the vanishing of $c_{(m,n)}$ for the Eisenstein series $E_{k,\Gamma^0(N)}^{a,0}$ with $a=A/N$, written in the form of Proposition \ref{sokal}. Note that 
$$ G \defeq \Gamma(N) \backslash \Gamma^0(N) = \{ \left( \begin{smallmatrix} \alpha & 0 \\ \beta & 
\alpha^{-1} \end{smallmatrix} \right) \colon \alpha \in (\Z/N\Z)^*, \beta \in \Z/N\Z \}.$$

Expanding the definition as trace, we find 
$ c_{(m,n)} =2 N c_N(nA)$ in terms of Ramanujan sums. Now recall that $c_N(nA) \neq 0$ precisely if $N/(N,nA)$ is squarefree. A local calculation show that this happens precisely when $M \defeq N/\mathrm{rad}(N)$ divides $nA$, i.e., if $n$ is a multiple of $M' \defeq M/(M,A))$.  (This is similar to the considerations in Remark \ref{remkappagamma0noconj}.) 

It follows that $c_{(\lambda,0)} \neq 0$ always holds, and $c_{(0,d \lambda)} \neq 0$ holds if and only if $M' \mid d \lambda$. Since $c_{(0,M')} \neq 0$, also the inequality $\lambda d \leq M'$ should hold. We conclude that we need $\lambda d=M'$. In particular, $d$ needs to divide $M=N/\mathrm{rad}(N)$.

Conversely, suppose that $d \mid M$. By choosing $A=M/d$, we see that $(M,A)=M/d$ and so $M'=M/(M,A)=d$. Let us now check that for this choice of $A$, all further inequalities are indeed satisfied automatically.  These are $|z|=d \leq |mz+n|$ for any $(m,n) \in I$ with $M' \mid n$: if $m=0$, it follows by the minimality property of $d$, and if $m \neq 0$, it is equivalent to $|z|^2 \leq m^2 |z|^2 + n^2 + 2\Re(z)mn$, which is equivalent to $(m^2-1)d^2+n^2+2\Re(z)mn \geq 0$, and with $|\Re(z)| \leq 1/2$, follows from $$(m^2-1)d^2+n^2-|mn| > (|m|d+|n|)^2-d^2 = ((|m|-1)d+|n|)((|m|+1)d+|n|) \geq 0$$ for $|m| \geq 1$. 
This finishes the proof.  
\end{proof} 

\section{Convergence speed and angular equidistribution}

Throughout this section, we consider the case where $\Gamma=\Gamma(N)$ for $N$ odd, and, as before, Eisenstein series $E_k^{a,b}$ for $k \geq 4$ even. 

\subsection{The main term for $\F_1$} 

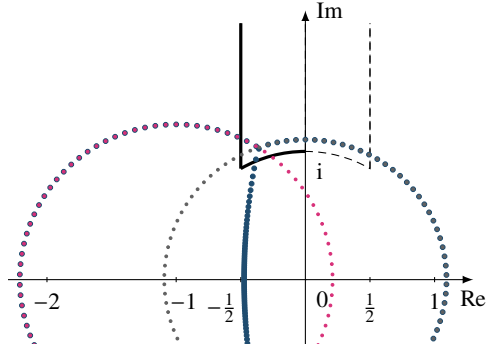
\begin{figure}[h!]
\centering
\begin{tikzpicture}[scale=1.7]
\clip (-2.3,-0.5) rectangle (2,2.2);
\draw[densely dashed] (-1/2,2)--(120:1) arc (120:60:1)--(1/2,2);   
\draw[very thick]  (120:1) arc (120:90:1);
\draw[very thick]  (120:1)--(-0.5,2);
\draw[densely dashed]  (0,1)--(0,2);
\draw[densely dashed]  (60:1)--(0.5,2);
\draw[-latex] (-2.4,0) -- (1.3,0)node[below]{$\scriptstyle \Re$};
\draw[-latex] (0,-1.3) -- (0,2.1)node[right]{$\scriptstyle \Im$};
\draw (.5,.02)--(.5,-.02)node[below]{$\scriptstyle \frac{1}{2}$};
\draw (0,0)--(0,0)node[below right]{$\scriptstyle 0$};
\draw (-.5,.02)--(-.5,-.02)node[below,xshift=-0.25cm]{$\scriptstyle -\frac{1}{2}$};
\draw (0,1)--(0,1)node[below right]{$\scriptstyle \ii$};
\draw (1,.02)--(1,-.02)node[below]{$\scriptstyle 1$};
\draw (-1,.02)--(-1,-.02)node[below,xshift=0.1cm]{$\scriptstyle -1$};
\draw (-2,.02)--(-2,-.02)node[below]{$\scriptstyle -2$};
\foreach \x/\y in {
-2.2100/-0.0405,-2.2100/0.0405,-2.2046/0.1212,-2.2046/-0.1212,-2.1938/0.2014,-2.1938/-0.2014,-2.1777/-0.2807,-2.1777/0.2807,-2.1563/0.3587,-2.1563/-0.3587,-2.1298/-0.4351,-2.1298/0.4351,-2.0982/-0.5096,-2.0982/0.5096,-2.0617/0.5818,-2.0617/-0.5818,-2.0205/0.6514,-2.0205/-0.6514,-1.9747/-0.7181,-1.9747/0.7181,-1.9245/-0.7816,-1.9245/0.7816,-1.8702/0.8416,-1.8702/-0.8416,-1.8121/-0.8979,-1.8121/0.8979,-1.7503/-0.9501,-1.7503/0.9501,-1.6852/-0.9981,-1.6852/0.9981,-1.6170/-1.0416,-1.6170/1.0416,-1.5460/-1.0805,-1.5460/1.0805,-1.4726/-1.1146,-1.4726/1.1146,-1.3971/-1.1437,-1.3971/1.1437,-1.3198/-1.1676,-1.3198/1.1676,-1.2411/-1.1864,-1.2411/1.1864,-1.1614/-1.1999,-1.1614/1.1999,-1.0809/-1.2080,-1.0809/1.2080,-1.000/-1.2107,-1.000/1.2107,-0.9191/-1.2080,-0.9191/1.2080,-0.8386/-1.1999,-0.8386/1.1999,-0.75885/-1.18639,-0.75885/1.18639,-0.68015/-1.16764,-0.68015/1.16764,-0.60288/-1.14367,-0.60288/1.14367,-0.52737/-1.11459,-0.52737/1.11459,-0.47428/0.00833,-0.47428/-0.00833,-0.47423/-0.02502,-0.47423/0.02502,-0.47411/0.04176,-0.47411/-0.04176,-0.47394/0.05860,-0.47394/-0.05860,-0.47370/-0.07556,-0.47370/0.07556,-0.47341/-0.09270,-0.47341/0.09270,-0.47304/0.11005,-0.47304/-0.11005,-0.47261/-0.12766,-0.47261/0.12766,-0.47211/0.14557,-0.47211/-0.14557,-0.47152/0.16383,-0.47152/-0.16383,-0.47086/-0.18249,-0.47086/0.18249,-0.47010/-0.20161,-0.47010/0.20161,-0.46924/0.22125,-0.46924/-0.22125,-0.46828/0.24147,-0.46828/-0.24147,-0.46719/0.26235,-0.46719/-0.26235,-0.46597/-0.28398,-0.46597/0.28398,-0.46460/0.30643,-0.46460/-0.30643,-0.46306/0.32982,-0.46306/-0.32982,-0.46134/-0.35425,-0.46134/0.35425,-0.45939/-0.37986,-0.45939/0.37986,-0.45720/0.40680,-0.45720/-0.40680,-0.45472/-0.43522,-0.45472/0.43522,-0.45398/-1.08055,-0.45398/1.08055,-0.45191/-0.46534,-0.45191/0.46534,-0.44872/-0.49736,-0.44872/0.49736,-0.44507/-0.53157,-0.44507/0.53157,-0.44088/0.56827,-0.44088/-0.56827,-0.43605/-0.60783,-0.43605/0.60783,-0.43044/-0.65070,-0.43044/0.65070,-0.42388/-0.69740,-0.42388/0.69740,-0.41616/-0.74860,-0.41616/0.74860,-0.40699/-0.80510,-0.40699/0.80510,-0.39600/-0.86789,-0.39600/0.86789,-0.38266/-0.93823,-0.38266/0.93823,-0.38009/-1.03943,-0.38009/1.03943,-0.35985/-1.02026,-0.35985/1.02026,-0.28852/-1.05354,-0.28852/1.05354,-0.21756/-1.07036,-0.21756/1.07036,-0.14558/-1.08250,-0.14558/1.08250,-0.07295/-1.08981,-0.07295/1.08981, 0/-1.09225, 0/1.09225,0.07295/-1.08981,0.07295/1.08981,0.14558/-1.08250,0.14558/1.08250,0.21756/-1.07036,0.21756/1.07036,0.28857/-1.05344,0.28857/1.05344,0.35828/-1.03181,0.35828/1.03181,0.42640/-1.00558,0.42640/1.00558,0.49261/0.97485,0.49261/-0.97485,0.55663/-0.93977,0.55663/0.93977,0.61815/-0.90049,0.61815/0.90049,0.67692/0.85719,0.67692/-0.85719,0.73266/-0.81007,0.73266/0.81007,0.78513/-0.75932,0.78513/0.75932,0.83410/0.70519,0.83410/-0.70519,0.87934/-0.64790,0.87934/0.64790,0.92065/-0.58772,0.92065/0.58772,0.95785/-0.52491,0.95785/0.52491,0.99077/-0.45976,0.99077/0.45976,1.01926/-0.39256,1.01926/0.39256,1.04321/-0.32361,1.04321/0.32361,1.06249/0.25320,1.06249/-0.25320,1.07703/-0.18167,1.07703/0.18167,1.08676/-0.10933,1.08676/0.10933,1.09164/-0.03650,1.09164/0.03650
}
{
\draw [myblue] plot [only marks, mark=*,mark size=0.52] coordinates {(\x,\y)};
}
\foreach \x/\y in {
-1.09164/-0.03650,1.09164/0.03650,0/-1.0922,0/1.0922,-1.08676/-0.10933,1.08676/0.10933,-1.07703/-0.18167,1.07703/0.18167,-1.06249/-0.25320,1.06249/0.25320,-1.04321/-0.32361,1.04321/0.32361,-1.01926/-0.39256,1.01926/0.39256,-0.99077/-0.45976,0.99077/0.45976,-0.95785/-0.52491,0.95785/0.52491,-0.92065/-0.58772,0.92065/0.58772,-0.87934/-0.64790,0.87934/0.64790,-0.83410/-0.70519,0.83410/0.70519,-0.78513/-0.75932,0.78513/0.75932,-0.73266/-0.81007,0.73266/0.81007,-0.67692/-0.85719,0.67692/0.85719,-0.61815/-0.90049,0.61815/0.90049,-0.55663/-0.93977,0.55663/0.93977,-0.49261/-0.97485,0.49261/0.97485,-0.42640/-1.00558,0.42640/1.00558,-0.35828/-1.03181,0.35828/1.03181,-0.28857/-1.05344,0.28857/1.05344,-0.21756/-1.07036,0.21756/1.07036,-0.14558/-1.08250,0.14558/1.08250,-0.07295/-1.08981,0.07295/1.08981,0.07295/-1.08981,-0.07295/1.08981,0.14558/-1.08250,-0.14558/1.08250,0.21756/-1.07036,-0.21756/1.07036,0.28857/-1.05344,-0.28857/1.05344,0.35828/-1.03181,-0.35828/1.03181,0.42640/-1.00558,-0.42640/1.00558,0.49261/-0.97485,-0.49261/0.97485,0.55663/-0.93977,-0.55663/0.93977,0.61815/-0.90049,-0.61815/0.90049,0.67692/-0.85719,-0.67692/0.85719,0.73266/-0.81007,-0.73266/0.81007,0.78513/-0.75932,-0.78513/0.75932,0.83410/-0.70519,-0.83410/0.70519,0.87934/-0.64790,-0.87934/0.64790,0.92065/-0.58772,-0.92065/0.58772,0.95785/-0.52491,-0.95785/0.52491,0.99077/-0.45976,-0.99077/0.45976,1.01926/-0.39256,-1.01926/0.39256,1.04321/-0.32361,-1.04321/0.32361,1.06249/-0.25320,-1.06249/0.25320,1.07703/-0.18167,-1.07703/0.18167,1.08676/-0.10933,-1.08676/0.10933,1.09164/-0.03650,-1.09164/0.03650
}
{
\draw [gray4] plot [only marks, mark=*,mark size=0.3] coordinates {(\x,\y)};
}
\foreach \x/\y in {
-1.0000/-1.2107,-1.0000/1.2107,-1.0809/1.2080,-0.9191/-1.2080,-1.0809/-1.2080,-0.9191/1.2080,-1.1614/1.1999,-0.8386/-1.1999,-1.1614/-1.1999,-0.8386/1.1999,-1.2411/1.1864,-0.75885/-1.18639,-1.2411/-1.1864,-0.75885/1.18639,-1.3198/1.1676,-0.68015/-1.16764,-1.3198/-1.1676,-0.68015/1.16764,-1.3971/1.1437,-0.60288/-1.14367,-1.3971/-1.1437,-0.60288/1.14367,-1.4726/1.1146,-0.52737/-1.11459,-1.4726/-1.1146,-0.52737/1.11459,-1.5460/1.0805,-0.45398/-1.08053,-1.5460/-1.0805,-0.45398/1.08053,-1.6170/1.0416,-0.38303/-1.04165,-1.6170/-1.0416,-0.38303/1.04165,-1.6852/0.9981,-0.31483/-0.99811,-1.6852/-0.9981,-0.31483/0.99811,-1.7503/0.9501,-0.24970/-0.95012,-1.7503/-0.9501,-0.24970/0.95012,-1.8121/0.8979,-0.18791/-0.89788,-1.8121/-0.8979,-0.18791/0.89788,-1.8702/0.8416,-0.12975/-0.84164,-1.8702/-0.8416,-0.12975/0.84164,-1.9245/0.7816,-0.07548/-0.78163,-1.9245/-0.7816,-0.07548/0.78163,-1.9747/0.7181,-0.02534/-0.71814,-1.9747/-0.7181,-0.02534/0.71814,-2.0205/0.6514,0.02045/-0.65143,-2.0205/-0.6514,0.02045/0.65143,-2.0617/0.5818,0.06168/-0.58182,-2.0617/-0.5818,0.06168/0.58182,-2.0982/0.5096,0.09817/-0.50961,-2.0982/-0.5096,0.09817/0.50961,-2.1298/0.4351,0.12976/-0.43512,-2.1298/-0.4351,0.12976/0.43512,-2.1563/0.3587,0.15630/-0.35869,-2.1563/-0.3587,0.15630/0.35869,-2.1777/0.2807,0.17767/-0.28065,-2.1777/-0.2807,0.17767/0.28065,-2.1938/0.2014,0.19379/-0.20137,-2.1938/-0.2014,0.19379/0.20137,-2.2046/0.1212,0.20457/-0.12118,-2.2046/-0.1212,0.20457/0.12118,-2.2100/0.0405,0.20998/-0.04045,-2.2100/-0.0405,0.20998/0.04045
}
{
\draw [mypink] plot [only marks, mark=*,mark size=0.3] coordinates {(\x,\y)};
}

\end{tikzpicture}
\caption{{For $k=94$, $N=10^8+1$, $A=\floor{N/4}$, and $B=A-10^3$, the dark blue points represent the zeros of $\cs{M}_k^{a,b}(z)$ (inside $\F_1$ we should think of this as an accurate approximation of the actual zeros of $E_k^{a,b}(z)$), the gray points represent the zeros of $f_1(z)$, and the magenta points represent the zeros of $f_2(z)$. When $z$ is close to $C_1$ and $C_2$, the zeros of $\cs{M}_k(z)$ are close to the zeros of $w_1/z^k+w_2/(z+1)^k$, which all lie on the third circle that appears in the figure.}} \label{fig: ZerosOfMainTerm}
\end{figure}

In view of the identity $E_k^{a,b}(-\ol{z})=\ol{E_k^{a,-b}(z)}$ from Lemma \ref{eisrel}, it suffices to study the region in the fundamental domain between $\Re(z)=-1/2$ and $\Re(z)=0$. We denote by $\F_1$ the aforementioned region without the imaginary axis (that itself contains at most a negligeable amount of zeros, as will follow from our theory), i.e., 
\[
\F_1\defeq \Big\{z\in \bb{H} \; : \; -\f12\le \Re(z) \lt 0 \quad\mbox{and}\quad |z|\ge 1\Big\}.
\]
Following the ideas from section \ref{section: Downward imaginary concentration}, we write 
$
E_k^{a,b}(z)=\cs{M}_k^{a,b}(z)+\cs{R}_k^{a,b}(z),
$
where $\cs{M}_k^{a,b}(z)$ is the `dominant term' and $\cs{R}_k^{a,b}(z)$ is an exponentially decreasing function on $k$. In order to choose our main term in $\F_1$, we pick up the four terms with $m^2+n^2=1$ together with the terms $(1,1)$ and $(-1,-1)$ and define
\begin{equation}\label{eq: MainTerm}
    \cs{M}_k^{a,b}(z)\defeq 2\Big(u_A+\frac{u_B}{z^k}+\frac{u_{A-B}}{(z+1)^k}\Big),
\end{equation}
where $u_j= \cos(2\pi j/N)$, $a=A/N$ and $b=B/N$. 
Note that $\cs{M}_k^{a,b}$ is essentially $M_k^{a,b,\le 1}$, but without the terms coming from $(-1,1), (1,-1), (\pm 2,0)$ and $(0,\pm 2)$, the reason being that these terms are exponentially decreasing in $\F_1$. 
Now, let $\cs{R}_k^{a,b}(z)\defeq E_k^{a,b}(z)-\cs{M}_k^{a,b}(z)$, so that
\begin{equation}
\cs{R}_k^{a,b}(z)=\frac{2u_{A+B}}{(z-1)^k}+\sum_{m^2+n^2 \ge 4}\frac{\me(an-bm)}{(mz+n)^k}
  \le {}  42\cdot 2^{-\f{k}{2}}, \label{eq: ExpDecreasing}
\end{equation} decreases exponentially in $k$ for $z \in \F_1$, by a method similar to the one used in Subsection \ref{secrem}. 
Let $Z_1^{a,b}$ denote the zeros of $E_k^{a,b}(z)$ in $\F_1$, and let $\tilde{Z}^{a,b}$ denote the zeros of $E_k^{a,b}$ in $\F$ lying on the imaginary axis, so that 
\begin{equation}\label{eq: IdentityForZk}
\#Z_{k}^{a,b}=\#Z_1^{a,b}+\#Z_1^{a,-b}+\#\tilde{Z}^{a,b}.
\end{equation}
By the valence formula,
\[
\sum_{0\le A,B\lt N}\#Z_{k}^{A/N,B/N}=\f{N^2k}{12}
\]
For each $a,b$, our objective is to find a lower bound for $\#Z_k^{a,b}$, and then use the valence formula to conclude that asymptotically, we have the correct count for $\#Z_k^{a,b}$. In order to do this, we start by noticing that since $N$ is odd, $u_A\neq 0$ for all $A$, so that $u_A^{-1}$ always exists. In order to simplify notation, for any $a,b$, we let $w_1\defeq u_Bu_A^{-1}$, and $w_2\defeq u_{A-B}u_A^{-1}$. Then, by \eqref{eq: MainTerm}, $\cs{M}_k^{a,b}(z)=0$ if and only if 
\begin{equation}\label{eq: w1,w2}
1+\Big(\frac{w_1^{1/k}}{z}\Big)^k+\Big(\frac{w_2^{1/k}}{z+1}\Big)^k=0.
\end{equation}
Consider the circles $C_1\colon |z|=|w_1|^{1/k}$, and $C_2\colon |z+1|=|w_2|^{1/k}$. Intuitively, by \eqref{eq: w1,w2}, if $z$ is away from $C_2$, then the zeros of $\cs{M}_k(z)$ are close to the zeros of 
\[
f_1(z)\defeq 1+\f{w_1}{z^k},
\]
which all lie on the circle $C_1$, and similarly, if $z$ is away from $C_2$, then the zeros of $\cs{M}_k(z)$ are close to the zeros of 
\[
f_2(z)\defeq 1+\f{w_2}{(z+1)^k},
\]
which all lie on $C_2$ (see Figure \ref{fig: ZerosOfMainTerm}; the three circles in the background are converging to the anharmonic orbit as in Remark \ref{curproj}; cf.\ also the discussion around Figure \ref{monroepic}).

The main idea of our approach is: for $|w_1|\lt 1$, since $C_1$ lies outside the fundamental domain, we don't expect to have any zeros for the Eisenstein series inside $\F_1$. In the case that $|w_1|=1$ (which only happens if $a=\pm b$ because $N$ is odd) we will use a change of sign argument similar to the one from Rankin and Swinnerton-Dyer \cite{RSD} to show that most of the zeros lie exactly on the unit circle. In the remaining case $|w_1|\gt 1$, we will use a Rouché type argument to approximate the zeros of $E_k^{a,b}$ by the zeros of $f_1(z)$. This suggests that we consider two separate cases.
\subsection{Case 1: $\boldsymbol{|w_1|\gt 1}$}\label{subsection: a different from plus or minus b}
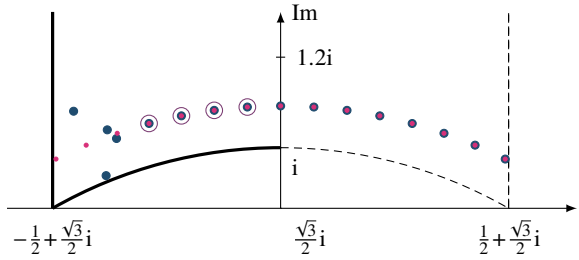
\begin{figure}[h!]
\centering
\begin{tikzpicture}[scale=6]

\draw[densely dashed] (90:1)--(90:1) arc (90:60:1)--(1/2,1.3);   
\draw[very thick]  (120:1) arc (120:90:1);
\draw[very thick]  (120:1)--(-0.5,1.3);
\draw[densely dashed]  (60:1)--(0.5,1.3);
\draw[-latex] (-0.6,{sqrt(3)/2}) -- (0.65,{sqrt(3)/2})node[below]{$\mbox{ }$};
\draw[-latex] (0,{sqrt(3)/2}) -- (0,1.3)node[right]{$\scriptstyle \Im$};
\draw (.5,{sqrt(3)/2})--(.5,{sqrt(3)/2})node[below]{$\scriptstyle \frac{1}{2}+\frac{\sqrt{3}}{2}\ii$};
\draw (0,{sqrt(3)/2})--(0,{sqrt(3)/2})node[below right]{$\scriptstyle \frac{\sqrt{3}}{2}\ii$};
\draw (-.5,{sqrt(3)/2})--(-.5,{sqrt(3)/2})node[below,]{$\scriptstyle -\frac{1}{2}+\frac{\sqrt{3}}{2}\ii$};
\draw (0,1)--(0,1)node[below right]{$\scriptstyle \ii$};
\draw (-0.01,1.2)--(0.01,1.2)node[right]{$\scriptstyle 1.2\ii$};
\foreach \x/\y in {
-0.45398/1.08055,-0.382662/0.938232,-0.380089/1.03943,-0.359849/1.02026,-0.288515/1.05354,-0.217561/1.07036,-0.145582/1.0825,-0.0729541/1.08981,0/1.09225,0.0729541/1.08981,0.145582/1.0825,0.21756/1.07036,0.288567/1.05344,0.358284/1.03181,0.426402/1.00558,0.492615/0.974851
}
{
\draw [myblue] plot [only marks, mark=*,mark size=0.25] coordinates {(\x,\y)};
}
\foreach \x/\y in {
-0.492615/0.974851,-0.426402/1.00558,-0.358284/1.03181,-0.288567/1.05344,-0.21756/1.07036,-0.145582/1.0825,-0.0729541/1.08981,0./1.09225,0.0729541/1.08981,0.145582/1.0825,0.21756/1.07036,0.288567/1.05344,0.358284/1.03181,0.426402/1.00558,0.492615/0.974851
}
{
\draw [mypink] plot [only marks, mark=*,mark size=0.125] coordinates {(\x,\y)};
}
\foreach \x/\y in {
-0.288515/1.05354,-0.217561/1.07036,-0.145582/1.0825,-0.0729541/1.08981}
{
\draw [mypurple] (\x,\y) circle (0.0182521); 
}
\end{tikzpicture}
\caption{For $k=94$, $N=10^8+1$, $A=\floor{N/4}$, and $B=A-10^3$, the dark blue points represent the zeros of $\cs{M}_k^{a,b}(z)$ (inside $\F_1$ we should think of this as an accurate approximation of the actual zeros of $E_k^{a,b}(z)$), and the magenta dots represent the zeros of $f_1(z)$. Around each zero of $f_1(z)$ with $\pi/2\lt \Arg(z)\lt t_{w_1,w_2}$ we open a small purple contour, in which we will show that $|f_1(z)|$ is dominant.} \label{fig: IllustrationRouche}
\end{figure}

First, we study the zeros of $f_1(z)$ on the circular sector of $C_1\colon z=|w_1|^{1/k}\me^{\ii t}$ with $\pi/2\le t\le t_{w_1,w_2}$, where 
\begin{equation}\label{eq: Definition of tw1w2}
t_{w_1,w_2}\defeq \arccos\Big(\frac{|w_2|^{2/k}-|w_1|^{2/k}-1}{2|w_1|^{1/k}}\Big)
\end{equation}
is the angle such that $|w_1|^{1/k}\me^{\ii t_{w_1,w2}}$ is the point of intersection inside $\F_1$ between $C_1$ and $C_2$; we show that the amount of zeros in said region is asymptotic (as $k\to\infty$) to $k/12$, and finally we conclude that the number of zeros of $E_k^{a,b}$ is also asymptotic to $k/12$ by considering a small closed contour around each zero of $f_1(z)$, and then showing that 
\begin{equation}\label{eq: fw1 dominant}
|f_1(z)|\gt \Big| \f{w_2}{(z+1)^k}+u_A^{-1}\cs{R}_k^{a,b}(z)\Big| 
\end{equation}
in that contour, and so, by Rouché's theorem, the number of zeros of $E_k^{a,b}(z)$ and $f_1(z)$ inside the contour is the same (see Figure \ref{fig: IllustrationRouche} for an illustration of this). 

With these ideas in mind, we start with the following lemma:
\begin{lemma}\label{lem: Zeros of fw1}
    Let $\tilde{C}_1$ denote the circular sector of $C_1\colon z=|w_1|^{1/k}\me^{\ii t}$ lying between $t=\pi/2$ and $t= t_{w_1,w_2}$. Denote by $Z_{\tilde{C}_1}$ denote the number of zeros of $f_{1}$ in $\tilde{C}_1$. Then, 
    \[
    Z_{\tilde{C}_1}= \Big\lfloor \frac{k}{2\pi}t_{w_1,w_2}-\frac{\delta_1}{2} \Big\rfloor -\Big\lceil\frac{k}{4}-\frac{\delta_1}{2}\Big\rceil+1,
    \]
    where 
    \[
    \delta_1\defeq \begin{cases}
  1, & \mbox{if } \sgn(w_1)=1 \\
  0, & \mbox{otherwise}.
\end{cases}
    \]
    In particular, $ Z_{\tilde{C}_1}\sim k/12$ as $k\to\infty$.
\end{lemma}
\begin{proof}
    Clearly all the zeros of $f_1(z)$ lie on the circle $C_1$, so letting $z=|w_1|^{1/k}\me^{\ii t}$, we see that $f_1(z)=0$ if and only if $\me^{-itk}=-\sgn(w_1)=\me^{-\ii\pi \delta_1}$, so that $t=(2m+\delta_1)\pi/k$ for some integer $m$, and since we want to count the zeros between $t=\pi/2$ and $t= t_{w_1,w_2}$, we see that $m$ must satisfy 
    \[
    \frac{k}{4}-\frac{\delta_1}{2}\le m \le \frac{k}{2\pi}t_{w_1,w_2}-\frac{\delta_1}{2}.
    \]
    Counting the number of such $m$ completes the proof for the count of $ Z_{\tilde{C}_1}$. For the asymptotic, note that $t_{w_1,w_2}\to 2\pi/3$ as $k\to\infty$. \end{proof}
Now we use a Rouché type argument: let $t_0=(2m+\delta_1)\pi/k \in[\pi/2, t_{w_1,w_2}]$ be such that $|w_1|^{1/k}\me^{\ii t_0}$ is a zero of $f_1(z)$, and consider the closed contour 
\begin{equation}\label{eq: Rouché contour gamma0}
\ga_{t_0}\colon z=|w_1|^{1/k}\me^{\ii(t_0+\f{c}{k}\me^{\ii t})}\quad (t\in[0,2\pi]),
\end{equation}
where $c$ is a positive constant to be chosen later. Heuristically, for big enough $k$, 
\begin{equation*}
  |z-|w_1|^{1/k}\me^{\ii t_0}|={}  |w_1|^{\frac{1}{k}}|\me^{\ii t_0}(\me^{\ii\frac{c}{k}\me^{\ii t}}-1)|
  \approx {}  |w_1|^{\frac{1}{k}} \Big|1+\ii\frac{c}{k}\me^{\ii t}-1\Big| 
  \approx{}  \frac{c}{k},
\end{equation*}
and so we can think of $\ga_{t_0}$ as a small `circle' of radius $c/k$ centered at $|w_1|^{1/k}\me^{\ii t_0}$ (it is not a true circle, but rather, a conformal image of such; we will make use of this later). We claim that when $c\lt \pi$, the curves do not intersect. In order to see this, it suffices to show that if $t_0=(2m_0+\delta_1)\pi/k$ and $t_0'=(2(m_0+1)+\delta_1)\pi/k$, then $\ga_{t_0}\cap \ga_{t_0'}=\emptyset$. Note that if, on the contrary,
\[
|w_1|^{\frac{1}{k}}\me^{\ii(t_0+\frac{c}{k}\me^{it_1})}=|w_1|^{\frac{1}{k}}\me^{\ii(t_0'+\frac{c}{k}\me^{it_2})}
\]
for some $t_1,t_2\in[0,2\pi]$, then 
\[
t_0'-t_0+\frac{c}{k}(\me^{it_2}-\me^{it_1})+2\pi n=\frac{2\pi}{k}+\frac{c}{k}(\me^{it_2}-\me^{it_1})+2\pi n= 0 
\]
for any $n\in\bb{Z}$. From this we would have the system of equations 
\[
\begin{cases}
  \frac{2\pi}{k}+\frac{c}{k}(\cos{t_2}-\cos{t_1})+2\pi n= 0  \\
  \sin{t_2}-\sin{t_1}=0
\end{cases}
\] 
with solution 
$
\cos{t_1}=\pm \sqrt{1-t_2^2}+{2\pi nk}/{c}+{2\pi}/{c}. 
$
Therefore, for big enough $k$, a solution can only exist when $n=0$ and $c\ge \pi$. This proves the claim, and in particular, we see that for $c=\frac{\pi}{2}$, the curves never intersect. Moreover, note that if $c\lt \log{|w_1|}$, and $z\in\ga_{t_0}$, then 
\[
|z|^k=|w_1|\me^{-c\sin{t}}\ge |w_1|\me^{-c}\gt |w_1|\me^{-\log{|w_1|}}=1
\]
In particular, this shows that the choice $c=\frac{1}{2}\min\{\pi,\log{|w_1|}\}$ guarantees that 
  \begin{enumerate}
        \item For $t_0\neq t_0'$, we have $\ga_{t_0}\cap \ga_{t_0'}=\emptyset$.
        \item For $z\in \ga_{t_0}$, we have $|z|\gt 1$, so all the closed contours lie entirely inside the fundamental domain, with at most one possible exception for the $t_0$ closest to $t_{w_1,w_2}$.
    \end{enumerate}
From now on, $c$ will always be $\frac{1}{2}\min\{\pi,\log{|w_1|}\}$ unless specified otherwise. After the above discussion, we are ready to show that \eqref{eq: fw1 dominant} holds as long as $t_0$ is sufficiently far from $t_{w_1,w_2}$:
\begin{lemma}\label{lem: Inequality for using Rouche}
For every sufficiently large $k$, there exists a non negative integer $\tilde{m}$ depending only on $N$ with the following property: if $t_0$ is such that $|w_1|^{1/k}\me^{\ii t_0}$ is a zero of $f_1$ with $|w_1|\gt 1$, then
    \[
    |f_1(z)|\gt \Big| \f{w_2}{(z+1)^k}+u_A^{-1}\cs{R}_k^{a,b}(z)\Big| 
    \]
 for $z\in \ga_{t_0}\cap\F_1$ with
 \begin{equation} \label{stareq}
        \f{\pi}{2}+\frac{2\pi}{k}\le t_0\le t_{w_1,w_2}-\frac{2\pi \tilde{m}}{k}.
\end{equation}
        Moreover, we can take 
\begin{equation} \label{starstareq}
    \tilde{m}\defeq \begin{cases}
                      0 & \mbox{if } c=\f{\pi}{2}, \\
                      \displaystyle\Bigg\lceil h\Bigg(\min_{\substack{0\le A,B\lt N \\ |u_B|\gt |u_A|}} \f{\log{|u_Bu_A^{-1}|}}{2}\Bigg)\Bigg\rceil+1 & \mbox{otherwise},
                    \end{cases}
\end{equation}
where $h\colon (0,\pi/2]\to \bb{R}$ is the function defined by 
    \begin{equation} \label{defh}
    h(x)\defeq \frac{1}{4\sqrt{3}\pi}\Big( \log\Big(\f{8}{(1-\me^{-\frac{\sqrt{2}}{2}x})^2}\Big)-4x\Big).
    \end{equation} 
    Furthermore, $\tilde{m}\ll \log{N}$ with an absolute implied constant.
\end{lemma}
\begin{proof}
Let $t_0\in[\f{\pi}{2}, t_{w_1,w_2}]$ be such that $|w_1|^{1/k}\me^{\ii t_0}$ is a zero of $f_1(z)$ (so $\me^{-it_0 k}=-\sgn(w_1)$), and assume further that $|w_1|^{1/k}\me^{\ii t_0}$ is not the zero closest to $\ii$ and also sufficiently far away from $|w_1|^{1/k}\me^{\ii t_{w_1,w_2}}$ by assuming that \eqref{stareq} holds 
         for some positive integer $\tilde{m}$ to be chosen later. This guarantees that all the closed contours lie entirely inside $\F_1$. Now, let $z\in \ga_{t_0}$. Our objective is to choose $\tilde{m}$ (not depending on $k$) so that for sufficiently large $k$, we have 
         \begin{equation}\label{eq: |fw1|^2 bigger than |w2/(z+1)^k|^2}
                 |f_1(z)|^2\gt \Big| \f{w_2}{(z+1)^k}+u_A^{-1}\cs{R}_k^{a,b}(z)\Big|^2
         \end{equation}
         Since $z=|w_1|^{1/k}\me^{\ii(t_0+c/k\me^{\ii t})}$ for some $t\in[0,2\pi]$, we have 
         \[
         f_1(z)=1+\frac{w_1}{z^k}=1+\sgn(w_1)\me^{-it_0 k+ic\me^{\ii t}}=1-\me^{ic\me^{\ii t}},
         \]
         so that 
         $
         |f_1(z)|^2=1-2\cos(c\cos{t})\me^{-c\sin{t}}+\me^{-2c\sin{t}}.
         $
    On the other hand, since $t_0+c/k\cos{t}\in (\pi/2,2\pi/3)$ for large enough $k$, $\cos(t_0+c/k\cos{t})\ge \cos(t_{w_1,w_2}-2\pi \tilde{m}/k + c/k\cos{t})$ since $\cos$ is decreasing in said interval. This shows that 
    \[
     \Big| \f{w_2}{(z+1)^k}\Big|^2\le \xi_k(t),
    \]
    where 
    \[
    \xi(t)=
\defeq |w_2|\Big(1+2|w_1|^{\f{1}{k}}\me^{-\f{c}{k}\sin{t}}\cos\Big(t_{w_1,w_2}-\f{2\pi \tilde{m}}{k} + \f{c}{k}\cos{t}\Big)+|w_1|^{\f{2}{k}}\me^{-\f{2c}{k}}\sin{t}\Big)^{-k}.
    \]
    Now, since 
$\lim\limits_{k\to\infty}\xi_k(t)=\exp({-2\sqrt{3}\tilde{m}\pi+\sqrt{3}c\cos{t}+c\sin{t}})\le \exp({-2\sqrt{3}\tilde{m}\pi+2c}),
    $
     we have 
    \[
    \xi_k(t)\lt 2\me^{-2\sqrt{3}\tilde{m}\pi+2c}
    \]
    for sufficiently large $k$, say $k\ge K_1$. 
    Using the inequality $|\al+\be|^2\le 2(|\al|^2+|\be|^2)$ and \eqref{eq: ExpDecreasing}, we see that 
    \[
    \Big| \f{w_2}{(z+1)^k}+u_A^{-1}\cs{R}_k^{a,b}(z)\Big|^2\le 2\xi_k(t)^2+2|u_A^{-1}|^2 42^2 2^{-k}\lt 4\me^{-4\sqrt{3}\tilde{m}\pi+4c}+|u_A^{-1}|^2 2^{12-k}
    \] 
    for $k\ge K_1$, and from this we deduce that if we choose $\tilde{m}$ such that the function $g\colon [0,2\pi] \to \bb{R}$ defined by 
    \begin{equation}\label{eq: EquationDefining m}
   g(t)\defeq |f_1(z)|^2-4\me^{-4\sqrt{3}\tilde{m}\pi+4c}= 1-2\cos(c\cos{t})\me^{-c\sin{t}}+\me^{-2c\sin{t}}-4\me^{-4\sqrt{3}\tilde{m}\pi+4c}
    \end{equation}
    is bounded away from zero, \eqref{eq: |fw1|^2 bigger than |w2/(z+1)^k|^2} holds for sufficiently large $k$. If $t\in[0,\pi/4]$, note that since $c\le \pi/2$, we have  
    $
    \cos(c\cos{t})\le \cos(\sqrt{2}/{2} c),
    $
    so that 
    \[
      1-2\cos(c\cos{t})\me^{-c\sin{t}}+\me^{-2c\sin{t}} \ge 1-2\cos\Big(\frac{\sqrt{2}}{2} c\Big)\me^{-c\sin{t}}+\me^{-2c\sin{t}}.
    \]
   It is  easy to see that in the interval $[0,\pi/4]$, the function on the right hand side attains it minimum value either at the endpoints $t=0, \pi/4$, or when $\me^{-c\sin{t}}=\cos\Big(\frac{\sqrt{2}}{2} c\Big)$. Therefore, 
    \begin{align*}
      1-&2\cos\Big(\frac{\sqrt{2}}{2} c\Big)\me^{-c\sin{t}}+\me^{-2c\sin{t}} \\ &\ge{} \min\Big\{2-2\cos\Big(\frac{\sqrt{2}}{2} c\Big),1-2\cos\Big(\frac{\sqrt{2}}{2} c\Big)\me^{-\f{\sqrt{2}}{2}c}+\me^{-\sqrt{2}c}, 1-\cos^2\Big(\frac{\sqrt{2}}{2} c\Big) \Big\}
      \ge{}  (1-\me^{-\frac{\sqrt{2}}{2}c})^2,
    \end{align*}
    where the last inequality follows from the fact that $c\in (0,\pi/2]$. If $t\in [\pi/4,\pi/2]$, then we use the trivial bound $\cos(c\cos{t})\le 1$ to obtain 
    \begin{align*}
      1-2\cos(c\cos{t})\me^{-c\sin{t}}+\me^{-2c\sin{t}}\ge{} & 1-2\me^{-c\sin{t}}+\me^{-2c\sin{t}} 
      ={} (1-\me^{-c\sin{t}})^2 
      \ge{}  (1-\me^{-\frac{\sqrt{2}}{2}c})^2.
    \end{align*}
    By symmetry, it is easy to see that the same bound holds in the whole interval $[0,2\pi]$, so that 
    $
    g(t)\ge (1-\me^{-\frac{\sqrt{2}}{2}c})^2-4\me^{-4\sqrt{3}\tilde{m}\pi+4c}.
    $
  Let 
    $
    \tilde{c}\defeq {(1-\me^{-\frac{\sqrt{2}}{2}c})^2}/{2}\gt 0.
    $
    Then, $g(t)\gt \tilde{c}$ (meaning $g$ is bounded away from zero) as long as we choose $\tilde{m}$ such that 
    $
    (1-\me^{-\frac{\sqrt{2}}{2}c})^2\gt 8\me^{-4\sqrt{3}\tilde{m}\pi+4c},
    $
    or equivalently, $\tilde{m} > h(c)$ with 
    the function $h\colon (0,\pi/2]\to \bb{R}$ defined as in Equation \ref{defh}, decreasing, with $h(\pi/2)\approx -0.15$. Therefore, upon recalling that 
    $
    c=\f12\min\{\pi,\log{|w_1|}\} 
    $
    and $|w_1|=|u_Bu_A^{-1}|$, we see that we may take $\tilde{m}$ as in \eqref{starstareq}. 
    
    To complete the proof, it remains to show that $\tilde{m}\ll \log{N}$. In order to do this, since $h$ is decreasing, to obtain an upper bound for $\tilde{m}$, we need a lower bound for 
    \[
  \min_{\substack{0\le A,B\lt N \\ |u_B|\gt |u_A|}} \Big|\f{u_B}{u_A}\Big|,
  \]
  and by symmetry, we may assume without loss of generality that $u_A, u_B \gt 0$, so that $0\le A,B\lt N/4$ with $B\gt A$. Then, by the mean value theorem, there is some $\al \in (A,B)$ such that 
  \[
  u_B-u_A =-\frac{2\pi}{N}\sin\Big(\frac{2\pi \al}{N}\Big) (A-B).
  \]
  Now, observe that for $0\le x\le 1/4$, we have $\sin(2\pi x)\ge 4x$, and since $\al\lt A \lt N/4$, then 
  \begin{equation}\label{eq: u_B/u_A lower bound}
  \frac{u_B}{u_A}=1+\frac{2\pi}{N}\sin\Big(\frac{2\pi \al}{N}\Big) (A-B) u_A^{-1}\ge 1+\frac{8\pi^2 \al}{N^2}(A-B).
  \end{equation}
  This leads us to consider two separate cases, one in which $\al$ could possibly be arbitrarily close to $0$ (e.g. when $B=0$), and one where $\al \ge 1$.

  \textbf{Case 1: $\boldsymbol{B=0}$.} In this case, since $u_B=1$, it is clear that the minimum occurs when $u_A$ is closest to 1, but $A\gt B=0$, so $A=1$. Then, 
  \[
h\Big(\f12 \log\Big(\f{u_0}{u_1}\Big)\Big)\le \log\Big(\frac{8}{1-\cos(\frac{2\pi}{N})^{\frac{\sqrt{2}}{4}}}\Big) \ll \log{N},
  \]
  showing that $\tilde{m}\ll \log{N}$.

  \textbf{Case 2: $\boldsymbol{B\ge 0}$.} In this case, since $\al\gt B$, and $B$ is an integer, $\al \ge 1$. Then, using \eqref{eq: u_B/u_A lower bound}, we have 
  $
  u_B/u_A\ge 1+8\pi^2/N^2,
  $
  and from this, a routine computation shows that
  \[
  h\Big(\f12 \log\Big(\f{u_B}{u_A}\Big)\Big) \le h\Big(\f12 \log\Big(1+\frac{8\pi^2}{N^2}\Big) \Big) \ll \log{N},
  \]
  again implying that $\tilde{m}\ll \log{N}$. This completes the proof.
\end{proof}
\begin{corollary}\label{cor: main result for a different than b}
  Let $\tilde{m}$ be as in Lemma \ref{lem: Inequality for using Rouche}, and let $Z_{\tilde{C_1}}$ be as in Lemma \ref{lem: Zeros of fw1}. Then, for every $a,b$ with $a\neq \pm b$ and $|w_1|\gt 1$, $E_k^{a,b}$ has at least $Z_{\tilde{C_1}}-(\tilde{m}+1)$ zeros in $\F_1$. Moreover, these zeros are all simple, and lie within Euclidean distance $\asymp 1/k$ from the unit circle.
\end{corollary}
\begin{proof}
  By Lemma \ref{lem: Inequality for using Rouche}, for all $t_0$ such that $f_1(|w_1|^{1/k}\me^{\ii t_0})=0$, we have 
  \[
    |u_A f_1(z)|\gt \Big| \f{u_A w_2}{(z+1)^k}+\cs{R}_k^{a,b}(z)\Big| 
    \]
    for all $z\in \ga_{t_0}\cap \F_1$ with \eqref{stareq}. 
     Therefore, by Rouché's theorem, $u_A f_1(z)$ and 
         $
     E_k^{a,b}(z)=u_A f_1(z)+{u_A w_2}/{(z+1)^k}+\cs{R}_k^{a,b}(z)
     $
     have the same number of zeros in the interior of each $\ga_{t_0}$. Observe that $f_1$ has exactly one simple zero inside each $\ga_{t_0}$ (because these closed contours do not overlap), so that $E_k^{a,b}$ has exactly one simple zero inside each $\ga_{t_0}$. The number of such $t_0$ is $Z_{\tilde{C_1}}-(\tilde{m}+1)$. Finally, since the diameter of the interior of $\ga_{t_0}$ is $\asymp 1/k$, by our choice of $c$, these zeros lie within Euclidean distance $\asymp 1/k$ from the unit circle.
\end{proof}

\subsection{Case 2: $\boldsymbol{a=\pm b}$}

Since the case $a=b=0$ is already known by the work of Rankin and Swinnerton-Dyer \cite{RSD}, we will assume that $a,b$ are not simultaneously 0. As $a=\pm b$, we have $u_A=u_B$. Recall from the previous subsection that $w_1= u_Bu_A^{-1}=1$ and $w_2= u_{A-B}u_A^{-1}$, so that $C_1$ is the unit circle. We will use a change of sign argument to show that most of the zeros of $E_k^{a,b}$ lie in the unit circle. Afterwards, we deal with the zeros close to $C_2$ by again using a Rouché type argument, and we shall see that our lower bound for $\#Z_k^{a,b}$ is better in this case. 
\subsubsection*{Zeros on the unit circle}
On the unit circle, $z=\me^{\ii t}$ with $\frac{\pi}{2}\le t \le \frac{2\pi}{3} $. The function $h_k^{a,b}\colon [\frac{\pi}{2},\frac{2\pi}{3}]\to \bb{C}$ defined by 
\begin{equation} \label{defht} 
h_k^{a,b}(t)\defeq \frac{u_A^{-1}}{2}\me^{ikt/2}E_k^{a,b}(\me^{\ii t})
\end{equation} 
is real valued because since for $a=b$ or $a=-b$, we have
\[
\me^{ikt/2}E_k^{a,b}(\me^{\ii t})=\sideset{}{'}\sum_{m,n \in \Z} \frac{\e(an - bm)}{(m \me^{it/2} + n\me^{-it/2})^k}=\sideset{}{'}\sum_{m,n \in \Z} \frac{\e(-an + bm)}{(m \me^{-it/2} + n\me^{it/2})^k}=\overline{\me^{ikt/2}E_k^{a,b}(\me^{\ii t})}.
\]
Therefore, a (real) zero of $h_k^{a,b}(t)$ corresponds uniquely to a zero of $E_k^{a,b}(z)$ on the unit circle. Using \eqref{eq: MainTerm}, we have 
\begin{equation}\label{eq: Definition of h(t)}
h_k^{a,b}(t)=2\cos(tk/2)+w_2(2\cos(t/2))^{-k}+ \frac{u_A^{-1}}{2}\me^{ikt/2}\cs{R}_k^{a,b}(\me^{\ii t}).
\end{equation}
We compare  the zeros of $h_k^{a,b}(t)$ with those of $2\cos(kt/2)$. The following lemma gives a lower bound for the number of zeros of $h_k^{a,b}(t)$ for $t\in [\frac{\pi}{2},\frac{2\pi}{3}]$.
\begin{lemma}\label{lem: Function h(t)}
  For every odd level $N$ and $A\in [0,N)\cap\bb{Z}$ and for any sufficiently large even integer $k$, there is a constant $\tilde{l}$ depending only on $A$ and $N$ and the congruence class of $k \mod{3}$ such that the function $h_k^{a,b}$ has at least 
\[
\Big\lfloor\frac{k}{12}\Big\rfloor-\tilde{l}-\lambda_k
\]
 zeros on $[\frac{\pi}{2},\frac{2\pi}{3}]$, where 
 \[
 \lambda_k=\begin{cases}
            1 & \mbox{if } k\equiv 0,2,4,8 \mod{12}, \\
            0 & \mbox{if } k\equiv 6,10 \mod{12}.
          \end{cases}
 \]
 Moreover, 
 \begin{equation} \label{deftildel}
 \tilde{l}=\begin{cases}
      1 & \mbox{if } k\equiv 0 \mod 3 \mbox{ and } |w_2|=2, \\
      0 & \mbox{if } k\not\equiv 0 \mod 3 \mbox{ and } |w_2|=2 \\ & \mbox{or } |w_2|\lt 2,\\
      \bigg\lceil \frac{\log|w_2|-\log{2}}{\sqrt{3}\pi}-\frac{\mu_k}{3}\bigg\rceil+ \mathds{1}_{\bb{Z}}\Big( \frac{\log|w_2|-\log{2}}{\sqrt{3}\pi}-\frac{\mu_k}{3}\Big)  & \mbox{otherwise};
    \end{cases}
 \end{equation}
 with $\mathds{1}_{\bb{Z}}$ as the characteristic function of $\bb{Z}$, and
    \begin{equation} \label{defmuk}
  \mu_k \in \{0,1,2\} \mbox{ with }\mu_k \equiv d \mod 3.  \end{equation}
\end{lemma}
\begin{proof}

  Let $\lambda_k$ be as above. Since the interval $I_k\defeq [\frac{2\pi}{k}(\floor{\frac{k}{4}}+1),\frac{2\pi}{k}\floor{\frac{k}{3}}]$ is contained in $[\frac{\pi}{2},\frac{2\pi}{3}]$, it suffices to show that $h_k^{a,b}$ changes sign at least $\floor{k/12}-\tilde{l}-\lambda_k$ times in said interval. Let $t=\frac{2\pi m}{k} \in I_k$ for some integer $m$. Since $\floor{k/3}-\floor{k/12}=\floor{k/4}+1-\lambda_k$, $m$ can be written as $m=\floor{k/3}-l$ for some integer $0\leq l \leq \floor{k/12}-\lambda_k$. Let $\tilde{l}\in [0,\floor{k/12}]$ be the smallest integer such that for every sufficiently large even integer $k$, we have
  \begin{equation}\label{eq: Condition<2}
  |w_2|\Big(2\cos\Big(\frac{\pi}{k}\Big(\Big\lfloor\frac{k}{3}\Big\rfloor-\tilde{l}\Big)\Big)\Big)^{-k}+21|u_A|^{-1} 2^{-k/2} \lt 2. 
  \end{equation}
  Intuitively, $\tilde{l}$ is the smallest integer for which $h_k^{a,b}(t)$ is well approximated by $2\cos(kt/2)$ in the interval $[\frac{2\pi}{k}(\floor{\frac{k}{4}}+1),\frac{2\pi}{k}\floor{\frac{k}{3}}-\tilde{l}]$ (see Figure \ref{fig: IllustrationOfNoProblems} for a comparison). The existence of such $\tilde{l}$ will depend on the value of $w_2$ and the weight $k$. We consider 3 separate cases: 
  
  \textbf{Case 1: $\boldsymbol{|w_2|=2}$.} 
  If $k\equiv 0\mod 3$, then
 $
  |w_2|(2\cos(\frac{\pi}{k}(\lfloor\frac{k}{3}\rfloor-0)))^{-k}=2,
  $
  so we must have $\tilde{l}\gt 0$. Now, since 
  $
  2(2\cos(\frac{\pi}{3}-\frac{\pi}{k}))^{-k}\lt 1
  $
  for all $k\ge 2$, it is clear that 
  \[
  2\Big(2\cos\Big(\frac{\pi}{3}-\frac{\pi}{k}\Big)\Big)^{-k} +21|u_A|^{-1} 2^{-k/2} \lt 2
  \]
  for sufficiently big $k$. Therefore, $\tilde{l}=1$ in this case.
  
  On the other hand, If $k\equiv 1,2\mod 3$, then a computation reveals that 
  $
  2(2\cos(\frac{\pi}{k}\lfloor\frac{k}{3}\rfloor))^{-k}\lt 1
  $
  for all $k\ge 2$, so similarly as before, for sufficiently large $k$, we have $\tilde{l}=0$.
  
  \textbf{Case 2: $\boldsymbol{|w_2|\lt 2}$.}
 Then, 
   \[
   |w_2|\Big(2\cos\Big(\frac{\pi}{k}\Big\lfloor\frac{k}{3}\Big\rfloor\Big)\Big)^{-k} +21|u_A|^{-1} 2^{-k/2}\le  |w_2|+21|u_A|^{-1} 2^{-k/2}\lt 2
   \]
   for sufficiently large $k$. Therefore, $\tilde{l}=0$. 
   
 \textbf{Case 3: $\boldsymbol{|w_2|\gt 2}$.}
 In this case, we solve \eqref{eq: Condition<2} for $\tilde{l}$, and we see that $\tilde{l}$ is the smallest integer such that 
 \begin{equation*}
 \tilde{l}\gt  \Big\lfloor\frac{k}{3}\Big\rfloor-\frac{k}{\pi}\arccos\Big(\frac{1}{2}\Big(\frac{|w_2|}{2-21|u_A|^{-1} 2^{-k/2}}\Big)^{1/k}\Big).
 \end{equation*}
Observe that for big enough $k$, the right hand side of the above equation is well defined because 
 \begin{equation}\label{eq: arccos<1}
 \frac{1}{2}\Big(\frac{|w_2|}{2-21|u_A|^{-1} 2^{-k/2}}\Big)^{1/k}\lt 1.
 \end{equation}
 This shows that for sufficiently big $k$, $\tilde{l}$ is always well defined. Furthermore, 
 \[
 \Big\lfloor\frac{k}{3}\Big\rfloor-\frac{k}{\pi}\arccos\Big(\frac{1}{2}\Big(\frac{|w_2|}{2-21|u_A|^{-1} 2^{-k/2}}\Big)^{1/k}\Big)=\frac{\log|w_2|-\log{2}}{\sqrt{3}\pi}-\frac{\mu_k}{3}+o(1),
 \]
  where $\mu_k$ is as in \eqref{defmuk}. This shows that 
  \[
  \tilde{l}=\bigg\lceil \frac{\log|w_2|-\log{2}}{\sqrt{3}\pi}-\frac{\mu_k}{3}\bigg\rceil+ \mathds{1}_{\bb{Z}}\Big( \frac{\log|w_2|-\log{2}}{\sqrt{3}\pi}-\frac{\mu_k}{3}\Big)
  \]
 The definition \eqref{eq: Condition<2} of $\tilde{l}$ and \eqref{eq: ExpDecreasing} shows that for any $m=\floor{k/3}-l$ with $\tilde{l}\le l\le \floor{k/12}-\lambda_k$,
 \begin{align*}
   \Big|w_2(2\cos(\pi m/k))^{-k} + \frac{u_A^{-1}}{2}\me^{ikt/2}\cs{R}_k^{a,b}(\me^{\ii t})\Big|  \le{} |w_2|(2\cos(\pi m/k))^{-k} +21|u_A|^{-1} 2^{-k/2} 
   \lt{}  2.
 \end{align*}
Hence, 
$
2(-1)^m-2\lt h_k^{a,b}(2\pi m/k)\lt 2(-1)^m+2,
$
so $h_k^{a,b}$ is strictly positive or negative according to whether $m$ is even or odd, and so $h_k^{a,b}$ changes sign at least $\floor{k/12}-\tilde{l}-\lambda_k$ times. 
\end{proof}
The previous lemma allows us to give a lower bound for the number of zeros of $h_k^{a,b}(t)$ in the interval $[\frac{2\pi}{k}(\floor{\frac{k}{4}}+1),\frac{2\pi}{k}(\floor{\frac{k}{3}}-\tilde{l})]$ (we also expect this count to be exact in said interval). 
\begin{example}\label{exmp: IllustrationNoProblems}
  As an illustration of the lemma, for $k=122, N=10^3+1, A=\floor{N/4}$, $B=-A$, we have $\lambda_k=1$, and a numerical computation gives $\tilde{l}=1$, so there are at least $\floor{k/12}-\tilde{l}-1=8$ zeros in the interval $[62\pi/122, 78\pi/122]\approx[1.596,2.008]$; see Figure \ref{fig: IllustrationOfNoProblems}.
  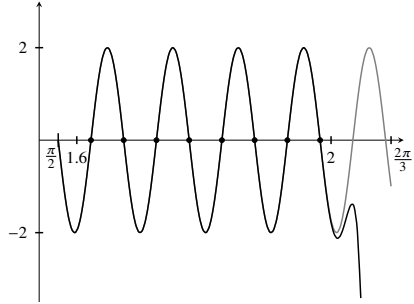
\begin{figure}[h!]
  \centering
\begin{tikzpicture}[scale=0.7]
  \begin{axis}[
    axis lines = middle,
    xmin=pi/2-0.03, xmax=2*pi/3+0.02,
     ymin=-3.5, ymax=3, 
    clip=false,
    xtick={1.6,2,2*pi/3},
    xticklabels={1.6, 2,},
    ytick={-2, 2,6},
    tick label style={font=\scriptsize},
    x tick style={black, thick},  
    y tick style={black, thick}
  ]
   \addplot[
     domain=pi/2:2*pi/3,
     samples=100,
     smooth,
     thick,
     color=black!50,
     restrict y to domain=-7:7
   ] {2*cos(61*deg(x))};
   \addplot[
     domain=pi/2:2*pi/3,
     samples=100,
     smooth,
     thick,
     black,
     restrict y to domain=-6.35:7 
   ] {2*cos(61*deg(x))-637.254* (2*cos(0.5 *deg(x)))^(-122)};

   \node[anchor=north west, font=\footnotesize] at (axis cs:2.084,0) {$\frac{2\pi}{3}$};
   \node[anchor=north west, font=\footnotesize] at (axis cs:pi/2-0.035,0) {$\frac{\pi}{2}$};
\addplot[
  only marks,
  mark=*,
  mark size=1.3pt,
  black
] coordinates {((2*31+1)*pi/122, 0) ((2*32+1)*pi/122, 0)((2*33+1)*pi/122, 0)((2*34+1)*pi/122, 0)((2*35+1)*pi/122, 0)((2*36+1)*pi/122, 0)((2*37+1)*pi/122, 0)((2*38+1)*pi/122, 0)};

\draw[black, thick] (axis cs:pi/2, -0.15) -- (axis cs:pi/2, 0.15);
  \end{axis}
\end{tikzpicture}
\caption{$h_{122}^{a,b}(t)$ in black and $2\cos(61t)$ in gray.}\label{fig: IllustrationOfNoProblems}
\end{figure}
\end{example}
Since a zero of $h_k^{a,b}(t)$ corresponds to a zero of $E_k^{a,b}(t)$ on the unit circle, we have the following immediate Corollary of Lemma \ref{lem: Function h(t)}:
\begin{corollary}\label{cor: Main theorem for a,a}
  Let $\tilde{l}$ and $\lambda_k$ be as in Lemma \ref{lem: Function h(t)}. Then, for all $a=\pm b$, and every sufficiently large $k$, $E_k^{a,b}$ has at least
  $
\lfloor\frac{k}{12}\rfloor-\tilde{l}-\lambda_k
$
 zeros on $[\frac{\pi}{2},\frac{2\pi}{3}]$ lying exactly on the unit circle.
\end{corollary}

\subsubsection*{Zeros close to $C_2$}
\begin{figure}[h!]
\centering
\begin{tikzpicture}[scale=6]

\draw[densely dashed] (90:1)--(90:1) arc (90:60:1)--(1/2,1.3);   
\draw[very thick]  (120:1) arc (120:90:1);
\draw[very thick]  (120:1)--(-0.5,1.3);
\draw[densely dashed]  (60:1)--(0.5,1.3);
\draw[-latex] (-0.6,{sqrt(3)/2}) -- (0.65,{sqrt(3)/2})node[below]{$\mbox{ }$};
\draw[-latex] (0,{sqrt(3)/2}) -- (0,1.3)node[right]{$\scriptstyle \Im$};
\draw (.5,{sqrt(3)/2})--(.5,{sqrt(3)/2})node[below]{$\scriptstyle \frac{1}{2}+\frac{\sqrt{3}}{2}\ii$};
\draw (0,{sqrt(3)/2})--(0,{sqrt(3)/2})node[below right]{$\scriptstyle \frac{\sqrt{3}}{2}\ii$};
\draw (-.5,{sqrt(3)/2})--(-.5,{sqrt(3)/2})node[below,]{$\scriptstyle -\frac{1}{2}+\frac{\sqrt{3}}{2}\ii$};
\draw (0,1)--(0,1)node[below right]{$\scriptstyle \ii$};
\draw (-0.01,1.2)--(0.01,1.2)node[right]{$\scriptstyle 1.2\ii$};
\foreach \x/\y in {
-0.4788/1.08228,-0.410528/1.04667,-0.344651/1.00675,-0.288359/0.957522,-0.255824/0.966723,-0.191122/0.981566,-0.127878/0.99179,-0.0640702/0.997945,0/1,0.0640702/0.997945,0.127877/0.99179,0.191159/0.981559,0.253655/0.967295,0.315108/0.949056,0.375267/0.926917,0.433884/0.900969,0.490718/0.871319}
{
\draw [myblue] plot [only marks, mark=*,mark size=0.25] coordinates {(\x,\y)};
}
\foreach \x/\y in {
-0.4788/1.08228,-0.410528/1.04667,-0.344679/1.00675,-0.281523/0.962694
}
{
\draw [mypink] plot [only marks, mark=*,mark size=0.125] coordinates {(\x,\y)};
}
\foreach \x/\y in {
-0.410528/1.04667,-0.344679/1.00675}
{
\draw [mypurple] (\x,\y) circle (0.019254); 
}
\end{tikzpicture}
 \caption{For $k=98$, $N=10^8+1$, $A=\floor{N/4}$, and $B=A$, the dark blue points represent the zeros of $\cs{M}_k^{a,b}(z)$, and the magenta dots represent the zeros of $f_2(z)$. Around each zero of $f_2(z)$ with $t_{w_1,w_2}\lt \Arg(z)\lt t'_{w_2}$ we open a small purple contour, in which we will show that $|f_2(z)|$ is dominant.}\label{fig: IllustrationRoucheC2}
\end{figure}
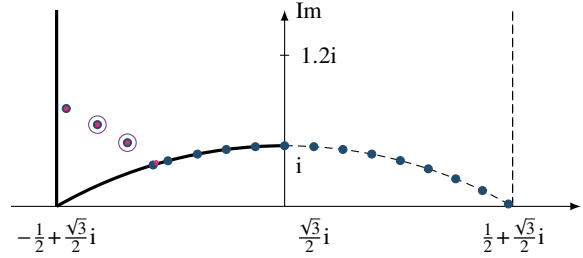

A completely analogous argument as the one used in the proof of Lemma $\ref{lem: Zeros of fw1}$ allows us to find the number of zeros of $f_2$ inside $\F_1$:
\begin{lemma}\label{lem: Zeros of fw2}
    Let $\tilde{C}_2$ denote the circular sector of $C_2\colon z+1=|w_2|^{1/k}\me^{\ii t}$ lying between $t= t_{w_1,w_2}$ and $t=t'_{w_2}$, where 
    $
    t'_{w_2}=\arccos(1/(2|w_2|^{1/k}))
    $
    is the angle such that $|w_2|^{1/k}\me^{\ii t'_{w2}}-1$ is the point of intersection inside $\F_1$ between $C_2$ and the line $\Re(z)=-1/2$. Now, denote by $Z_{\tilde{C}_2}$ denote the number of zeros of $f_{2}$ in $\tilde{C}_2$. Then, 
    \[
    Z_{\tilde{C}_2}= \Big\lfloor \frac{k}{2\pi}t'_{w_2}-\frac{\delta_2}{2} \Big\rfloor -\Big\lceil \frac{k}{2\pi}t_{w_1,w_2}-\frac{\delta_2}{2}\Big\rceil+1,
    \]
    where 
    \[
    \delta_2\defeq \begin{cases}
  1 & \mbox{if } \sgn(w_2)=1, \\
  0 & \mbox{otherwise}.
\end{cases}
    \]
\end{lemma}
Let $z_0=|w_2|^{1/k} \me^{\ii t_0}-1$ be a zero of $f_2(z)$ in $\tilde{C}_2$, so $t_0=(2m_0+\delta_{2})\pi/k$ for some integer $m_0$ such that $t_{w_1,w_2}\le t_0 \le t'_{w_2}$. For each $t_0$, consider the contour 
\begin{equation}\label{eq: CircleAroundTheAngleC2}
  \ga'_{t_0}\colon z=|w_2|^{\frac{1}{k}}\me^{\ii(t_0+\frac{c'}{k}\me^{\ii t})}-1, \quad (t\in[0,2\pi])
\end{equation}
where $c'$ is a constant to be chosen later. By a similar argument as the one from Subsection \ref{subsection: a different from plus or minus b}, it is easy to see that the choice $c'=\pi/2$ guarantees that the curves do not intersect. From this, it follows that if $t_0\neq t_{w_1,w_2}, t'_{w_2}$, then $\ga'_{t_0}$ lies entirely inside $\F_1$. Our next objective is to show that for each such $t_0$, we have
\[
  |f_2(z)|\gt \Big|\frac{1}{z^k}+\frac{u_A^{-1}}{2}\cs{R}_k^{a,b}(z)\Big|
  \] 
for all $z\in \ga'_{t_0}$, and then use Rouché's theorem to conclude that the number of zeros of $E_k^{a,b}$ and $f_2$ is the same in the interior of each contour (see Figure \ref{fig: IllustrationRoucheC2}).

We are now ready to study the zeros of $E_k^{a,b}$ that are close to $C_2$:
\begin{lemma}\label{lem: ZerosCloseToC2}
  Let $Z_{\tilde{C}_2}$ be as in Lemma \ref{lem: Zeros of fw2}. If $\#Z_{\tilde{C}_2}\ge 3$, then for big enough $k$, $E_k^{a,a}$ and $E_k^{a,-a}$ have at least $\#Z_{\tilde{C}_2}-2$ zeros on $\F_1$ lying outside the unit circle. Moreover, these zeros lie at Euclidean distance of order $O(1/k)$ from the circle $C_2$.
\end{lemma}
\begin{proof}
  For each $t_0=(2m_0+1)\pi/k\in (t_{w_1,w_2},t'_{w_2})$, let $\ga'_{t_0}$ be as in \eqref{eq: CircleAroundTheAngleC2}. Since $t_0\neq t_{w_1,w_2}, t'_{w_2}$, we may assume that 
  \[
  t_{w_1,w_2}+\frac{2\pi}{k}\le t_0 \le t'_{w_2}-\frac{2\pi}{k}.
  \]
  Let $z=|w_2|^{\frac{1}{k}}\me^{\ii(t_0+\frac{\pi}{2k}\me^{\ii t})}-1\in \ga'_{t_0}$. Then, 
  \[
  f_2(z)= 1+\frac{w_2}{(z+1)^k}= 1+\sgn(w_2)\me^{-ic\me^{\ii t}},
  \]
  so that 
  \[
  |f_2(z)|=\sqrt{1-2\sgn(w_2)\me^{\frac{\pi}{2}\sin{t}}\cos\Big(\frac{\pi}{2}\cos{t}\Big)+\me^{\pi\sin{t}}},
  \]
  and, since $\cos\big(\frac{\pi}{2}\cos{t}\big)\ge 0$, 
    \[
  |f_2(z)|\ge \sqrt{1-2\me^{\frac{\pi}{2}\sin{t}}\cos\Big(\frac{\pi}{2}\cos{t}\Big)+\me^{\pi\sin{t}}}.
  \]
  On the other hand, we have
  \[
  \Big|\frac{1}{z^k}\Big|=\Big((|w_2|\me^{-\frac{\pi}{2}\sin{t}})^{\frac{2}{k}}-2(|w_2|\me^{-\frac{\pi}{2}\sin{t}})^{\frac{1}{k}}\cos\Big(t_0+\frac{\pi}{2k}\cos{t}\Big)+1\Big)^{-\frac{k}{2}}.
  \] 
  Since $t_0\ge t_{w_1,w_2}+\frac{2\pi}{k}$, then by letting $R\defeq |w_2|\me^{-\frac{\pi}{2}\sin{t}}$, we see that 
  \[
  \Big|\frac{1}{z^k}\Big|\le \Big(R^{\frac{2}{k}}-2R^{\frac{1}{k}}\cos\Big(t_{w_1,w_2}+\frac{2\pi}{k}+\frac{\pi}{2k}\cos{t}\Big)+1\Big)^{-\frac{k}{2}}.
  \]
  Consider the function $\beta\colon \bb{Z}_{\ge 4}\times [0,2\pi]\to \bb{R}$ defined by 
  \begin{align*}
  \beta(k,t)\defeq \sqrt{1-2\me^{\frac{\pi}{2}\sin{t}}\cos\Big(\frac{\pi}{2}\cos{t}\Big)+\me^{\pi\sin{t}}}-{}&\Big(R^{\frac{2}{k}}-2R^{\frac{1}{k}}\cos\Big(t_{w_1,w_2}+\frac{2\pi}{k}+\frac{\pi}{2k}\cos{t}\Big)+1\Big)^{-\frac{k}{2}} \\ {}& -21|u_A|^{-1}2^{-k/2}.
  \end{align*}
  We claim that for big enough $k$, $\be(k,t)\gt 0$ for all $t\in[0,2\pi]$; recall the definition \eqref{eq: Definition of tw1w2} of $t_{w_1,w_2}$ to see that in this case,
  $
  t_{w_1,w_2}= \arccos((|w_2|^{2/k}-2)/2). 
  $
  Therefore,
  \begin{align*}
  \beta(k,t)\xrightarrow[k\to\infty]{} {}&\sqrt{1-2\me^{\frac{\pi}{2}\sin{t}}\cos\Big(\frac{\pi}{2}\cos{t}\Big)+\me^{\pi\sin{t}}}-\exp\Big(-\frac{\pi}{4}(4\sqrt{3}+\sqrt{3}\cos{t}-\sin{t})\Big) \\
  \gt {}& 0.79, 
  \end{align*}
  where the last inequality follows from a numerical computation. We just showed that the $\#Z_{\tilde{C}_2}-2$ curves $z=|w_2|^{\frac{1}{k}}\me^{\ii(t_0+\frac{\pi}{2k}\me^{\ii t})}-1$ satisfy 
  \[
  |f_2(z)|\gt \Big|\frac{1}{z^k}+\frac{u_A^{-1}}{2}\cs{R}_k^{a,b}(z)\Big|,
  \] 
  where $a=b$ or $a=-b$. By Rouché's theorem, $E_k^{a,b}(z)=2u_A(f_2(z)+{1}/{z^k}+\cs{R}_k^{a,b}(z))$ and $2u_A f_2(z)$ have the same number of zeros inside $\ga'_{t_0}$. Finally, since the diameter of $\ga'_{t_0}$ is $O(1/k)$, every zero of $E_k^{a,b}$ inside $\ga'_{t_0}$ is within $O(1/k)$ distance from the circle $C_2$. 
\end{proof}
We are ready to combine our results regarding the count of zeros into the following proposition.  
\begin{proposition}[Total number of zeros in $\F$]\label{prop: Total number of zeros in F}
Let $N$ be a fixed odd level.  Then, for sufficiently large even weight $k$, we have 
\begin{enumerate} 
\item for all $a,b$ with $a\neq \pm b$, and $|\cos(2\pi b)|\gt |\cos(2\pi a)|$, there is some positive constant $\al$ such that $\#Z_k^{a,b} \ge k/6 - 2\al\log N$;
\item for all $a$ and either choice of sign, $\# Z_k^{a,\pm a} \geq k/12-5$.
\end{enumerate}
\end{proposition}
\begin{proof}
    \begin{enumerate} 
\item Let $\tilde{m}$ be as in Lemma \ref{lem: Inequality for using Rouche}, and let $Z_{\tilde{C_1}}$ be as in Lemma \ref{lem: Zeros of fw1}. By Corollary \ref{cor: main result for a different than b}, for every $a,b$ with $a\neq \pm b$ with $|w_1|\gt 1$, $E_k^{a,b}$ has at least $Z_{\tilde{C_1}}-(\tilde{m}+1)$ zeros in $\F_1$. Now, upon recalling the definition of $ t_{w_1,w_2}$ from Equation \eqref{eq: Definition of tw1w2} 
     and using the fact that 
    \[
    \frac{k}{2\pi}t_{w_1,w_2} -\frac{k}{4}-\frac{k}{12}\xrightarrow[k\to\infty]{} \frac{2\log|w_1|-3\log{|w_2|}}{2\sqrt{3}\pi},
    \]
    we see that (with absolute constants) 
    $$  Z_{\tilde{C}_1}={}  \frac{k}{2\pi}t_{w_1,w_2} -\frac{k}{4}+O(1)
      ={}  \frac{k}{12}+O(2\log|w_1|-3\log{|w_2|}). $$
    Now, since $|u_j|\le 1$ for any $j$, it is clear that 
 $
    2\log|w_1|-3\log{|w_2|}\ll \log{|u_A|^{-1}},
   $
    and since $N$ is odd, we see that $|u_A|\ge 1/N$ by \eqref{eq: LowerBound uA}. This shows that 
    \begin{equation}\label{eq: UpperBound log(1/uA)}
    \log{|u_A|^{-1}}\le \log{N}, 
    \end{equation}
   and since $\tilde{m}\ll \log{N}$ it follows that there is some positive constant $\al$ such that $E_k^{a,b}$ has at least $k/12-\al\log{N}$ zeros in $\F_1$. Finally, we recall that by \eqref{eq: IdentityForZk}, 
$
    \#Z_{k}^{a,b}=\#Z_1^{a,b}+\#Z_1^{a,-b}+\#\tilde{Z}^{a,b},
$
    where $Z_1^{a,b}$ is the set of zeros of $E_k^{a,b}$ inside $\F_1$. Therefore, 
$
    \#Z_{k}^{a,b}\ge {k}/{6}- 2 \al \log{N}.
$
\item Let $\tilde{l}$ and $\lambda_k$ be as in Lemma \ref{lem: Function h(t)}. Then, by Corollary \ref{cor: Main theorem for a,a} for all $a$, $E_k^{a,\pm a}$ has at least
  \begin{equation}\label{eq: NumberOfZerosUnitCircle}
\Big\lfloor\frac{k}{12}\Big\rfloor-\tilde{l}-\lambda_k\ge \frac{k}{12}-\bigg\lceil \frac{\log|w_2|-\log{2}}{\sqrt{3}\pi}-\frac{\mu_k}{3}\bigg\rceil-2
\end{equation}
 zeros on $[\frac{\pi}{2},\frac{2\pi}{3}]$ lying exactly on the unit circle. Let $Z_{\tilde{C}_2}$ be as in Lemma \ref{lem: Zeros of fw2}. Then, by Lemma \ref{lem: ZerosCloseToC2}, $E_k^{a,\pm a}$ has at least $\#Z_{\tilde{C}_2}-2$ zeros on $\F_1$ lying outside the unit circle. Upon recalling that in this case 
$t_{w_1,w_2}= \arccos((|w_2|^{2/k}-2)/2)$ and $t'_{w_2}=\arccos({1}/(2|w_2|^{1/k}))$,
    we see that 
    \begin{align*}
       Z_{\tilde{C}_2}={}  \Big\lfloor \frac{k}{2\pi}t'_{w_2}-\frac{\delta_2}{2} \Big\rfloor -\Big\lceil \frac{k}{2\pi}t_{w_1,w_2}-\frac{\delta_2}{2}\Big\rceil+1 
      \ge{}  \frac{k}{2\pi}(t'_{w_2}-t_{w_1,w_2}) 
      \xrightarrow[k\to\infty]{}{}  \frac{\log|w_2|}{\sqrt{3}\pi}.
    \end{align*}
    In particular, since $\frac{k}{2\pi}(t'_{w_2}-t_{w_1,w_2})$ is decreasing for $|w_2|\gt 1$, for sufficiently large $k$, 
$
    Z_{\tilde{C}_2}\ge ({\log|w_2|}-{\log{2}}))/(\sqrt{3}\pi).
$
    This together with \eqref{eq: NumberOfZerosUnitCircle} shows that 
    \[
      \#Z_k^{a,\pm a}\ge \frac{k}{12}-\frac{\log|w_2|-\log{2}}{\sqrt{3}\pi}+\frac{\mu_k}{3}-3+\#Z_{\tilde{C}_2}-2\ge \frac{k}{12}-5. \qedhere
   \]
    
\end{enumerate}
\end{proof}
\subsection{Zeros (exponentially) close to $\ii$ and $\rho$}

In this section we analyze when $E_k^{a,b}$ has zeros exponentially close to $\ii$ and $\rho$. We will frequently make use of the Taylor expansion of a holomorphic function with precise error term determined by Cauchy's estimates (e.g., \cite[p. 190]{BrownChurchill}): 
if a domain $U$ contains the closed disk $\{ z \in \Cc \, | \, |z - z_0| \leq r \}$ and $f \colon U \to \Cc$ is a holomorphic function, then
    \[
    \Big| f(z) - \sum_{n = 0}^{B-1} \frac{f^{(n)}(z_0)}{n!} (z-z_0)^n \Big| \leq U_{r} \left(\frac{|z-z_0|}{r}\right)^{B} \mbox{ with } U_{r} = \sup_{|w-z_0| = r} |f(w)|.
 \]

\subsubsection*{Zeros close to $\ii$}

Recall the main term $\cs{M}_k^{a,b}(z)$ from \eqref{eq: MainTerm}. 
When $z$ is close to $i$, the term $u_{A-B}/(z+1)^k$ is exponentially small so that $u_A+u_B/z^k$ becomes dominant. If this dominant term is non-zero at $\ii$ we do not expect zeros exponentially close to $\ii$. This is made precise in the following lemma.
\begin{lemma}
\label{lem: LemExpoI}
Let $d \in \{0 , 2\}$.
Suppose $k \equiv d \mod 4$ and suppose
$
u_A + (-1)^{d/2} u_B \neq 0.
$
Then there exits a constant $c$ depending only on $a$ and $b$ such that the region
$
\{ z \in \mathcal{F} \, | \, |z - \ii| < \frac{c}{k} \}
$
contains no zeros of $E_k^{a,b}(z)$ for sufficiently large $k$. 
\end{lemma}
\begin{proof}

If, on the contrary, we can find a subsequence of zeros $z_k = i ( 1 + \frac{\epsilon_k}{k})$ of $E_k^{a,b}$ with $\lim\limits_{k \to \infty} \epsilon_k = 0$, the value of the main term is given by
\begin{align*}
\cs{M}_k^{a,b}(z_k) &=  2\Big(u_A+(-1)^{d/2}\f{u_B}{(1 + \frac{\epsilon_k}{k})^k}+\f{u_{A-B}}{(z_k+1)^k}\Big).
\end{align*}
We apply the Taylor expansion  to $f(z) = 1/(1+z)^{k}$ with $z = \frac{\epsilon_k}{k}$, $r 
=  1/k$ and $z_0 = 0$, to find
\[
    \frac{1}{(1 + \frac{\epsilon_k}{k})^k} = 1 + s_1\left(\frac{\epsilon_k}{k}\right)
\mbox{ with }
|s_1\left(\frac{\epsilon_k}{k}\right)| \leq U_r |\epsilon_k|.
\]
By the triangle inequality we have 
$
U_r = \sup \{|f(w)| \colon {|w| = \frac{1}{k}} \} \leq (1-{1}/{k})^{-k}.
$
It follows that $U_r$ is bounded as $k \to \infty$, and we conclude that $\lim\limits_{k \to \infty} s_1\left(\frac{\epsilon_k}{k}\right) = 0$. 
Since $\lim\limits_{k \to \infty} z_k = \ii$, it is clear that
$
\lim\limits_{k \to \infty} {u_{A-B}}/{(z_k+1)^k} = 0.
$
We rewrite the expression for $E_k^{a,b}$, putting terms with vanishing limits in $k$ on the right hand side.
\begin{align*}
    E_k^{a,b}(z_k) -2 \left( u_A + (-1)^{d/2} u_B \right) =\frac{2}{\ii^k} s_1\left(\frac{\epsilon_k}{k}\right) + 2\frac{u_{A-B}}{(z_k + 1)^k} + \cs{R}_k^{a,b}(z_k).
\end{align*}
Taking the limit as $k \to \infty$ on both sides, using $E_k^{a,b}(z_k) = 0$ and the bound \eqref{eq: ExpDecreasing}, we find
$
u_A +(-1)^{d/2} u_B = 0,
$
contradicting the assumption.  
\end{proof}

Lemma \ref{lem: LemExpoI} implies that zeros of $E_k^{a,b}$ can be exponentially close to $\ii$ only when $u_{A} = -u_{B}$ if $k \equiv 0 \mod 4$, or $u_{A} = u_{B}$ if $k \equiv 2 \mod 4$. We only treat the case $k \equiv 2 \mod 4$, since the equation $u_A = -u_B$ has no solutions for odd level $N$. Recall the description of zeros of the Eisenstein series $E^{a,b}_k$ on $\A$ in terms of the real-valued function $h^{a,b}_k$. 

\begin{lemma}
\label{lem: signhkfori}
    Suppose $k \equiv 2 \mod 4$, $|a| = |b|$, and $a \neq \pm \frac{1}{4}$. Let $1 < \beta < \sqrt{2}$. Then for sufficiently large $k$, 
    \begin{equation} \label{eq:signlem}
    \sgn h_k^{a,b} \left(\frac{\pi}{2} + \beta^{-k}\right) = (-1)^{(k+2)/4}.
    \end{equation}
\end{lemma}

\begin{proof}
    Write $t_k = \frac{\pi}{2} + \beta^{-k}$ and plug this into $h_k^{a,b}$ as in \eqref{eq: Definition of h(t)}.
Suppose $\beta < \eta < \sqrt{2}$. The limit $\lim\limits_{k \to \infty} 2 \cos(t_k/2) = \sqrt{2}$ and the bound \eqref{eq: ExpDecreasing} imply that 
\begin{equation}
\label{eq: bound1}
    \left|h_k^{a,b}(t_k) - 2\cos(t_k k/2) \right| \leq \eta^{-k}
\end{equation}
for sufficiently large $k$.
Since $k \equiv 2 \mod 4$, we have $\cos(t_k k/2) = (-1)^{(k+2)/4} \sin (k \beta^{-k}/2)$. Moreover, $\sin(x) \geq \frac{2}{\pi} x$ for $x \in [0, \frac{\pi}{2}]$. This means that for large enough $k$
\begin{equation}
\label{eq: bound2}
|\cos(t_k k/2)| \geq \frac{k \beta^{-k}}{\pi}.     
\end{equation}
Finally, the bounds \eqref{eq: bound1} and \eqref{eq: bound2} imply that \eqref{eq:signlem} holds  for large enough $k$. 
\end{proof}

\begin{lemma}
\label{lem: snghkati}
    Suppose $k \equiv 2 \mod 4$, $|a| = |b|$ and $a \not \in \{0, \frac{1}{4}, \frac{1}{2}, \frac{3}{4} \}$.
    Then
    \[
    \sgn h^{a,b}_k\left( \frac{\pi}{2}\right) = \sgn u_A^{-1}\left( u_{A-B}-u_{A+B} \right).
    \]
\end{lemma}
\begin{proof}
In this case
\begin{equation}\label{eq: E_k(i)}
    E_k^{a,b}(\ii) = 2(2\ii)^{-k/2}(u_{A-B}-u_{A+B}) +  \sum_{m^2+n^2\ge 4}\frac{\me(an-bm)}{(m \ii+n)^k}
\end{equation}
and
\begin{equation}\label{eq: R_k(i)}
\sum_{m^2+n^2\ge 4}\frac{\me(an-bm)}{(m \ii+n)^k} \le \int_{4}^{\infty}(\sqrt{x}+1)x^{-k/2}\dx \leq C \cdot 2^{-k}
\end{equation}
for some positive constant $C$.
Using the definition of $h^{a,b}_k(t)$ from \eqref{defht}, this implies the result. \end{proof}
A direct consequence of Lemmas \ref{lem: signhkfori} and \ref{lem: snghkati} is the following proposition, showing the existence of zeros exponentially close to $\ii$.

\begin{proposition}
\label{thm: expozerosclosetoi}
Let $d \in \{2 , 6 \}$.
    Suppose $k \equiv d \mod 8$, $|a| = |b|$, and $a \not \in \{0, \frac{1}{4}, \frac{1}{2}, \frac{3}{4} \}$. Furthermore assume
    $
    (-1)^{(d-2)/4} = \sgn u_A^{-1}\left( u_{A-B}-u_{A+B} \right).
    $
    Let $1 < \beta < \sqrt{2}$. Then for sufficiently large $k$, 
        $E_k^{a,b}(\me^{\ii t})$ has a zero  with
    $t \in (\frac{\pi}{2}, \frac{\pi}{2} + \beta^{-k}).
    $ \qed
\end{proposition}

\subsubsection*{Zeros close to $\rho$} The reasoning here is somewhat similar to the one at $\ii$, so we will leave out some of the computations. 
We first show that $E_k^{a,b}$ has no zeros exponentially close to $\rho$ if the main term does not vanish. The proof is very similar to that of Lemma \ref{lem: LemExpoI} and is left out. 

\begin{lemma}
\label{lem: LemExpoRho}
Let $d \in \{0, 1, 2\}$. Suppose $k \equiv d \mod 3$ and suppose
$
    \rho^{2d} u_A + \rho^d u_B + u_{A-B} \neq 0.
$
Then there exits a constant $c$ depending only on $a$ and $b$ such that the region
$
\{ z \in \mathcal{F} \, | \, |z - \rho| < \frac{c}{k} \}
$
contains no zeros of $E_k^{a,b}(z)$ for sufficiently large $k$. \qed
\end{lemma}

Lemma \ref{lem: LemExpoRho} implies that zeros exponentially close to $\rho$ can only come from the pairs $(A,B)$ satisfying the equation
\begin{equation} \label{eqrhou}
    \rho^{2d} u_A + \rho^d u_B + u_{A-B} = 0.
\end{equation}
that we now study in detail. 

\begin{lemma}
\label{Lem:rho}
Suppose $k \equiv d \mod 3$. 
    $\cs{M}_k^{a,b}(\rho)$ vanishes if and only if 
    \[
    (a,b) \in 
    \begin{cases}
         \{ \pm (\frac{1}{3}, \frac{1}{3}), (0,\pm \frac{1}{3}), (\pm \frac{1}{3}, 0) \} & \text{ if } d = 0, \\
        \{ (0,0), \pm (-\frac{1}{3}, \frac{1}{3})\} & \text{ if } d = 1 \mbox{ or } d=2. 
            \end{cases}
    \]
\end{lemma}

\begin{proof}  
We split $k$ into congruence classes modulo $3$.

\textbf{Case 1:}  $k \equiv 0 \mod 3$.  The equation becomes
    $
     u_A +  u_B + u_{A-B} = 0.
    $
    The only solutions to this equation are the six solutions $(a,b) \in \{\pm ( \frac{1}{3}, \frac{1}{3}), (0,\pm \frac{1}{3}), (\pm \frac{1}{3}, 0) \}$, see \cite{Cox}.

\textbf{Case 2:}  $k \equiv 1 \mod 3$.  The equation becomes
$
     \rho^2 u_A + \rho  u_B + u_{A-B} = 0.
$
Taking both the real and imaginary parts on both sides, we find
\begin{align*}
    -\frac{3}{2} u_A + \frac{1}{2} u_B + u_{A-B} = 0 \mbox{ and }
    -\frac{\sqrt{3}}{2} u_A + \frac{3}{2} u_B = 0.
\end{align*}
Hence we need to solve $u_A = u_B = u_{A-B}$. Now $u_A = u_B$ if and only if $a = b$ or $a = -b \mod 1$. An equality
$a = b \mod 1$ implies $u_A = u_B = 0$ and therefore $A = B = 0$.
An equality $a = -b \mod 1$ implies $u_{2A} = u_A$ and this implies $a = 0 \mod 1$ or $3 a = 0 \mod 1$. Hence $a = \pm  1/3 \mod 1$. We conclude that there are the  three solutions $(a,b) \in \{ (0,0), \pm (\frac{1}{3}, -\frac{1}{3})\} $.

\textbf{Case 3:}  $k \equiv 2 \mod 3$.  The equation becomes
    $     \rho u_A + \rho^2  u_B + u_{A-B} = 0.
    $
This is the same as Case $2$, interchanging $A$ and $B$. The result follows. 
\end{proof}

\begin{remark}
The question of solving \eqref{eqrhou} can be rephrased as finding cyclotomic points on the elliptic curve
$
    \rho^{2d} (X+1/X) + \rho^d (Y+1/Y) + X/Y + Y/X = 0.
$
In \cite{BeukSmyth} an algorithm is presented to find these points for general curves over $\Cc$. For our specific family, explicit results are obtained in \cite{Cox}, as used in the proof of Lemma \ref{Lem:rho}.
\end{remark}

Like we previously focussed on the real-valued version of the Eisenstein series on the arc close to $\ii$, we now focus on the real-valued function given by the series $
E_k^{0, 1/3}(z)
$ on the line $\Re(z)=-1/2$ close to $\rho$. 

\begin{lemma}
\label{lem: real-valued}
For $t > 0$, 
$
E_k^{0, 1/3}(-\frac{1}{2} + \ii t) \in \bb{R}.
$
\end{lemma}
\begin{proof}
An easy application of $\eqref{eisrel2}$ and $\eqref{eisrel4}$ shows that $E_k^{0, 1/3}(-\frac{1}{2} + \ii t)$ is invariant under complex conjugation.  
\end{proof}

\begin{lemma}
\label{lem: signchangelemrho}
    Suppose $k \equiv 0 \mod 3$ and $(a,b) = (0, \frac{1}{3})$. For $1 < \beta < \sqrt{2}$ 
    and sufficiently large $k$, 
    $
    \sgn E_k^{0, 1/3}(\rho + i \beta^{-k}) = 1.
    $
\end{lemma}
\begin{proof}
In this situation, the main term $\cs{M}_k^{0, 1/3}(z)$ vanishes at $z = \rho$. This means that to closely analyze the main term near $\rho$, we need to use the Taylor series. We first apply it to $f(z) = (1+z)^{-k}$ with $z = \frac{\ii}{\rho \beta^k}$, $c = 0$ and $r = \beta^{-k/4}$, to find
\[
\frac{1}{(1 + \frac{\ii}{\rho\beta^{k}})^k} =1 - k\frac{\ii}{\rho \beta^{k}} + s_2 \left( \frac{\ii}{\rho \beta^{k}} \right),
\mbox{ where }
\Big|s_2 \left( \frac{\ii}{\rho \beta^{k}} \right)\Big| \leq U_r \beta^{-3k/2}.
\]
with $U_r$ is bounded as $k \to \infty$.
We do a similar analysis for the term \[
\frac{1}{(\rho + 1 + \frac{\ii}{\beta^k})^k} = \frac{1}{\rho^{2k}} \frac{1}{(1 - \frac{\ii}{\rho^2 \beta^k})^k}
\]
and, collecting everything and using $k \equiv 0 \mod 3$, we find 
\[
|\cs{M}_k^{0, 1/3}(\rho + \ii \beta^{-k}) - 2 k \sqrt{3} \beta^{-k} | \leq 4 U_r \beta^{-3k/2}. 
\]
Let $\beta < \eta < \sqrt{2}$. The bound \eqref{eq: ExpDecreasing}
implies
$
|\cs{R}_k^{0, {1}/{3}}(\rho + \ii \beta^{-k})| \leq \eta^{-k}
$
for sufficiently large $k$.
We conclude that 
   $
    \sgn E_k^{0,1/3}(\rho + \ii \beta^{-k}) = \sgn  2 k \sqrt{3} \beta^{-k} = 1
   $
for sufficiently large $k$.
\end{proof}

\begin{lemma}
\label{lem: leadingrhoterm}
Suppose $k \equiv 0 \mod 3$ and 
$
u_A + u_B + u_{A-B} = 0. 
$
Then
\[
|E_k^{a,b}(\rho) + 3 \cdot (-27)^{k/6}| \leq 2 \cdot 7^{-k/2} + C \cdot 2^{-k}
\]
for some positive constant $C$. 
\end{lemma}
\begin{proof}
In this case we choose as new `main term' at $\rho$
\[
\widetilde{\cs{M}_k}^{a,b} = \sum_{(m,n) \in I} \frac{\e(an - bm)}{(m \rho + n)^k},
\]
where $I = \{ (\pm 1, 0), (0, \pm 1), \pm(1,1), \pm(1,-1), \pm(1,2), \pm(2,1) \}$. The remainder can be estimated as 
$
|E_k^{a,b}(\rho) - \widetilde{\cs{M}_k}^{a,b}| \leq 2 \cdot 7^{-k/2} + C \cdot 2^{-k},
$
for some constant $C>0$. 
Using $u_A + u_B + u_{A-B} = 0$, we can write the main term as the non-vanishing quantity 
$
\widetilde{\cs{M}_k}^{a,b} = 2 (-27)^{k/6} \left( u_{A+B} + u_{2A-B} + u_{A-2B}\right) = -3 \cdot (-27)^{k/6}. $
\end{proof}

\begin{proposition}
\label{thm: exporhozero}
Suppose $k \equiv 0 \mod {12}$ and $1 < \beta < \sqrt{2}$.
Then $E_k^{0, 1/3}(-1/2 + \ii t)$ has a zero for some $t \in (\sqrt{3}/2, \beta^{-k})$. 
\end{proposition}

\begin{proof}
    By Lemmas \ref{lem: real-valued} and \ref{lem: leadingrhoterm} it follows that for sufficiently big $k$, the sign of $E_k^{0, 1/3}(\rho)$ is equal to the sign of $-3 \cdot (-27)^{k/6}$, viz., $- (-1)^{k/6}.$ The result follows by combining this with Lemma \ref{lem: signchangelemrho}.
\end{proof}

\subsubsection*{Zeros precisely at $\ii$ or $\rho$} 

\begin{proposition}[Vanishing of $E_k^{a,b}$ at $\ii$ and $\rho$] \label{zerosat} 
    For all sufficiently large even $k$, and odd level $N$, we have
    \begin{enumerate}
        \item \label{vani} $E_k^{a,b}(i)=0$ if and only if $(a,b)=(0,0)$ and $k\equiv 2\mod{4}$;
        \item \label{vanishingrho} $E_k^{a,b}(\rho)=0$ if and only if $(a,b)=(0,0)$ and $k\equiv 1,2\mod{3}$, or $(a,b)=\pm (\f13,-\f13)$ and $k\equiv 2,4 \mod{6}$.
    \end{enumerate}
\end{proposition}
\begin{proof}
The fact that $E_k^{0,0}(\ii)=0$ if and only if $k\equiv 2\mod{4}$, and $E_k^{0,0}(\rho)=0$ if and only if $k\equiv 1,2\mod{3}$ is a classical result, and can be found in \cite{GekArch}, so we assume that $(a,b)\neq (0,0)$.

For case   \eqref{vani}, by Lemma \ref{lem: signhkfori}, if $u_A+(-1)^{\f{k}{2}}u_B\neq 0$, then $E_k^{(a,b)}(i)\neq 0$. If $u_A+(-1)^{\f{k}{2}}u_B= 0$, since $N$ is odd, the equation $u_A=-u_B$ has no solutions, so the only possibility is $u_A=u_B$, in which case $A=\pm B$. Then, since $A,B$ are not simultaneously zero, we have $u_{A-B}-u_{A+B}\neq 0$. Therefore, by \eqref{eq: E_k(i)} and \eqref{eq: R_k(i)}, for large enough $k$, we have 
        \[
|E_k^{a,b}(\ii)-2(2\ii)^{-k/2}(u_{A-B}-u_{A+B}) |\le C \cdot 2^{-k} \lt |2(2\ii)^{-k/2}(u_{A-B}-u_{A+B})|,
    \]    
    showing that $E_k^{(a,b)}(\ii)\neq 0$ in this case. 

For case \eqref{vanishingrho}, we claim that that $E_{k}^{1/3,-1/3}(\rho)=0$ when $k\equiv 2,4 \mod{6}$. %
Indeed, 
\begin{align*}
    E_{k}^{1/3,-1/3}(&\rho) = \sideset{}{'}\sum_{m,n \in \Z} \frac{\e((n + m)/3)}{(m \rho + n)^k} 
    = \frac{1}{\rho^k} \sideset{}{'}\sum_{m,n \in \Z} \frac{\e((n + m)/3)}{(m -n (\rho + 1))^k} \\
     &\overset{\left\{\begin{subarray}{l}\tilde{m} = -n\\ \tilde{n} = m - n\end{subarray}\right.}{=} \frac{1}{\rho^k} \sideset{}{'}\sum_{\tilde{m},\tilde{n} \in \Z} \frac{\e((\tilde{n} + \tilde{m})/3) \e(-3\tilde{m}/3)}{(\tilde{m}\rho + \tilde{n})^k} = \frac{1}{\rho^k} E_{k}^{1/3,-1/3}(\rho).
\end{align*}
This implies that $E_{k}^{1/3,-1/3}(\rho) = 0$ whenever $k \not \equiv 0 \mod 3$. The fact that $E_{k}^{-1/3,1/3}(\rho)=0$ then follows  from \eqref{eisrel1}.
Now, assume that $(a,b)\neq \pm (\f13,-\f13)$. If $\cs{M}_k^{a,b}(\rho)\neq 0$, then by Lemma \ref{lem: LemExpoRho}, $E_k^{a,b}(\rho)\neq 0$. If $\cs{M}_k^{a,b}(\rho)$ vanishes, as in the proof of Lemma \ref{Lem:rho}, we have $k\equiv 0\mod{3}$ and $u_A+u_B+u_{A-B}=0$. Then, by Lemma \ref{lem: leadingrhoterm}, for sufficiently large $k$, we have
    \[
|E_k^{a,b}(\rho) + 3 \cdot (-27)^{k/6}| \leq 2 \cdot 7^{-k/2} + C \cdot 2^{-k} \lt 3\cdot 27^{k/6},
\]
showing that $E_k^{a,b}(\rho)\neq 0$. \qedhere
\end{proof}
\subsection{Proof of Theorem \ref{main2}}
The results of the previous subsections can now be combined to prove Theorem \ref{main2}. 

\begin{proof}[Proof of Theorem \ref{main2}] \mbox{ }
  \begin{enumerate}[leftmargin=*] 
\item Let 
$
\cs{B}\defeq \{ (a,b) \; : \; |w_1|\gt 1\}.
$
By Corollary \ref{cor: main result for a different than b}, we know that for $(a,b)\in\cs{B}$, $E_k^{a,b}$ has at least $Z_{\tilde{C_1}}-(\tilde{m}+1)$ zeros in $\F_1$ that are simple, and at Euclidean distance $\asymp 1/k$ from the unit circle. Let $L_1^{a,b}$ denote the set of these zeros for fixed $a,b$. In particular, since these zeros do not lie on the unit circle, they are different from $\ii$ and $\rho$, so by Theorem \ref{tran}, all of these zeros are transcendental for sufficiently large $k$. Let $$\eta(z)\defeq -\ol{z},$$ and observe that in view of the identity $E_k^{a,b}(-\ol{z})=\ol{E_k^{a,-b}(z)}$, $\eta(L_1^{a,-b})$ consists of simple, transcendental zeros of $E_k^{a,b}$ at Euclidean distance $\asymp 1/k$ from the unit circle lying in $\F\setminus \F_1$. Set 
\begin{equation}\label{eq: DefinitionHatZk}
\widehat{Z}_k^{\text{off}}\defeq \bigcup_{(a,b)\in\cs{B}}L_1^{a,b}\cup \eta(L_1^{a,-b}).
\end{equation}
(There might also be further zeros on the imaginary axis, i.e., in $\F \setminus (\F_1 \cup \eta(\F_1))$, but we will not have to count those.) By construction, it is clear that 
$
d_E(\widehat{Z}_k^{\text{off}}, \A \cup \A') \asymp 1/k.
$
Since $|w_1|=|u_B u_A^{-1}|$, the fact that $(a,b)\in \cs{B}$ implies $(b,a)\notin \cs{B}$. Also, note that $(a,\pm a)\notin \cs{B}$. Therefore
\begin{equation*}
\sum_{\substack{a\neq \pm b \\ |w_1|\gt 1}}1=\frac{1}{2}\sum_{a\neq \pm b}1= \frac{(N-1)^2}{2}.
\end{equation*} 
Now recall that in the proof of Proposition \ref{prop: Total number of zeros in F} we showed that 
\begin{equation}\label{eq: LowerBoundForL1(a,b)}
\#L_1^{a,b}\ge \f{k}{12} -\al \log N
\end{equation}
for some absolute constant $\al$, and so
\begin{equation}\label{eq: lower bound for size of hat Z_k}
    \# \widehat{Z}_k^{\text{off}} \ge \frac{(N-1)^2 k}{12} -c_{\textup{off}} N^2 \log{N}.
\end{equation}
for some absolute constant $c_{\textup{off}}$.
\item This follows from Corollary \ref{cor: Main theorem for a,a}, noting that by \eqref{eq: UpperBound log(1/uA)}, we have $\log{|w_2|}\ll \log{|u_A|^{-1}}\le \log{N}$ with an absolute implied constant, which we can use to bound $\tilde l$ in Equation \eqref{deftildel}.

\item This follows immediately from Propositions  \ref{thm: expozerosclosetoi} and \ref{thm: exporhozero}. \qedhere  
\end{enumerate}
\end{proof}

\subsection{Angular equidistribution} We now study the distribution of the arguments of the zeros of Eisenstein series for $\Gamma(N)$ with fixed odd level $N$. For an even integer $k \geq 4$, we define the \emph{angular discrepancy} in $\F_1$ as 
$$ D_{k} \defeq \sup_{\frac{\pi}{2} \leq \theta_1 \leq \theta_2 \leq \frac{2\pi}{3}} \Big| \frac{\#(\arg(Z_{k}) \cap [\theta_1,\theta_2])}{\# (Z_{k}\cap \F_1)} - \frac{6}{\pi} (\theta_2-\theta_1)\Big|,$$
where the zeros are counted with multiplicity and weights. We call the set of zeros $Z_{\Gamma(N)}$ \emph{angularly equidistributed in $\A$} if $D_{k} \rightarrow 0$ with $k \rightarrow +\infty$. By the usual argument of approximating continuous functions by step functions (compare, e.g., \cite[Thm.~1.1]{KN}), 
this implies a weak convergence of measures 
$$ \frac{1}{\# (Z_{k}\cap \F_1)} \sum_{z \in Z_{k}} \delta_{\arg(z)} \rightarrow \frac{6}{\pi}\d\theta$$
on $[\pi/2,2\pi/3]$. 

\begin{proof}[Proof of Proposition \ref{equidist}] We continue to use the notation introduced in the previous proofs. In this proof, all $O(\cdot)$ terms carry absolute constants. 
{
	For $N=1$, the result follows from \cite{RSD}: their
	sign-change argument gives at least one zero between
	consecutive angles $2\pi m/k$ in $[\pi/2,2\pi/3]$.
	Together with the valence formula, this shows that in level $N=1$
	the number of zeros in any angular subinterval
	$[\theta_1,\theta_2]$ is
	\[
	\frac{k}{2\pi}(\theta_2-\theta_1)+O(1),
	\]
	uniformly in $\theta_1,\theta_2$, and hence the
	discrepancy is $O(1/k)$. We may therefore assume
	$N\geq 3$.}

By Corollary \ref{cor: Main theorem for a,a}, for either choice of sign, $E_k^{a,\pm a}$ has at least $\floor{k/12}-\tilde{l}-\lambda_k$ zeros in $\F$, where $\lambda_k\in \{0,1\}$, and $\tilde{l}\ll \log{N}$ by \eqref{eq: UpperBound log(1/uA)}. Let $L_2^{a,\pm a}$ denote the set of such zeros, and define 
\[
\widehat{Z}_k^{\text{on}}\defeq \bigcup_{a} L_2^{a,a}\cup L_2^{a,-a},
\]
so that 
\begin{equation}\label{eq: LowerBoundZk'}
    \#\widehat{Z}_k^{\text{on}}\ge \frac{(2N-1)k}{12}-c_{\textup{on}} N\log{N}
\end{equation}
for some absolute constant $c_{\textup{on}}$. Now we can write $\#Z_k=\#\widehat{Z}_k^{\text{off}} + \# \widehat{Z}_k^{\text{on}}+ Q$, where $\widehat{Z}_k^{\text{off}}$ is as in \eqref{eq: DefinitionHatZk}, and $Q\defeq \#Z_k-\#\widehat{Z}_k^{\text{off}} -\# \widehat{Z}_k^{\text{on}}$. By Theorem \ref{main2} (i), and \eqref{eq: LowerBoundZk'}, we have $Q\ll N^2\log{N}$. Hence
\begin{equation}\label{eq: arg(Zk)}
    \#(\arg(Z_{k}) \cap [\theta_1,\theta_2])=\#(\arg(\widehat{Z}_k^{\text{off}}) \cap [\theta_1,\theta_2])+\#(\arg(Z_{k}^{\text{on}}) \cap [\theta_1,\theta_2])+O(N^2\log{N}).
\end{equation}

\subsubsection*{Angles for $\widehat{Z}_k^{\text{off}}$}
By the definition of $\widehat{Z}_k^{\text{off}}$, since $\eta(L_1^{a,-b})\notin \F_1$, we have
\begin{equation}\label{eq: arg(hat(Zk))}
\#(\arg(\widehat{Z}_k^{\text{off}}) \cap [\theta_1,\theta_2])=\f12 \sum_{a\neq \pm b } \#(\arg(L_1^{a,b}) \cap [\theta_1,\theta_2]).
\end{equation}
Recall that in the proof of Corollary \ref{cor: main result for a different than b}, we showed that each zero of $L_1^{a,b}$ lies uniquely inside $\ga_{t}$ for some $t$ such that $f(|w_1|^{1/k}\me^{\ii t})=0$, so $t=(2m+\de_1)\pi/k$ for some integer $m$, where $\de_1=1$ if $\sgn(w_1)=1$, and $0$ otherwise. This shows that for each $\hat{\theta}\in \arg(L_1^{a,b})$ such that $z_0=R\me^{i\hat{\theta}}\in L_1^{a,b}$, there exists a $t_0=(2m_0+\de_1)\pi/k$ for some integer $m_0$ such that $z_0\in \mbox{int}(\ga_{t_0})$, where $\ga_{t_0}$ is as in \eqref{eq: Rouché contour gamma0}. Let $F\colon \bb{C} \to \bb{C}$ be defined as $F(w)\defeq |w_1|^{1/k}\me^{\ii(t_0+\f{c}{k}w)}$, and observe that since $F$ is a conformal map, $\mbox{int}(\ga_{t_0})=F(D)$, where $D$ is the unit disk. This shows that since $z_0\in \mbox{int}(\ga_{t_0})$, we have $z_0=|w_1|^{1/k}\me^{\ii(t_0+\f{c'}{k}\me^{\ii t})}$ for some $t\in [0,2\pi]$ and some $c'\in [0,c)$. Since $z_0=R\me^{i\hat{\theta}}$, it follows that 
$
\hat{\theta}=t_0+{c'}/{k}\cdot \cos{t}+2\pi m
$
for some integer $m$, but since $t_0\in [\pi/2, 2\pi/3]$, it follows that for sufficiently large $k$, it must be the case that $m=0$. Therefore, to each $\hat{\theta}\in \arg(L_1^{a,b})$ corresponds a unique integer $m_0$ such that 
$
\hat{\theta}={\pi}/{k}(2m_0+\delta_1)+\frac{c''}{k}
$
for some $c''\in (-c,c)$. 
Counting the number of integers $m_0$ such that $\hat{\theta}\in[\theta_1,\theta_2]$ shows that 
\[
\#(\arg(L_1^{a,b}) \cap [\theta_1,\theta_2])=\frac{k}{2\pi}(\theta_2-\theta_1)+O(1).
\]
After plugging this bound into \eqref{eq: arg(hat(Zk))}, we obtain 
\begin{equation}\label{eq: AsymptoticFor hat(Zk)}
    \#(\arg(\widehat{Z}_k^{\text{off}}) \cap [\theta_1,\theta_2])=\frac{k}{2\pi}(\theta_2-\theta_1)\frac{(N-1)^2}{2}+O(N^2).
\end{equation}
\subsubsection*{Angles for $\widehat{Z}_k^{\text{on}}$} Recall that in the proof of Lemma \ref{lem: Function h(t)}, we showed that for each zero $z\in L_2^{a,a}$ (which lies on the unit circle) there is a unique integer $m_1$ such that $2\pi(m_1-1)/k\lt \arg(z)\lt 2\pi m_1/k$. This shows that $\arg(z)=2\pi m_1/k - c_1/k$ for some $c_1\in (0,2\pi)$. Counting the number of $m_1$ such that $\arg(z)\in [\theta_1,\theta_2]$ gives us 
\[
\#(\arg(L_1^{a,a}) \cap [\theta_1,\theta_2])=\frac{k}{2\pi}(\theta_2-\theta_1)+O(1),
\]
and the same bound holds for $L_1^{a,-a}$. Hence
\begin{equation}\label{eq: AsymptoticFor Zk'}
    \#(\arg(\widehat{Z}_k^{\text{on}}) \cap [\theta_1,\theta_2])=\frac{k}{2\pi}(\theta_2-\theta_1)(2N-1)+O(N).
\end{equation}

\subsubsection*{Angles for $\widehat{Z}_k$} Plugging \eqref{eq: AsymptoticFor hat(Zk)} and \eqref{eq: AsymptoticFor Zk'} into \eqref{eq: arg(Zk)} shows that
\begin{equation}\label{eq: AsymptoticFor arg(Zk)}
 \#(\arg(Z_{k}) \cap [\theta_1,\theta_2])=\frac{k}{2\pi}(\theta_2-\theta_1)\Big(\frac{(N-1)^2}{2}+2N-1\Big)+O(N^2\log{N}).
\end{equation}
Since all the zeros of $\widehat{Z}_k^{\text{on}}$ lie in $\F_1$, and since all the zeros of $L_1^{a,b}$ with $(a,b)\in\cs{B}$ also lie in $\F_1$, by \eqref{eq: LowerBoundForL1(a,b)}, and \eqref{eq: LowerBoundZk'}, we have
\[
\#(Z_k\cap\F_1)\ge \frac{k}{12}\Big(\frac{(N-1)^2}{2}+2N-1\Big)+O(N^2\log{N}).
\]
For the corresponding upper bound, we use the fact that all the zeros of $\eta(L_1^{a,-b})$ lie in $\F\setminus \F_1$, together with \eqref{eq: LowerBoundForL1(a,b)}, and the valence formula in order to obtain
\[
\#(Z_k\cap\F_1)\le \# Z_k - \#(\widehat{Z}_k^{\text{off}}\cap (\F\setminus \F_1))\le \frac{k}{12}\Big(N^2-\frac{(N-1)^2}{2}\Big)+O(N^2\log{N}).
\]
Hence 
\[
\#(Z_k\cap\F_1)=\frac{k}{12}\Big(\frac{(N-1)^2}{2}+2N-1\Big)+O(N^2\log{N}).
\]
Together with \eqref{eq: AsymptoticFor arg(Zk)}, {and incorporating simultaneously the case $N=1$,} this shows that 
\[
\frac{\#(\arg(Z_{k}) \cap [\theta_1,\theta_2])}{\# (Z_{k}\cap \F_1)}= {\frac{6}{\pi}(\theta_2-\theta_1) + O\Big(\frac{1 + \log{N}}{k}\Big)} \qedhere 
\]
\end{proof} 

{ 

Let $N$ be odd,  and set $Z'_k \defeq Z_{k,\Gamma(N)}\cap \F_1$ for the zeros in the left half of the fundamental domain.  
Hausdorff convergence and angular equidistribution can be combinedto imply weak convergence for the entire normalized zero counting measure of $\mathcal N_k$ in $\F_1$, given as 
\[\mu_k \defeq \frac{1}{\#Z'_k} \sum_{z \in Z'_k} \delta_z.\] 

\begin{corollary} \label{cortotconv} As $k \rightarrow + \infty$ for even $k$, we have a weak convergence of measures
\[ \mu_k \rightarrow \frac{6}{\pi} |\d \arg z|\big|_{\A}, \]
the argument measure restricted to the arc $\A$. 
\end{corollary} 

\begin{proof} Write a complex number $z$ in polar coordinates as $ z =r_z e^{\ii \theta_z}$. By Hausdorff convergence, zeros in $Z_k'$ are uniformly close to $\A$ as $k \rightarrow + \infty$ ($k$ even): there exists $\varepsilon_k \rightarrow 0$ such that for all $z\in Z_k'$ we have  $|r_z-1| \leq \varepsilon_k.$
In particular, for $k$ sufficiently large, $Z_k' \subseteq \mathbb A \defeq \A \cap \{ 1/2 \leq |z| \leq 2\}$, a compact set. By the Heine--Cantor theorem, any $f \in C_b(\Hh,\R)$ is uniformly continuous on $\mathbb A$, so 
\[ \sup_{z \in Z'_k} |f(z)-f(e^{\ii \theta_z})| \leq  \omega_f(|r_z e^{\ii \theta_z}-e^{\ii \theta_z}|) \leq \omega_f(\varepsilon_k), \]
where the modulus of continuity $\omega_f$ satisfies $\omega_f(\varepsilon_k) \rightarrow 0$. Thus, 
\begin{equation} \label{equi1} \int f \d \mu_k = \frac{1}{\#Z'_k} \sum_{z \in Z'_k} f(e^{\ii \theta_z}) + O(\omega_f(\varepsilon_k)), \end{equation} 
and equidistribution of angles means 
\begin{equation} \label{equi2} \frac{1}{\#Z'_k} \sum_{z \in Z'_k} f(e^{\ii \theta_z}) \rightarrow \frac{6}{\pi} \int_{\pi/2}^{2 \pi/3} f(e^{\ii \theta}) \d \theta. \end{equation}
The result follows from \eqref{equi1} and \eqref{equi2}.
\end{proof} 

\begin{remark} A similar result holds for the zeros in the right half of the fundamental domain, switching $\F_1$ and $\A$ for $-\overline{\F_1}$ and $-\overline{\A}$, respectively. Using the pushforward of measures, this implies weak convergence of the zero counting measure on $Y(1)$ to the argument measure on the geodesic segment corresponding to $\A$. 
In a follow-up paper \cite{CM}, potential-theoretic methods are used to study the limit measure of the zero counting measure of Eisenstein series on the modular orbifold $Y(1)$ for general congruence groups $\Gamma$ as the even weight tends to infinity. The limit measure is concentrated on the image of the geodesic segments occurring in the Hausdorff limit, and equal, up to a normalizing constant, to the argument measure on those arcs. 
\end{remark}

}

\section{Open problems} 

We list some open problems. 

\begin{enumerate} 
\item Study Kluyver sums and find the optimal constants $\kappa_\Gamma$ for congruence groups of any type in $\mathrm{SL}_2(\Z/N\Z)$, e.g., for $N$ prime, when they are either $1$ or $N$; 
See Example \ref{exotic} for one case. 
\item Study lower bounds on non-vanishing Kluyver sums in terms of the level $N$. This is relevant for determining a good value of the constant $c_\Gamma$ (polynomial in $\log(N)$?) in Theorem \ref{mainIm} or a general group $\Gamma$ of level $N$ (via a lower bound on the implied constant in \eqref{uOmega}). 
\item Describe the exact limit geodesic configuration, or the limit of the zero set as the weight tends to infinity, for non-principal congruence groups, such as $\Gamma_0(N)$, in terms of $N$ alone. 
\item Improve the estimates of the cardinalities of sets of zeros with given convergence speed in Theorem \ref{main2}, or give the exact numbers; find all zeros converging at exponential speed. Discuss the exact convergence speed of zeros for non-principal congruence groups, or for even $N$. 
\item Discuss zeros of Eisenstein series of odd weight for congruence groups.
\item Find the Hausdorff limit of zeros of the `Hecke' Eisenstein series $G_k^{A,B}$ in $\F$ as the (even) weight increases. Experiments suggest that these zeros are all located on the boundary arc, or on specific vertical lines defined by fixed rational numbers depending on $N$. For example, for odd $N$, the zeros appear to be on $|z|=1$ and $\Re(z)=\pm 1/(2j)$ for $j=1,\dots,(N-1)/2$. 
\item Study the arithmetic of the polynomial whose $j$-invariants are the zeros of the norm modular forms $\mathcal N_{\Gamma,k}$. How does the discriminant factor, what is the Galois group, and does it satisfy interesting congruences modulo primes related to the weight? Compare with \cite{GC} for the case of full level in function fields. 
\end{enumerate}

\bibliographystyle{amsplain}
\providecommand{\bysame}{\leavevmode\hbox to3em{\hrulefill}\thinspace}
\providecommand{\MR}{\relax\ifhmode\unskip\space\fi MR }
\providecommand{\MRhref}[2]{%
  \href{http://www.ams.org/mathscinet-getitem?mr=#1}{#2}
}
\providecommand{\href}[2]{#2}

\end{document}